\documentclass[11pt,a4paper]{article}
\usepackage[latin1]{inputenc}
\usepackage{mathtools}
\usepackage{amsfonts}
\usepackage{amssymb,amsthm}
\usepackage{graphicx}
\usepackage[scale=0.775]{geometry}
\usepackage[dvipsnames]{xcolor}
\usepackage{multirow}
\usepackage{array}
\usepackage{booktabs}
\usepackage{hyperref}
\usepackage{pdflscape}
\usepackage{subcaption}

\newtheorem{theorem}{Theorem}[section]
\newtheorem{lemma}[theorem]{Lemma}

\title{An all Mach number finite volume method \\
    for isentropic two-phase flow}

\author{
	M\'aria Luk{\'a}{\v{c}}ov{\'a}-Medvid'ov{\'a}\footnote{Institut f\"ur Mathematik, Johannes-Gutenberg-Universit\"at Mainz, Staudingerweg 9, 55099 Mainz, Germany ({lukacova@mathematik.uni-mainz.de},{athomann@uni-mainz.de})},
	Gabriella Puppo\footnote{Dipartimento di Matematica, La Sapienza Universit\`a di Roma, Piazzale Aldo Moro 5, 00185 Roma, Italy ({gabriella.puppo@uniroma1.it})}, Andrea Thomann\footnotemark[1]~\footnote{Corresponding author}
	
}

\date{\today}
\begin{document}
	\maketitle
	
	% equation numbers
	\numberwithin{equation}{section}
	
	\section*{Abstract}
	We present an implicit-explicit finite volume scheme for isentropic two phase flow in all Mach number regimes. 
	The underlying model belongs to the class of symmetric hyperbolic thermodynamically compatible models.
	The key element of the scheme consists of a linearisation of pressure and enthalpy terms at a reference state.
	The resulting stiff linear parts are integrated implicitly, whereas the non-linear higher order and transport terms are treated explicitly. 
	Due to the flux splitting, the scheme is stable under a CFL condition which determined by the resolution of the slow material waves and allows large time steps even in the presence of fast acoustic waves.
	Further the singular Mach number limits of the model are studied and the asymptotic preserving property of the scheme is proven. 
	In numerical simulations the consistency with single phase flow, accuracy and the approximation of material waves in different Mach number regimes are assessed.  
	
	\paragraph{Keywords.} All-speed scheme, RS-IMEX, two-phase flow, asymptotic preserving, Symmetric Hyperbolic Thermodynamically Compatible models

	\section{Introduction}
	   
    Multi-phase flows are omnipresent in environmental and industrial processes.
    The broad range of applications poses an intrinsic problem of modeling two-phase flows.
    A widely used model was introduced by Baer \& Nunziato \cite{BaeNun1986} and still forms the basis of many models used to describe compressible two-phase flows.  
    Since then a wide range of modifications and extensions towards different applications have been proposed, see \cite{AmbChaRav2012,AndWar2004,KapMenBdzSonSte2001,KreKor2010,MueHanRich2016,SauAbg1999} and contributions mentioned therein.
    These models are based on conservation laws of mass, momentum and energy for each phase. 
    However the system cannot be written in a flux conservative form which causes problems in predicting correct shock speeds and the formulation of Rankine-Hugoniot conditions \cite{GouGav1999}.
    Therefore special techniques in the numerical treatment of non-conservative products are required as proposed, for example, in \cite{AndSauWar2003}.
    
    Here, we consider an alternative to the Baer \& Nunziato formulation of two-phase mixtures, namely a symmetric hyperbolic model in conservation form proposed in \cite{RomDriTor2010,RomTor2004}.
    It is based on the theory of Symmetric Hyperbolic Thermodynamically Compatible (SHTC) systems \cite{God1961,GodRom2003,GodRom1995}.
   	The latter class of equations is derived from thermodynamics \cite{GodRom2003,GodRom1995,Romenski1998} and variational principles \cite{PesPavRomGrm2018}. 
   	The approach is versatile and not
restricted to the modeling of two-phase flows. 
    It constitutes a monolithic mathematical framework
that encompasses the evolution of all considered materials and provides a unified mathematical
description of multi-physics systems, see e.g. \cite{RomBelPes2016} for a generalization of the two-phase flow model
given in \cite{RomDriTor2010,RomTor2004} to an arbitrary number of phases, \cite{DumPesRomZan2016,PesDumBosRomChiIor2021} for applications of the SHTC theory to fluid and solid mechanics modeling, or \cite{RomResPesDum2020} for recent advances in the description of
poroelastic fluid-saturated media.
    
    In this work, we focus on the isentropic setting of the two-phase flow model given in \cite{RomDriTor2010,RomTor2004}.
    Due to the conservative model formulation, the characteristic fields and wave relations were recently analysed in \cite{TheRomDum2022unpubl}.  
	Here, we are interested in the numerical simulation of gas-liquid interactions as they occur in air-water mixtures in form of droplets in air or dispersed bubbles in water. 
	Other possible applications include pipe flows where the transported medium can exist in its liquid and gas state due to depressurization events.
	Thereby the considered phases exhibit different behaviour with respect to flow properties ranging from compressible for gases to almost incompressible for some liquids.   
	Depending on the application this can imply a significant difference in the propagation speed of acoustic waves.
    Consequently, the Mach numbers that characterize the flow regime of each phase can differ in several orders of magnitude.  
      
    The construction of schemes that are designed for applications in low Mach number regimes is an active field of research, for models based on the Baer-Nunziato model see \cite{DemScaPelBra2022,ReAbg2021arXiv,Pelanti2017,DarLeQDulLeM2010}.
    In \cite{ZouAudTen2020} the low Mach limit is considered for the compressible gas phase only, where the liquid phase is described by incompressible flow.
    The novelty in our work lies in the fact that in the non-dimensionalisation of the here considered model two Mach numbers are considered. 
    They are given by the ratio between the local flow velocity of the mixture and the respective sound speeds of the phases.
	
	A severe difficulty in the construction of a numerical scheme applied to weakly compressible flow regimes is posed by the scale differences between acoustic waves and the material wave. 
	The focus of the numerical simulation usually lies on the evolution of the slower material waves for which a time step oriented towards the local flow speed suffices. 
	The time step of an explicit scheme, as proposed in \cite{RomDriTor2010,RomTor2004} for compressible two-phase flow, is bounded by the smallest appearing Mach number.
	This leads to vanishingly small time steps in low Mach number regimes and consequently to long computational times, especially when long time periods are considered. 
	This defect can be overcome by considering implicit-explicit (IMEX) time integrators, where fast waves are treated implicitly leading to the Courant-Friedrichs-Levy (CFL) condition that is restricted only be the local flow velocity.
    This allows larger time steps while keeping the material wave well resolved. 
	Additionally, an implicit treatment of the associated stiff pressure terms, which trigger fast acoustic waves, has the advantage that centered differences can be applied without loss of stability.
    This is important for obtaining the correct numerical viscosity in low Mach number flows as established in \cite{Dellacherie2010,GuiVio1999}, see also \cite{FeiDolKuc2007,FeiKuc2008} for a particular successive linearisation approach.
	In particular, the upwind schemes suffer from an excessive numerical diffusion \cite{GuiVio1999,Klein1995} and are therefore not applicable. 
	Indeed, the correct amount of numerical diffusion is an integral part to obtain so-called asymptotic preserving (AP) schemes \cite{Jin1999}.
	Since the flow regime of the two-phase flow considered here is characterized by two potentially distinct phase Mach numbers, different singular Mach number limits can be obtained which depend on the constitution of the mixture. 
	For their formal derivation we apply asymptotic expansions, as done for the (isentropic) Euler equations \cite{CorDegKum2012,DegTan2011,Dellacherie2010,GuiVio1999,KlaiMaj1981,Klein1995}, see also references therein.
	To obtain physically admissible solutions, especially in the weakly compressible
	flow regime, the numerical scheme has to preserve these asymptotics.
	This means a consistent discretization of the limit equations as the Mach numbers tend to zero. 
	To the best of our knowledge, this is the first systematic study of the effect of the Mach number scalings in two phase flows indicating the important terms in weakly compressible regimes.
	
	The profound knowledge of the structure of well-prepared initial data can be used to construct an AP scheme by applying a reference solution (RS)-IMEX approach. 
	This approach was successfully applied to construct AP schemes for the (isentropic) Euler equations \cite{BisLukYel2017,KaiSchSchNoe2017,KucLukNoeSch2021,ZeiSchKaiBecLukNoe2020}.
	Here, we only linearise the nonlinear pressure based terms around a reference state given by well-prepared data. 
	The stiff linear part is then treated implicitly whereas the non-linear higher order terms are integrated explicitly respecting the asymptotics in the low Mach number limit. 
	By doing this, nonlinear implicit solvers can be avoided which are computationally costly. 
	
	The paper is organized as follows. 
	In Section \ref{sec:Model} we revisit the model formulation from \cite{RomTor2004} and give the non-dimensional formulation for which afterwards well-prepared initial conditions and singular Mach number limits are derived. 
	Motivated by the structure of the well-prepared data, the numerical scheme is constructed in Section \ref{sec:Scheme} based on the reference solution approach. 
	The derivation of the time semi-discrete scheme is discussed in detail and the fully discrete RS-IMEX scheme based on a finite volume framework is formulated.
	The subsequent section is dedicated to the AP proof. 
	The scheme is validated by numerical tests in different Mach number regimes in Section \ref{sec:NumRes}. 
	In particular, the consistency with single phase flow, first order accuracy and the behaviour of the scheme for different Riemann problems are assessed. 
	Final conclusions are drawn in Section \ref{sec:Concl}.

	\section{Isentropic two-phase flow}
	\label{sec:Model}
	\subsection{The compressible model}
	The one-dimensional isentropic two phase model as introduced in \cite{RomDriTor2010} is given by 
	\begin{subequations}
	\label{eq:Comp_model}
    \begin{align}
		\partial_t \rho + \partial_x (\rho u) &= 0, \label{eq:Comp_model_mix_mass}\\
		\partial_t (\alpha_1 \rho) + \partial_x (\alpha_1 \rho u) &= -\frac{1}{\tau} \left(p_2 - p_1\right), \label{eq:Comp_model_vol_frac}\\
		\partial_t (\alpha_1 \rho_1) + \partial_x (\alpha_1 \rho_1 u_1) & = 0, \label{eq:Comp_model_phase_mass}\\
		\partial_t (\rho u) + \partial_x \big(\alpha(\rho_1 u_1^2 + p_1) + \alpha_2 (\rho_2 u_2^2 + p_2)\big) &= 0, \label{eq:Comp_model_mix_mom}\\
		\partial_t (u_1 - u_2) + \partial_x \left( \frac{u_1^2}{2} + h_1 - \frac{u_2^2}{2} - h_2\right) &= - \zeta\chi_1 \chi_2 (u_1 - u_2). \label{eq:Comp_model_rel_vel}
	\end{align}
	\end{subequations}
	System \eqref{eq:Comp_model} consists of the conservation of mixture mass $\rho$ \eqref{eq:Comp_model_mix_mass} and partial mass $\alpha\rho_1$ \eqref{eq:Comp_model_phase_mass}. 
	Equation \eqref{eq:Comp_model_vol_frac} gives the balance of the evolution of the volume fraction $\alpha$ with respect to a pressure relaxation source term, where $\tau$ denotes the relaxation rate. 
	Further, equation \eqref{eq:Comp_model_mix_mom} gives the conservation of mixture momentum $\rho u$ and equation \eqref{eq:Comp_model_rel_vel} the balance of the relative velocity $u_1 - u_2$ against a friction source term with a friction coefficient $\zeta$. 
	The densities of the respective phases are denoted by $\rho_1$ and $\rho_2$ and 
	the mixture density is given by 
	\begin{equation*}
		\rho = \alpha_1 \rho_1 + \alpha_2 \rho_2,
	\end{equation*}
	where the volume fraction $\alpha_1 \in (0,1)$ is associated to phase one and obeys $\alpha_1 + \alpha_2 = 1$. 
	The phase velocities are given by $u_1$ and $u_2$, respectively, and the mixture velocity is defined as 
	\begin{equation*}
		u = \chi_1 u_1 + \chi_2 u_2.
	\end{equation*}
	where $\chi_1 = \frac{\alpha_1\rho_1}{\rho}$ denotes the mass fraction of phase one and obeys the relation $\chi_1 + \chi_2 = 1$. 
	To close the system, we consider two different equations of states (EOS), the ideal gas law
	\begin{equation}
	\label{EOS:gas}
		p(\rho) = \kappa \left(\frac{\rho}{\rho_0}\right)^\gamma, \quad e(\rho) = \frac{p}{\rho (\gamma - 1)},
	\end{equation}
	and the stiffend gas equation to model liquids
	\begin{equation}
		\label{EOS:water}
		p(\rho) = \kappa \left(\frac{\rho}{\rho_0}\right)^\gamma - p_\infty, \quad e(\rho) = \frac{p + \gamma p_\infty}{\rho (\gamma - 1)}.
	\end{equation}
	The parameter $\gamma$ denotes the adiabatic constant and $p_\infty$ and $\kappa$ denote positive constants describing the considered medium. 
    The internal energy of the mixture is given by a linear combination of the internal energies $e_1,e_2$ of the respective phases
    \begin{equation}
    \label{eq:Mix_internal_E}
        e(\rho_1,\rho_2) = \chi_1 e_1 + \chi_2 e_2.
    \end{equation}
    The mixture pressure $P$ is obtained from \eqref{eq:Mix_internal_E} by
    \begin{equation}
    \label{eq:Mix_pressure}
        P = \rho^2 \frac{\partial e}{\partial \rho} = \alpha_1 p_1 + \alpha_2 p_2,
    \end{equation}
    which is a linear combination of the phase pressures $p_1,p_2$.
    The phase enthalpies $h_1,h_2$ in equation \eqref{eq:Comp_model_rel_vel} are defined as $h_1 = e_1 + \frac{p_1}{\rho_1}$, $h_2 = e_2 + \frac{p_2}{\rho_2}$ and are determined by the respective EOS \eqref{EOS:gas},\eqref{EOS:water}. 
    Further, we can define the mixture sound speed from \eqref{eq:Mix_pressure} by
    \begin{equation}
    \label{eq:Mix_sound}
        c^2 = \frac{\partial P}{\partial \rho} = \chi_1 c_1^2 + \chi_2 c_2^2,
    \end{equation}
    where the phase sound speeds $c_1, c_2$ are given by 
    \begin{equation*}
    c_1^2 = \frac{\partial p_1(\rho_1)}{\partial \rho_1}, \quad c_2^2 = \frac{\partial p_2(\rho_2)}{\partial \rho_2}.
    \end{equation*}
	We distinguish the following types of variables, state variables 
	\begin{equation}
	\label{eq:State_var}
		W = (\rho, \alpha_1 \rho, \alpha_1 \rho_1, \rho u, u_1 - u_2)^T,
	\end{equation}
    mixture variables 
    \begin{equation}
    \label{eq:Mix_var}
        Q = (\alpha_1, \rho, \chi_1, u, u_1 - u_2)^T
    \end{equation}
	and the phase variables 
	\begin{equation*}
		V = (\alpha_1, \rho_1, u_1, \rho_2, u_2)^T.
	\end{equation*}
	Using the state variables $W$, we can write \eqref{eq:Comp_model} in a compact form as
	\begin{equation*}
		\partial_t W + \partial_x f(W) = r(W),
	\end{equation*}
	where $f$ denotes the nonlinear flux function and $r$ the physical relaxation source terms. 
	Following \cite{RomTor2004}, model \eqref{eq:Comp_model} is strictly hyperbolic with the following eigenvalues
	\begin{equation*}
		\lambda^u = u, \quad \lambda_1^\pm = u_1 \pm c_1, \quad \lambda_2^\pm = u_2 \pm c_2.
	\end{equation*}
	They can be obtained by diagonalising the Jacobian $\partial_W f(W)$. 
	Thus, model \eqref{eq:Comp_model} exhibits two acoustic waves for each phase $\lambda_1^\pm, ~\lambda_2^\pm$ and the mixture velocity $\lambda^u$ which is also referred to as material velocity.
	Especially, when the sound speeds $c_1$ and/or $c_2$ are large, the acoustic waves $\lambda_1^\pm$ and/or $\lambda_2^\pm$ travel consistently faster than the material wave which introduces different scales into the model. 
	
	\subsection{Non-dimensional formulation}
	\label{sec:NonDim}
	
	To obtain a better understanding of the scales that are present in the model, we rewrite \eqref{eq:Comp_model} in non-dimensional form. 
	Let us denote the non-dimensional quantities by $\left(~\widetilde{\cdot}~\right)$ and the corresponding reference value by $(\cdot)^r$.
	We assume that the convective scales are of the same order, ie. $u_1^r = u_2^r = u^r = x^r/t^r$ which can be expressed through the reference length $x^r$ and time $t^r$.
	The ratio of the phase densities however can be large especially when considering mixtures of light gases and liquids. 
	To take this potentially large difference into account, we define two different reference densities which are connected to the reference mixture density by $\rho_1^r = \varrho_1 \rho^r, ~ \rho_2^r = \varrho_2 \rho^r$, where $\varrho_1, \varrho_2 \in \mathbb{R}$ are scaling constants.
	Further we define two different reference pressures $p_1^r$, $p_2^r$ from which we can compute the reference sound speeds given by
	\begin{equation*}
		\left(c_1^r\right)^2 = \frac{p_1^r}{\rho_1^r} = \frac{1}{\varrho_1}\frac{p_1^r}{\rho^r}, \quad \left(c_2^r\right)^2 = \frac{p_2^r}{\rho_2^r} = \frac{1}{\varrho_2}\frac{p_2^r}{\rho^r}.
	\end{equation*}
	For each phase, we can define a reference Mach number which is given by the ratio of the reference velocities and sound speeds
	\begin{equation*}
	M_1 = \frac{u^r}{c_1^r} = \sqrt{\varrho_1} ~\frac{u^r}{\sqrt{p_1^r/\rho^r}} = \sqrt{\varrho_1} ~M_1^\ast, \quad M_2 = \frac{u^r}{c_2^r} = \sqrt{\varrho_2}~\frac{u^r}{\sqrt{p_2^r/\rho^r}} = \sqrt{\varrho_2}~M_2^\ast.
	\end{equation*}
	We note that the Mach numbers $M_1, M_2$ are defined using reference phase densities, whereas $M_1^\ast, M_2^\ast$ are defined using the reference mixture density only.  
	Summarizing, we can express the dimensional variables as the product of non-dimensional quantity and reference value as follows
	\begin{align}
		\rho &= \left(\alpha_1 \varrho_1 \widetilde\rho_1 + \alpha_2 \varrho_2 \widetilde\rho_2\right) \rho^r, &\quad \rho_1 &= \widetilde{\rho_1} \varrho_1 ~\rho^r, &\quad \rho_2 &= \widetilde{\rho_2} \varrho_2 ~\rho^r, \notag\\
		u &= \left(\frac{\alpha_1 \varrho_1 \widetilde\rho_1}{\widetilde \rho}\tilde u_1 + \frac{\alpha_2 \varrho_2 \widetilde\rho_2}{\widetilde{\rho}}\tilde u_2\right) u^r, &\quad u_1 &= \widetilde{u}_1 ~u^r, &\quad u_2 &= \widetilde{u}_2 ~u^r, 	\label{eq:Decomp_var}\\
		&&p_1 &= \widetilde{p_1} ~p_1^r, &\quad p_2 &= \widetilde{p_2} ~p_2^r. \notag
	\end{align}
	In addition, we have also the following reference values for the pressure relaxation rate and friction coefficient given respectively by
	\begin{equation}
	\label{eq:Ref_tau_fric}
		\tau^r = (u^r)^2 t^r \quad \text{ and } \quad \zeta^r = \frac{1}{t^r} \quad \text{ defining } \quad \tau = \widetilde{\tau} ~\tau^r, \quad \zeta = \widetilde{\zeta} ~\zeta^r.
	\end{equation}
	Inserting expressions \eqref{eq:Decomp_var}, \eqref{eq:Ref_tau_fric} into \eqref{eq:Comp_model} and dropping the $(\,\widetilde{\cdot} \,)$, we obtain the following non-dimensional formulation
	\begin{subequations}
		\label{eq:NonDim_model}
	\begin{align}
	\partial_t \rho + \partial_x (\rho u) &= 0, \label{eq:ND_dens} \\
	\partial_t (\alpha_1 \rho) + \partial_x (\alpha_1 \rho u) &= -\frac{1}{\tau} \left(\frac{\varrho_2 p_2}{M_2^2} - \frac{\varrho_1 p_1}{M_1^2}\right), \label{eq:ND_alphadens}\\
	\partial_t (\alpha_1 \rho_1) + \partial_x (\alpha_1 \rho_1 u_1) & = 0, \label{eq:ND_alphadens1}\\
	\partial_t (\rho u) + \partial_x \left(\alpha_1 \varrho_1\rho_1 u_1^2 + \frac{\alpha_1 \varrho_1 p_1}{M_1^2} + \alpha_2 \varrho_2\rho_2 u_2^2 + \frac{\alpha_2 \varrho_2 p_2}{M_2^2}\right) &= 0, \label{eq:ND_mom}\\
	\partial_t (u_1 - u_2) + \partial_x \left( \frac{u_1^2}{2} - \frac{u_2^2}{2} + \frac{h_1}{M_1^2} - \frac{h_2}{M_2^2}\right) &= - \zeta \chi_1 \chi_2 (u_1 - u_2) \label{eq:ND_relvel}
	\end{align}
	\end{subequations}
	with 
	\begin{equation*}
	    \rho = \alpha_1 \varrho_1 \rho_1 + \alpha_2 \varrho_2 \rho_2, \quad u = \chi_1 u_1 + \chi_2 u_2, \quad 
	\chi_1 = \varrho_1 \frac{\alpha_1 \rho_1}{\rho}.
	\end{equation*}
	Note that the mass fraction still obeys $\chi_1 + \chi_2 = 1$.
	Analogously to the dimensional formulation \eqref{eq:Comp_model}, the non-dimensional model \eqref{eq:NonDim_model} is strictly hyperbolic and exhibits 5 waves given by
	\begin{equation*}
	\lambda^u = u, \quad \lambda_1^\pm = u_1 \pm \frac{c_1}{M_1}, \quad \lambda_2^\pm = u_2 \pm \frac{c_2}{M_2}.
	\end{equation*} 
	We see that the the acoustic waves scale with the phase Mach numbers $M_1,M_2$ respectively and propagate significantly faster than the material wave for low Mach numbers. 
	To obtain more insight on the behaviour of the solution in the low Mach number limit, we perform an asymptotic analysis of the non-dimensional model \eqref{eq:NonDim_model}. 
	This is subject of the next section.
	
	\subsection{Well-prepared data and the low Mach limit}
	\label{sec:WP}
	According to \cite{Dellacherie2010} an initial condition of \eqref{eq:NonDim_model} is called well-prepared, if it is close to a solution of a limit model for small Mach numbers. 
	For the derivation of well-prepared initial data and the low Mach limit model of \eqref{eq:NonDim_model}, we focus on the scales induced by the Mach numbers $M_1, M_2$.
	The scaling factors of the densities $\varrho_1, \varrho_2$ are considered to be fixed values independent of the Mach regimes. 
    For simplicity of notation, we will neglect the scaling parameters $\varrho_1, \varrho_2$ in the following analysis. 
%    \cite{CorDegKum2012,Dellacherie2010} 
	We consider the subsequent cases in detail: 
	\begin{itemize} 
        \item \textbf{Case 1:} Phase one is compressible, i.e. $M_1 = 1$, and phase two is characterized by a low Mach number $M_2 \ll 1$. 
	    \item \textbf{Case 2:} Both phases are in the same Mach regime, i.e. $M_1 = M_2 = M$.
%	    \item The second case consists of two phases being in different Mach regimes  $M_1, M_2 < 1$, where without loss of generality the first phase is described by $M_1 = M$ and the second can be expressed by terms of $M$ by $M_2 = M^\delta$ with $\delta > 1$, {\color{blue} where $\delta \in \mathbb{N}.$} {\color{red} with $\delta \in \mathbb{Q}$?, if in $\mathbb{N}$ do we need to use $M_2 = \mathcal{C} M^\delta$?}. 
	\end{itemize}
	\paragraph{Case 1}
	Let phase one be compressible with $\rho_1 = \mathcal{O}(1)$ and $u_1 = \mathcal{O}(1)$. 
    The variables of the weakly incompressible phase two are expanded with respect to a Mach number $M$ in the following way
    \begin{equation}
    \label{eq:R2_U2_exp}
        \rho_2 = \rho_2^{(0)} + M \rho_2^{(1)} + M^2 \rho_2^{(2)} + \mathcal{O}(M^3), \quad u_2 = u_2^{(0)} + M u_2^{(1)} + \mathcal{O}(M^2).
    \end{equation}
    The Mach number expansion of the pressure can be obtained from the density expansion \eqref{eq:R2_U2_exp} via the respective EOS. 
    With $p_2(\rho_2^{(0)}) = p_2^{(0)}$ and $(c_2^{(0)})^2 = \gamma_2 \frac{p_2(0)}{\rho_2^{(0)}}$ we obtain the following expansion
    \begin{align}
    %		p_1(x,t) = p_1^{(0)} + \frac{p_1^{(0)} \gamma_1 \rho_1^{(1)}}{\rho_1^{(0)}} M + \frac{p_1^{(0)} \gamma_1 \left(2 \rho_1^{(0)} \rho_1^{(2)} + (\gamma_1 -1)\left(\rho_1^{(1)}\right)^2\right)}{\left(\rho_1^{(0)}\right)^2} M^2 + \mathcal{O}(M^3)
    p_2(x,t) = p_2^{(0)} + (c_2^{(0)})^2\rho_2^{(1)} M + \left(\frac{1}{2}(1-\gamma_2)\frac{(c^{(0)}_2)^2}{\rho_2^{(0)}} (\rho_2^{(1)})^2 + (c^{(0)}_2)^2\rho_2^{(2)}\right) M^2 + \mathcal{O}(M^3). \label{eq:Exp_p}
    \end{align}
    Further we assume a Mach number expansion of the volume fraction given by 
    \begin{equation}
    \label{eq:M_exp_alpha}
        \alpha = \alpha^{(0)} + \mathcal{O}(M).
    \end{equation}
    Inserting the expansions into the non-dimensional equations \eqref{eq:NonDim_model} and sorting by orders of the Mach number we find from \eqref{eq:ND_alphadens}, \eqref{eq:ND_mom} and \eqref{eq:ND_relvel} for the $\mathcal{O}(M^{-2})$ order terms
    \begin{equation*}
        p_2^{(0)} = 0, \quad 
        \partial_x \left(\alpha^{(0)}p_2^{(0)}\right) = 0, \quad
        \frac{\partial_x p_2^{(0)}}{\rho_2^{(0)}} = 0. 
    \end{equation*}
    We immediately find from \eqref{eq:Exp_p} and the EOS that $\rho_2^{(0)}$ is a non-negative constant. 
    Especially for an ideal gas follows $\rho_2^{(0)} = 0$ which means phase two is vanishing or in vacuum at leading order. 
    Since we are interested in obtaining a mixture of two phases also in the limit, we assume phase two to be associated with the stiffened gas equation, thus $\rho_2^{(0)}$ is positive and constant.
    Analogously, we find for the order $\mathcal{O}(M^{-1})$ terms
    \begin{equation*}
    p_2^{(1)} = (c_2^{(0)})^2\rho_2^{(1)} = 0
    \end{equation*}
    which implies $\rho_2^{(1)} = 0$ according to \eqref{eq:Exp_p}. 
    Summarizing, the well-prepared data for phase density two is given by
    \begin{equation*}
        \rho_2 = \rho_2^{(0)} + \mathcal{O}(M^2), \quad \rho_2^{(0)} \text{ const. }
    \end{equation*}
    For $\mathcal{O}(1)$ terms we find for the mixture density and velocity the following relations
    \begin{equation*}
        \rho^{(0)} = \alpha^{(0)}\rho_1 + (1-\alpha^{(0)})\rho_2^{(0)}, \quad u^{(0)} = \frac{\alpha^{(0)}\rho_1}{\rho^{(0)}} u_1 + \frac{(1-\alpha^{(0)}\rho_2^{(0)})}{\rho^{(0)}}u_2^{(0)} = \chi_1^{(0)} u_1 + \chi_2^{(0)} u_2^{(0)}.
    \end{equation*}
    Using this notation, we find
    \begin{align}
        \partial_t \alpha_2^{(0)} + \partial_x (\alpha_2^{(0)} u_2^{(0)}) & = 0 \label{eq1}\\
        \partial_t \alpha^{(0)} + u^{(0)}\partial_x \alpha^{(0)} &= -\frac{1}{\tau \rho^{(0)}} (p_2^{(2)} - p_1) \label{eq2}\\
        \partial_t (\alpha^{(0)}\rho_1) + \partial_x (\alpha^{(0)}\rho_1 u_1) &= 0 \label{eq:3} \\
        \partial_t (\rho^{(0)} u^{(0)}) + \partial_x (\alpha^{(0)} \rho_1 u_1^2 + \alpha^{(0)}p_1 + \alpha_2^{(0)}\rho_2^{(0)}(u_2^{(0)})^2 + \alpha^{(0)}_2 p_2^{(2)} &= 0 \label{eq:4}\\
        \partial_t (u_1 - u_2^{(0)}) + \partial_x \left(\frac{1}{2}u_1^2 - \frac{1}{2}(u_2^{(0)})^2\right) + \frac{\partial_x p_1}{\rho_1} - \frac{\partial_x p_2^{(2)}}{\rho_2^{(0)}} &= \zeta \chi_1^{(0)} \chi_2^{(0)} (u_1 - u_2^{(0)}). \label{eq:5}
    \end{align}
    Integrating \eqref{eq1} over the domain $\Omega$ with periodic or no-flux boundary conditions, as done in \cite{Dellacherie2010,DegTan2011}, we obtain that the average of $\alpha^{(0)}$ is constant in time. 
    From this it follows immediately that 
    \begin{equation}
    \label{eq:Limit1_alpha}
        u^{(0)} \partial_x \alpha^{(0)} = - \frac{1}{\rho^{(0)}}(p_2^{(2)} - p_1)
    \end{equation}
    and we obtain
    \begin{equation*}
        \partial_x u_2^{(0)} = \frac{u_2^{(0)}}{\alpha_2^{(0)}}\partial_x \alpha^{(0)}.
    \end{equation*}
    Further we get from \eqref{eq:4} and \eqref{eq:5}
    \begin{equation}
    \label{eq:Limit1_u2}
        \partial_t u_2^{(0)} + u_2^{(0)}\partial_x u_2^{(0)} + \frac{\partial_x p_2^{(2)}}{\rho_2^{(0)}} + \frac{p_2^{(2)} - p_1}{\rho^{(0)}}\partial_x \alpha^{(0)} = \zeta (\chi_1^{(0)})^2\chi_2^{(0)} (u_1 - u_2^{(0)}).
    \end{equation}
    For the compressible phase we obtain from \eqref{eq:3}, \eqref{eq:4} and  \eqref{eq:Limit1_u2}
    \begin{align}
   	\label{eq:Limit1_r1_u1}
    \begin{split}
        \partial_t (\alpha^{(0)} \rho_1) + \partial_x (\alpha^{(0)}\rho_1 u_1) &= 0 \\
        \partial_t u_1 + u_1 \partial_x u_1 + \frac{\partial_x p_1}{\rho_1} + \frac{p_2^{(2)} - p_1}{\rho}\partial_x \alpha^{(0)} &= -\zeta \chi_1^{(0)}(\chi_2^{(0)})^2 (u_1 - u_2^{(0)}).
    \end{split}
	\end{align}
    Note that $\partial_x \alpha^{(0)} = 0$, i.e. constant $\alpha^{(0)}$, yields $\partial_x u_2^{(0)} = 0$ and in the limit $M \to 0$ we formally obtain the incompressible Euler equations with friction for phase two and the compressible isentropic Euler equations with friction for phase one. 
    In particular, in the limit we have only one pressure given by $p_1$ and both phases are decoupled. 
    
    Summarizing, we have derived the following result regarding well-prepared initial data for two phase flow consisting of a compressible and a weakly incompressible phase.
    \begin{lemma}[Well-prepared data for compressible/weakly compressible flow]
    	Let phase one be compressible, i.e. characterized by $M_1 = 1$ and phase two be weakly compressible, i.e. characterized by $M_2 = M \ll 1$.
    	Let $\alpha, \rho_2, u_2$ be given by the Mach number expansions \eqref{eq:M_exp_alpha}, \eqref{eq:R2_U2_exp} and the set of well-prepared initial data be defined as 
    	\begin{equation}
    	\label{eq:wp_1}
    	\Omega_1^{wp} = \left\lbrace W \in \mathbb{R}^{2d + 3}: \rho_2^{(0)} \text{ const. },~ p_2(\rho_2^{(0)}) = 0, ~\rho_1^{(1)} = 0, ~\partial_t \alpha^{(0)} = 0, ~\partial_x u_2^{(0)} = \frac{u_2^{(0)}}{\alpha_2^{(0)}} \partial_x \alpha^{(0)} \right\rbrace.
    	\end{equation}
    	Then formally for $W \in \Omega_1^{wp}$ for $M \to 0$ the limit equations are given by \eqref{eq:Limit1_alpha},\eqref{eq:Limit1_u2},\eqref{eq:Limit1_r1_u1}.
    	If in addition holds $\partial_x \alpha^{(0)} = 0$, i.e. $\alpha^{(0)}$ is constant, then the limit equations are given by 
    	\begin{align}
    		\label{eq:Limit1_cons}
    	\begin{split}
    	\partial_t u_2^{(0)} + \frac{\partial_x p_1}{\rho_2^{(0)}} &= \zeta (\chi_1^{(0)})^2 \chi_2^{(0)} (u_1 - u_2^{(0)}), \quad \partial_x u_2^{(0)} = 0\\
    	\partial_t \rho_1 + \partial_x (\rho_1 u_1) & = 0, \\
    	\partial_t (\rho_1 u_1) + \partial_x (\rho_1 (u_1)^2 + \partial_x p_1) &= - \zeta \chi_1^{(0)}(\chi_2^{(0)})^2 (u_1 - u_2^{(0)}).
    	\end{split}
    	\end{align}
    \end{lemma}
	\paragraph{Case 2}
    When the flow of both phases can be characterized by the same Mach number $M$, i.e. $M_1 = \mathcal{O}(M)$ and $M_2 = \mathcal{O}(M)$, we consider a Mach number expansion of the variables of both phases. 
	Denoting the different phases by $k=1,2$, we obtain for the densities and velocities
	\begin{align}
		\rho_k(x,t) &= \rho_k^{(0)}(x,t) + \rho_k^{(1)}(x,t) M + \rho_k^{(2)}(x,t) M^2 + \mathcal{O}(M^3), \label{eq:Exp_rho} \\
		u_k(x,t) &= u_k^{(0)}(x,t) + u_k^{(1)}(x,t) M + \mathcal{O}(M^2). \label{eq:Exp_u}
	\end{align}
    As in the previous case, we can write a Mach number expansion of the associated phase pressures given by \eqref{eq:Exp_p} for each phase respectively and we consider a Mach number expansion of the volume fraction given by \eqref{eq:M_exp_alpha}.

%	Further we need the Mach number expansion of $p/\rho$ for example in the enthalpy which is given by 
%	\begin{equation}
%		\frac{p_1}{\rho_1}(x,t) = \frac{p_1^{(0)}}{\rho_1^{(0)}} + \frac{p_1^{(1)} \rho_1^{(0)} - p_1^{(0)} \rho_1^{(1)}}{\left(\rho_1^{(0)}\right)^2} M + \frac{\rho_1^{(0)}\left(p_1^{(2)} \rho_1^{(0)} - p_1^{(1)} \rho_1^{(1)} \right)+ p_1^{(0)} \left((\rho_1^{(1)})^2 - \rho_1^{(0)} \rho_1^{(2)}\right) }{\left(\rho_1^{(0)}\right)^3}M^2 + \mathcal{O}(M^3)
%	\end{equation}
%	The expansions of the mixture density and velocity follows directly from \eqref{eq:Exp_rho} and \eqref{eq:Exp_u} and is given by
%	\begin{equation}
%		\rho = \rho^{(0)} + M \rho^{(1)} + M^2 \rho^{(2)} + \mathcal{O}(M^3), \quad u = u^{(0)} + M u^{(1)} + M^2 u^{(2)} + \mathcal{O}(M^3)
%	\end{equation} 
%	where the zero order terms are, respectively, given by 
%	\begin{equation}\rho^{(0)} = \alpha_1^{(0)} \varrho_1 \rho_1^{(0)} + \alpha_2^{(0)} \varrho_2 \rho_2^{(0)}, \qquad u^{(0)} = \varrho_1\frac{\alpha_1^{(0)}\rho_1^{(0)}}{\rho^{(0)}}u_1^{(0)} + \varrho_2 \frac{\alpha_2^{(0)} \rho_2^{(0)}}{\rho^{(0)}} u_2^{(0)} .\label{eq:Exp_rho_u}
%	\end{equation}
	Inserting the expansions \eqref{eq:Exp_p} into the balance law for the volume fraction \eqref{eq:ND_alphadens} and sorting the terms by orders of the Mach number, we find for  $\mathcal{O}(M^{-2})$ the following relation
	\begin{equation}
    \label{eq:MachExpRelaxP-2}
p_2^{(0)} = p_1^{(0)}. 
	\end{equation}
	Using \eqref{eq:MachExpRelaxP-2} in the momentum equation \eqref{eq:ND_mom}, we find 
	\begin{equation*}
    0 = \partial_x (\alpha_1^{(0)} p_1^{(0)} + \alpha_2^{(0)} p_2^{(0)}) = \partial_x p_1^{(0)}.
	\end{equation*}
	Analogously, we obtain $p_2^{(1)} = p_1^{(1)}$ and $\partial_x p_1^{(1)} = 0$ for the $\mathcal{O}(M^{-1})$ terms. 
    Summarizing, the pressure expansions are given by 
	\begin{equation}
		p_1(x,t) = p^{(0)}(t) + M^2 p_1^{(2)}(x,t) + \mathcal{O}(M^3), \quad 
	    p_2(x,t) = p^{(0)}(t) + M^2 p_2^{(2)}(x,t) + \mathcal{O}(M^3) \label{eq:Exp_p_final}
    \end{equation}	 
    with spatially constant component $p^{(0)}$.
    Comparing \eqref{eq:Exp_p_final} with \eqref{eq:Exp_p}, we obtain the following expansions for the phase densities
    \begin{equation*}
        \rho_1(x,t) = \rho_1^{(0)}(t) + M^2 \rho_1^{(2)}(x,t) + \mathcal{O}(M^3), \quad         \rho_2(x,t) = \rho_2^{(0)}(t) + M^2 \rho_2^{(2)}(x,t)  + \mathcal{O}(M^3),
    \end{equation*}
    where the zero order components only depend on time. 
    Following \cite{DegTan2011,Dellacherie2010}, we obtain by integrating \eqref{eq:ND_dens} and \eqref{eq:ND_alphadens1} on a domain $\Omega$ with periodic or no-flux boundary conditions the following conditions on the densities
    \begin{equation}
    	\partial_t \rho^{(0)} = 0, \quad \partial_t \left(\alpha_1^{(0)} \rho_1^{(0)}\right) = 0. \label{eq:Limit_density_same}
    \end{equation} 
    From \eqref{eq:Limit_density_same} we derive, using $\partial_x \rho_1^{(0)} = 0$ and $\partial_x \rho_2^{(0)} = 0$, the following conditions on the velocity derivatives
    \begin{equation}
    \label{eq:Limit_veloc_deriv}
        \partial_x u_1^{(0)} = - \frac{u_1^{(0)}}{\alpha_1^{(0)}} \partial_x \alpha^{(0)}, \quad \partial_x u_2^{(0)} = \frac{u_2^{(0)}}{\alpha_2^{(0)}} \partial_x \alpha^{(0)}.
    \end{equation}
    Looking at the $\mathcal{O}(1)$ terms we obtain the following limit system when $M$ goes to zero
    \begin{subequations}
    	\label{eq:Limit_M_same}
    \begin{align}
        \partial_t \alpha^{(0)} + u^{(0)} \partial_x \alpha^{(0)} & = -\frac{1}{\tau \rho^{(0)}} (p_2^{(2)} - p_1^{(2)}) \label{eq:Limit_M_same_alpha}\\
%        \partial_t \rho_1^{(0)} - \frac{\rho_1^{(0)}}{\alpha^{(0)}}\partial_x \alpha^{(0)} &= \frac{\rho_1^{(0)}}{\rho^{(0)}\alpha^{(0)}} \frac{1}{\tau} (p_2^{(2)} - p_1^{(2)}) \\
%        \partial_t \rho_2^{(0)} + \frac{\rho_2^{(0)}}{\alpha_2^{(0)}}\partial_x \alpha^{(0)} &= \frac{\rho_2^{(0)}}{\rho^{(0)}\alpha_2^{(0)}} \frac{1}{\tau} (p_2^{(2)} - p_1^{(2)}) \\
        \partial_t u_1^{(0)} + u_1^{(0)} \partial_x u_1^{(0)} + \frac{\partial_x p_1^{(2)}}{\rho_1^{(0)}} + \frac{p_2^{(2)} - p_1^{(2)}}{\rho^{(0)}}\partial_x \alpha^{(0)} &= -\zeta \chi_1^{(0)}(\chi_2^{(0)})^2 (u_1^{(0)} - u_2^{(0)}) \label{eq:Limit_M_same_u1}\\
         \partial_t u_2^{(0)} + u_2^{(0)}\partial_x u_2^{(0)} + \frac{\partial_x p_2^{(2)}}{\rho_2^{(0)}} + \frac{p_2^{(2)} - p_1^{(2)}}{\rho^{(0)}}\partial_x \alpha^{(0)} &= \zeta (\chi_1^{(0)})^2\chi_2^{(0)} (u_1^{(0)} - u_2^{(0)}) \label{eq:Limit_M_same_u2}
    \end{align}
\end{subequations}
    where $\partial_t \rho^{(0)} = 0$. 
    Analogously to Case 1, we see from \eqref{eq:Limit_veloc_deriv} that for $\partial_x \alpha^{(0)} = 0$ the derivatives of the velocities vanish and the phases are only coupled via the pressure relaxation and friction source terms.
	Moreover, equations \eqref{eq:Limit_M_same_u1} and \eqref{eq:Limit_M_same_u2} are two coupled incompressible Euler equations with variable densities. 
    For constant $\alpha^{(0)}$ however, we immediately find from \eqref{eq:Limit_density_same} that $\rho_1^{(0)}, \rho_2^{(0)}$ are constant and for the pressures holds $p_2^{(2)} = p_1^{(2)}$ in \eqref{eq:Limit_M_same_alpha}. 
    In particular, this leads to a system of incompressible Euler equations coupled via the friction source term with a single pressure $p^{(2)}$.

    Summarizing, we have derived the following result regarding well-prepared initial data for two weakly compressible phases in the same Mach regime. 
    
    \begin{lemma}[Well-prepared data for weakly compressible flow]
    Let both phases be weakly compressible in the same Mach number regime, i.e. characterized by the same Mach number $M$. 
    Let the phase variables $V \in \mathbb{R}^{2d+3}$ be given in the Mach number expansions \eqref{eq:Exp_rho}, \eqref{eq:Exp_u} and \eqref{eq:M_exp_alpha} and the set of well-prepared data be defined as 
    \begin{multline}
    \label{eq:wp_M}
    \Omega_M^{wp} = \left\lbrace w \in \mathbb{R}^{2d+3}: \partial_t \rho^{(0)} = 0,~ \partial_t (\alpha^{(0)}\rho_1^{(0)}) = 0, ~\partial_x \rho_1^{(0)} = 0, ~\partial_x \rho_2^{(0)} = 0, ~\rho_1^{(1)} = 0, ~\rho_2^{(1)} = 0,    \right. \\
    \left. p_1^{(0)} = p_2^{(0)}, ~\partial_x u_1^{(0)} = -\frac{u_1^{(0)}}{\alpha^{(0)}} \partial_x \alpha_1^{(0)}, ~\partial_x u_2^{(0)} = \frac{u_2^{(0)}}{\alpha_2^{(0)}} \partial_x\alpha^{(0)} \right\rbrace.
    \end{multline}
    Then formally for $W \in \Omega_M^{wp}$ for $M \to 0$ the limit equations are given by \eqref{eq:Limit_M_same}. 
    If in addition $\alpha^{(0)}$ is constant, we obtain constant $\rho_1^{(0)}, \rho_2^{(0)}$ and the following limit equations
    \begin{align}
    	\label{eq:Limit2_const_dens}
    \begin{split}
    	\partial_t u_1^{(0)} + \frac{\partial_x p^{(2)}}{\rho_1^{(0)}} &= -\zeta \chi_1^{(0)} (\chi_2^{(0)})^2 (u_1^{(0)} - u_2^{(0)}), \quad \partial_x u_1^{(0)} = 0\\
    	\partial_t u_2^{(0)} + \frac{\partial_x p^{(2)}}{\rho_2^{(0)}} &= \zeta (\chi_1^{(0)})^2 \chi_2^{(0)} (u_1^{(0)} - u_2^{(0)}), \quad \partial_x u_2^{(0)} = 0
    \end{split}
    \end{align}
    with a single pressure $p^{(2)} = p_1^{(2)} = p_2^{(2)}$. 
    \end{lemma}
	For completeness, we shortly mention the case where the two phases are weakly compressible and in different Mach number regimes. 
	Without loss of generality, we assume $M_2 \ll M_1 \ll 1$. 
	To obtain the simultaneous limit for $M_2 \to 0$ and $M_1 \to 0$ we consider the following Mach number expansions in $M_1, M_2$ of the phase variables
	\begin{equation*}
		V = \sum_{j,l=0}^{\infty} M_1^j M_2^l V^{(j,l)}.
	\end{equation*}
	Via the EOS we obtain an analogous Mach number expansion of the pressures. 
	Following the above procedure for Cases 1 and 2, we obtain the following set of well-prepared data given by
	\begin{multline}
	\label{eq:wp_2}
		\Omega_2^{wp} = \left\lbrace W \in \mathbb{R}^{2d+3}: \alpha^{(0,0)} \text{ const. }, ~\rho_1^{(0,0)} \text{ const. }, ~ \right.\\
		\left.\rho_2^{(0,0)} = \text{ const. }, \rho_1^{(j,l)} = 0, ~\rho_2^{(j,l)} = 0 \text{ for } j+l = 1, j=l=1, \right.\\
		\left. ~ p_1(\rho_1^{(0,0)}) = 0, ~p_2(\rho_2^{(0,0)}) = 0, ~p_1^{(2,0)} = p_2^{(0,2)},~\partial_x u_1^{(0,0)} = 0, ~ \partial_x u_2^{(0,0)} = 0. \right\rbrace
	\end{multline}
	In particular we obtain 
	\begin{equation*}
		\rho_1 = \rho_1^{(0,0)} + \mathcal{O}(M_1^2), \quad \rho_2 = \rho_2^{(0,0)} + \mathcal{O}(M_2^2).
	\end{equation*}
    The limit equations with constant $\alpha^{(0,0)}$ are given by 
        \begin{align*}
        \partial_t u_1^{(0,0)} + \frac{\partial_x p^{(2)}}{\rho_1^{(0,0)}} &= -\zeta \chi_1^{(0,0)} (\chi_2^{(0,0)})^2 (u_1^{(0,0)} - u_2^{(0,0)}), \quad \partial_x u_1^{(0,0)} = 0\\
        \partial_t u_2^{(0,0)} + \frac{\partial_x p^{(2)}}{\rho_2^{(0,0)}} &= \zeta (\chi_1^{(0,0)})^2 \chi_2^{(0,0)} (u_1^{(0,0)} - u_2^{(0,0)}), \quad \partial_x u_2^{(0,0)} = 0
        \end{align*}
    with a single pressure $p^{(2)} = p_1^{(2,0)} = p_2^{(0,2)}$. 

	\section{The numerical scheme}
	\label{sec:Scheme}
	The main objective for constructing a numerical scheme for the two phase model \eqref{eq:NonDim_model} is to achieve stability under a time step restriction independently of the Mach numbers $M_1$ and $M_2$.
	This allows to follow the material wave $\lambda^u$ while neglecting the resolution of the acoustic waves. 
	Therefore we desire a CFL condition of the type
	\begin{equation}
	\label{eq:MatCFL}
		\Delta t \leq \nu_{u} \frac{\Delta x}{\max|\lambda_u|}.
	\end{equation}
%	In the following we assume without loss of generality, that $1 \geq M_1 \geq M_2$. 
	An explicit scheme requires for stability a quite severe time step restriction of 
	\begin{equation*}
		\Delta t \leq \nu_{ac}\frac{\Delta x}{\max(|\lambda_1^\pm|,|\lambda_2^\pm|)} \leq \nu_{ac} \min(M_1,M_2) \frac{\Delta x}{\max (|M_2 u_2 \pm c_2|, |M_1 u_1 \pm c_1|)}
	\end{equation*}
	which vanishes as one of the Mach numbers tends to 0.
	To avoid the costs that arise from being forced to use small time steps in low Mach number regimes, we construct a numerical scheme based on an implicit-explicit (IMEX) approach where the fast waves are integrated implicitly thus do not contribute to the CFL condition.
	
	\subsection{ Reference solution approach}
	To determine which terms should be treated implicitly, we analyse the eigenstructure of the model in non-dimensional formulation \eqref{eq:NonDim_model}. 
	Thereby we find that the fast acoustic components of the eigenvalues $\lambda_1^\pm, \lambda_2^\pm$ stem from the respective phase pressure and enthalpy terms. 
	Therefore it is necessary to treat those terms implicitly in order to obtain a CFL condition that is independent of the Mach number regimes. 
	Considering the EOS given in \eqref{EOS:gas} and \eqref{EOS:water}, both pressures $p_1, p_2$ and enthalpies $h_1, h_2$ are non-linear functions of the densities $\rho_1, \rho_2$ which would require a nonlinear solver. 
	To avoid this, we linearise the phase pressures and enthalpies with respect to a known reference solution $\rho_1^{RS}$ and $\rho_2^{RS}$, respectively, as proposed in \cite{KaiSchSchNoe2017}. 
	Those reference states are motivated by the leading order terms of the Mach number expansions obtained in the derivation of well-prepared initial data given in \eqref{eq:wp_1}, \eqref{eq:wp_M} or \eqref{eq:wp_2}. 
	Here we focus on the case where the reference states are constant throughout the simulation thus considering the transition from weakly compressible to incompressible flow in the limit.
%	In case of time dependent reference states, for example when $\partial_t \rho^{(0)} \neq 0$, the reference states can be updated after each iteration using strategies discussed in \cite{}. 
	
	In the following, we detail the computations for the first phase only. 
	The formulations for the second phase are similar. 
	To linearise pressure and enthalpy, we consider the Taylor expansion with respect to the reference state $\rho^{RS}_1$ which reads
	\begin{align*}
		p_1(\rho_1) = p_1^{RS} + \left(c_1^{RS}\right)^2 (\rho_1 - \rho_1^{RS}) + \mathcal{O} \left(\left(\rho_1 - \rho_1^{RS}\right)^2\right), \\
		h_1(\rho_1) = h_1^{RS} + \frac{\left(c_1^{RS}\right)^2}{\rho_1^{RS}} \left(\rho_1 - \rho_1^{RS}\right) + \mathcal{O} \left(\left(\rho_1 - \rho_1^{RS}\right)^2\right),
	\end{align*}
	We split $p_1$ and $h_1$ into a part linear in $\rho_1$ 
	\begin{align}
	\hat p_1(\rho_1) = p_1^{RS} + \left(c_1^{RS}\right)^2 (\rho_1 - \rho_1^{RS}), \label{eq:Def_p_hat}\\
	\hat h_1(\rho_1) = h_1^{RS} + \frac{\left(c_1^{RS}\right)^2}{\rho_1^{RS}} \left(\rho_1 - \rho_1^{RS}\right) \label{eq:Def_h_hat}
	\end{align}
	and non-linear higher order terms
	\begin{align}
	\bar p_1(\rho_1) = p_1(\rho_1) - \hat p_1(\rho_1) = \mathcal{O} \left(\left(\rho_1 - \rho_1^{RS}\right)^2\right), \label{eq:p_bar}\\
	\bar h_1(\rho_1) = h_1(\rho_1) - \hat h_1(\rho_1) = \mathcal{O} \left(\left(\rho_1 - \rho_1^{RS}\right)^2\right). \label{eq:h_bar}
	\end{align} 
	Especially, if $\rho_1$ is well-prepared, ie. $\rho_1 = \rho_1^{RS} + \mathcal{O}(M_1^2)$, we obtain $\bar p_1 = \mathcal{O}(M_1^4)$ and $\bar h_1 = \mathcal{O}(M_1^4)$, whereas $\hat p_1 = \mathcal{O}(1)$ and $\hat h_1 = \mathcal{O}(1)$.
	We will refer to $\hat{p}_1, \hat{h}_1$ as the fast acoustic pressure and enthalpy, respectively, as they are the main contributors to the acoustic waves, whereas $\bar{p}_1, \bar{h}_1$ vanish as $M_1$ tends to 0. 
	Therefore, in the first steps of the construction of the numerical scheme, we will neglect $\bar p_1, \bar p_2$ and $\bar{h}_1, \bar{h}_2$ obtaining a pure low Mach number scheme. 
	The higher order pressure and enthalpy terms will be included afterwards leading to an all-speed scheme.
	
	\subsection{Time semi-discrete scheme}
	\label{sec:TimeSemi}
	
	We first consider the homogeneous model 
	\begin{equation*}
		\partial_t W + \partial_x \hat f(W) = 0,
	\end{equation*}
	where $\hat f(W)$ denotes the flux function with truncated pressure and enthalpy terms and is given by 
	\begin{equation}
	\label{def:hom_hat_flux}
		\hat f(W) = 
		\begin{pmatrix}
		\rho u \\ \alpha_1 \rho u \\ 
		\displaystyle\alpha_1 \rho_1 u_1 \\
		\displaystyle \alpha_1 \varrho_1\rho_1 u_1^2 + \frac{\alpha_1 \varrho_1 \hat{p}_1}{M_1^2} + \alpha_2 \varrho_2\rho_2 u_2^2 +\frac{\alpha_2\hat p_2}{M_2^2} \\ 
		\displaystyle \frac{u_1^2}{2} - \frac{u_2^2}{2} + \frac{\hat h_1}{M_1^2} - \frac{\hat h_2}{M_2^2}
		\end{pmatrix}
	\end{equation}
	The pressure relaxation and friction source terms will be added in the last step of the construction of the all-speed scheme. 
	In total the equations \eqref{eq:NonDim_model} will be split in four parts since they contain two fast scales connected to the pressures and the enthalpies, stiff relaxation terms and order one terms connected to the material velocity. 
	The challenge in identifying suitable subsystems of the homogeneous system lies firstly in their well-posedness, i.e. constituting hyperbolic systems on their own, and secondly they should be easily solvable.  
	The order in which they are described in the next subsections follows also their order in the final numerical scheme given in Section \ref{sec:FullDisc}.
	
	\subsubsection{Fast acoustics}
	As the pressure terms in the momentum equation are stiff for small Mach numbers, they are treated implicitly. 
	From the flux function \eqref{def:hom_hat_flux}, we see that they couple with the mixture mass flux $\rho u$ which is therefore also treated implicitly. 
	This strategy is also used in the construction of IMEX schemes for single phase isentropic Euler equations \cite{DegTan2011}.
	Therefore we propose to solve the following system implicitly
	\begin{subequations}
		\label{sys:Fast_ac}
	\begin{align}
		\partial_t \rho + \partial_x (\rho u) &= 0 \label{eq:Acc_density}\\
		\partial_t (\rho u) + \partial_x \left(\alpha \varrho_1 \frac{\hat p_1}{M_1^2} + \alpha_2 \varrho_2 \frac{\hat p_2}{M_2^2}\right) &= 0. \label{eq:Acc_mom} 
	\end{align}
	\end{subequations}
    The system is strictly hyperbolic with the following eigenvalues 
    \begin{equation}
    \label{eq:Acc_cRS}
    \lambda^\pm = \pm \sqrt{\frac{\partial\hat{P}}{\partial \rho}} = \pm \sqrt{ \chi_1 \frac{(c_1^{RS})^2}{M_1^2} + \chi_2 \frac{(c_2^{RS})^2}{M_2^2}} = \pm c^{RS}(\chi).
    \end{equation}
    Thereby $c^{RS}(\chi)$ denotes the mixture sound speed with the reference phase sound speeds. 
    The associated eigenvectors are given by 
    \begin{equation*}
        \mathbf{v}^\pm = 
        \begin{pmatrix}
        \displaystyle\pm \frac{1}{c^{RS}(\chi)} \\ 1
        \end{pmatrix}.
    \end{equation*}
    Note that $\hat p_1$ and $\hat p_2$ are linear functions of $\rho_1, \rho_2$, respectively, but are nonlinear in $\rho$ when rewritten in state variables $W$ as follows
	\begin{equation}
	\label{eq:Phase_dens_in_state}
		\rho_1 = \frac{\alpha \rho_1}{\alpha \rho} \rho = \frac{W_3}{W_2} W_1, \quad \rho_2 = \frac{\rho - \alpha \varrho_1 \rho_1}{\varrho_2(\rho - \alpha \rho)} \rho = \frac{W_1 - W_3}{\varrho_2(W_1 - W_2)} W_1.
	\end{equation}
    Considering the phase densities in terms of mixture variables $Q$ defined in \eqref{eq:Mix_var} however leads to a linear dependency on $\rho$ given by 
	\begin{equation}
    \label{eq:Phase_dens_in_mix}
\rho_1 = \frac{\chi}{\varrho_1\alpha} \rho = \frac{Q_3}{\varrho_1 Q_2} Q_1, \quad \rho_2 = \frac{1 - \chi}{\varrho_2(1-\alpha)} \rho = \frac{1 - Q_3}{\varrho_2(1 - Q_2)} Q_1.
    \end{equation}
    Note that the formulation of $\rho_1$ in state or mixture variables coincides whereas the structure of $\rho_2$ differs with respect to $\rho$. 
%    Since system \eqref{sys:Fast_ac} is written in conservative form for $\rho$ which is contained in both the state variables $W$ and mixture variables $Q$, 
    We choose to the linear formulation \eqref{eq:Phase_dens_in_mix} in mixture variables for the phase densities. 
    Then we can write \eqref{eq:Acc_mom} as follows
    \begin{equation}
    \label{eq:Mom_trunc_mix_var}
        \partial_t (\rho u) + \partial_x \left(\frac{\alpha \varrho_1}{M_1^2}\left(p_1^{RS} - \rho_1^{RS}(c_1^{RS})^2\right) +
         \frac{(1-\alpha)\varrho_2}{M_2^2} \left(p_2^{RS} - \rho_2^{RS}(c_2^{RS})^2\right) + \rho ~(c^{RS}(\chi))^2 \right) = 0. 
    \end{equation}
    We discretize \eqref{eq:Acc_density} and \eqref{eq:Mom_trunc_mix_var} implicitly by applying a backward Euler scheme in time. 
    The volume fraction $\alpha$ and mass fraction $\chi$ are not evolved in this step \eqref{sys:Fast_ac} and thus treated explicitly.
    Discrete time increments are defined as $t^{n+1} = t^n + \Delta t$, where $\Delta t$ is is the time dependent time step and obeys a time step restriction given by a CFL condition. 
    Since the backward Euler scheme is unconditionally stable, we do not obtain a CFL restriction from this step of the numerical scheme. 
    Further, we substitute $(\rho u)^{n+1}$ in the density flux by the relation obtained from the momentum update. 
    Let $\rho^\star$ be the update of $\rho$ at the end of the integration of \eqref{eq:Mom_trunc_mix_var} with time step $\Delta t$ starting from data at the previous time $t^n$.
    Then we obtain the following linear implicit equation for $\rho^\star$ given by
    \begin{equation}
    \label{eq:Acc_dens_update_lin}
        \rho^{\star} - \Delta t^2 \partial_x^2 \left( (c^{RS}(\chi^n))^2 ~\rho^{\star}\right) = \rho^n - \Delta t \partial_x (\rho u)^n + \Delta t^2 \partial_x^2\left(\frac{\alpha^n \varrho_1}{M_1^2} \eta_1^{RS} +
        \frac{(1-\alpha^n)\varrho_2}{M_2^2} \eta_2^{RS}\right),
    \end{equation}
    where in $\eta_1^{RS} = p_1^{RS} - \rho_1^{RS}(c_1^{RS})^2$ the reference terms of phase one are collected. 
    Analogously we define $\eta_2^{RS} = p_2^{RS} - \rho_2^{RS}(c_2^{RS})^2$.
    Note that $c^{RS}(\chi^n)$ is always positive and the coefficient matrix after applying a centred differences on the space derivates is a symmetric positive definite matrix.
    For details on the space discretisation we refer to Section \ref{sec:FullDisc}
    Thus, the linear system is well-defined and can be solved efficiently with standard linear solvers. 
    The momentum is then updated in state variables by 
    \begin{equation*}
        (\rho u)^{\star} = (\rho u)^n - \Delta t \partial_x \left(\frac{(\rho\alpha)^n}{\rho^\star} \frac{ \varrho_1\hat{p}_1 (\rho_1^\star)}{M_1^2} + \frac{\rho^\star-(\rho\alpha)^n}{\rho^\star} \frac{\varrho_2\hat{p}_2(\rho_2^{ \star})}{M_2^2}\right),
    \end{equation*}
    $\rho_1^{\star}, \rho_2^\star$ are calculated according to \eqref{eq:Phase_dens_in_state}. 
    The updated state variables after the first stage of the numerical scheme are given by $W^\star = \left(\rho^\star, (\alpha \rho)^n, (\alpha \rho_1)^n, (\rho u)^\star, (u_1 - u_2)^n\right)^T$.

	\subsubsection{Nonlinear transport}
	
%	One is tempted, now that the stiff pressure terms are treated implicitly, to discretize the remaining flux terms explicitly. 
%	In order to verify that this is feasible, we have to check whether the equations with the remaining flux terms are hyperbolic by calculating the Eigenstructure of the resulting flux Jacobian. 
%	In the un-split system \eqref{eq:NonDim_model} the eigenvalues were obtained by transforming the equations into a non-conservative form using the phase variables $v$. 
	
%	Here, instead, we found it impossible to obtain eigenvalues from the formulation in phase variables as well as in mixture variables, as the Jacobian is very complicated (Mathematica, Maxima, Maple were not possible to find eigenvalues).
%	Therefore, we decided to split away the mixture terms leaving us with solely the transport terms and the difference of two Burger's fluxes in the equation for the relative velocity. 
    In the next step we identify nonlinear transport terms of the flux function $\hat{f}$ defined in \eqref{def:hom_hat_flux}.
    In order to introduce the mixture velocity $u$ into the flux function, we rewrite the flux \eqref{def:hom_hat_flux} in terms of state variables which reads
	\begin{equation}
		\label{def:hom_hat_flux_state}
		\hat f(W) = 
		\begin{pmatrix}
		\rho u \\ \alpha_1 \rho u \\ 
		\displaystyle\alpha_1 \rho_1 u + \alpha_1 \varrho_1 \rho_1 \left(1 - \frac{\alpha_1 \varrho_1 \rho_1}{\rho}\right)(u_1 - u_2) \\
		\displaystyle \rho u^2 + \frac{\alpha_1 \varrho_1 \hat{p}_1}{M_1^2} + \frac{\alpha_2\hat p_2}{M_2^2} + \alpha_1 \varrho_1 \rho_1 \left(1 - \frac{\alpha_1 \varrho_1 \rho_1}{\rho}\right)(u_1 - u_2)^2\\ 
		\displaystyle u (u_1 - u_2) + \left( \frac{\hat h_1}{M_1^2} - \frac{\hat h_2}{M_2^2}\right) + \left(1 - 2 \frac{\alpha_1 \varrho_1 \rho_1}{\rho}\right)\frac{(u_1 - u_2)^2}{2}
		\end{pmatrix}.
	\end{equation}
	Considering that the mass flux was already treated in the previous step, the remaining transport terms are given by
	\begin{align}
	\label{sys:Transport}
	\begin{split}
	\partial_t \rho  &= 0, \\
	\partial_t (\alpha \rho) + \partial_x (\alpha \rho u) &= 0, \\
	\partial_t (\alpha \rho_1) + \partial_x (\alpha \rho_1 u) & = 0, \\
	\partial_t (\rho u) + \partial_x (\rho u^2) &= 0, \\
	\partial_t (u_1 - u_2) + \partial_x \left( u (u_1 - u_2) + \left(1 - 2 \frac{\alpha_1 \varrho_1 \rho_1}{\rho}\right)\frac{(u_1 - u_2)^2}{2}\right)&= 0.
	\end{split}
	\end{align}
	This system is hyperbolic with the following eigenvalues
	\begin{equation}
	\label{eq:EW_transport}
		\lambda_1 = 0, \quad \lambda_{2,3} = u, \quad \lambda_4 = 2u, \quad \lambda_5 = u_1 + u_2 - u.
	\end{equation}
	The associated eigenvectors read
	\begin{equation*}
		\mathbf{v}_1 = 
		\begin{pmatrix}
		1 \\ \frac{\alpha}{2} \\ \frac{u}{2} \\ \frac{\chi_1}{2} \\
		\frac{\varepsilon_1}{2}
		\end{pmatrix}, \quad
		\mathbf{v}_2 =
		\begin{pmatrix}
		0 \\ 1 \\ 0 \\ 0 \\ 0
		\end{pmatrix}, \quad
		\mathbf{v}_3 =
		\begin{pmatrix}
		0 \\ 0 \\ 0 \\ 1 \\ \varepsilon_2
		\end{pmatrix}, \quad
		\mathbf{v}_4 = 
		\begin{pmatrix}
		0 \\ 1 \\ \frac{\chi_1}{\alpha_1} \\ \frac{u}{\alpha_1} \\ \frac{\varepsilon_3}{\alpha_1}
		\end{pmatrix}, \quad
		\mathbf{v}_5 = 
		\begin{pmatrix}
		0\\0\\0\\0\\1
		\end{pmatrix}
	\end{equation*}
	with 
	\begin{equation}
		\varepsilon_1 = - \frac{(u_1 - u_2)(\chi_1 (u_1 - u_2) - u)}{-(1-2\chi_1)\rho (u_1 - u_2) + \rho u}, \quad \varepsilon_2 = \frac{(u_1 - u_2)}{(1-2 \chi_1)\rho}, \quad \varepsilon_3 = \frac{(u_1 - u_2)( \chi_1 (u_1 - u_2) - u)}{-(1 - 2\chi_1) \rho (u_1 - u_2) - \rho u}.
	\end{equation}
	Note that the double eigenvalue $u$ has two linearly independent eigenvectors and all eigenvectors are linearly independent. 
	Since all waves exhibited by system \eqref{sys:Transport} are of the order of the material velocity, we discretize the equations \eqref{sys:Transport} explicitly.
	Let $W^{\star\star}$ denote the new state after the advection step.
	Then the time discretization of \eqref{sys:Transport} is given by 
	\begin{align}
	\label{eq:Time_semi_transport}
		\begin{split}
		\rho^{\star\star}  &= \rho^\star, \\
		(\alpha \rho)^{\star\star} &= (\alpha\rho)^n - \Delta t \partial_x \left(\frac{(\alpha \rho)^n (\rho u)^\star}{\rho^\star}\right), \\
	    (\alpha \rho_1)^{\star\star} &= (\alpha\rho_1)^n - \Delta t \partial_x \left(\frac{(\alpha \rho_1)^n (\rho u)^\star}{\rho^\star}\right), \\
		(\rho u)^{\star\star} &= (\rho u)^\star -\Delta t \partial_x \left(\frac{(\rho u)^\star (\rho u)^\star}{\rho^\star}\right), \\
		(u_1 - u_2)^{\star\star} &= (u_1 - u_2)^n - \Delta t \partial_x \left(\frac{ (\rho u)^\star (u_1 - u_2)^n}{\rho^\star} + \left(1 - 2\varrho_1 \frac{(\alpha\rho_1)^n}{\rho^\star}\right)\frac{((u_1 - u_2)^n)^2}{2}\right).
		\end{split}
	\end{align}
	Therein we have used the values at time $t^n$ for the state variables $(\alpha \rho), (\alpha_1\rho_1)$ and $u_1 - u_2$ since they were not evolved in the acoustic stage \eqref{sys:Fast_ac}.
	Due to the explicit time integration, the following CFL condition 
		\begin{equation}
		\label{eq:CFL_mat_transport}
			\Delta t \leq \nu_{u} \frac{\Delta x}{\max(2|u|,|u_1 + u_2 - u|)}
		\end{equation}
	has to be met which is independent of the Mach numbers and follows the material wave. 
	 
%	\begin{remark}
%	From the perspective of eigenvalues only, the difference of the enthalpy could be included into the transport part. 
%	However, this will not lead to a stable scheme in practise, as the enthalpy terms are not coupled to the rest of the equations unlike in the original equations \eqref{eq:NonDim_model}. (Explain better!)
%	\end{remark}
	
	\subsubsection{Stiff mixture terms}
	Collecting the remaining terms in $\hat{f}$, that were not yet treated, results in the following subsystem 
	\begin{subequations}
	\label{eq:mixture_terms}
	\begin{align}
	\partial_t (\alpha \rho_1)  + \partial_x \left(\alpha \varrho_1 \rho_1 \left(1 - \frac{\alpha_1 \varrho_1 \rho_1}{\rho}\right)(u_1 - u_2)\right) & = 0, \label{eq:Mixture_terms_alphaRho1}\\
	\partial_t (\rho u) + \partial_x \left(\alpha \varrho_1 \rho_1 \left(1 - \frac{\alpha_1 \varrho_1 \rho_1}{\rho}\right)(u_1 - u_2)^2\right)&= 0, \label{eq:Mixture_terms_Mom}\\
	\partial_t (u_1 - u_2) + \partial_x \left( \frac{\hat h_1}{M_1^2} - \frac{\hat h_2}{M_2^2}\right)&= 0. \label{eq:Mixture_terms_relVel}
	\end{align}
	\end{subequations}
	This system is strictly hyperbolic with the eigenvalues 
	\begin{equation}
	\label{eq:EW_mixture}
	\lambda^0 = 0, \quad \lambda^\pm = \frac{\varrho_1}{2}\left( 1-2 \chi_1\right)(u_1 - u_2) \\
	\pm\frac{1}{2}\sqrt{\mathcal{A}~}
	\end{equation}
	with 
	\begin{align*}
		\mathcal{A} &= \varrho_1^2\left( 1 - 2 \chi_1\right)^2 {(u_1 - u_2)^{2}} + 4 \rho \chi \left( 1- \chi\right) \left(\frac{1}{\alpha_1 M_1^2}\frac{(c_1^{RS})^2}{\rho_1^{RS}} + \frac{1}{\alpha_2 M_2^2}\frac{(c_2^{RS})^2}{\rho_2^{RS}}\right).
	\end{align*}
	Since $\mathcal{A}$ is positive, the eigenvalues are real.
	The associated eigenvectors are given by 
	\begin{equation*}
		\mathbf{v}^0 = 
		\begin{pmatrix}
		0 \\ 1 \\ 0
		\end{pmatrix}, \quad 
		\mathbf{v}^\pm = 
		\begin{pmatrix}
		\displaystyle\frac{1}{2 \rho \mathcal{C}^{RS}} \left(\varrho_1 (u_1- u_2)\alpha_1 \alpha_2 \rho_1^{RS}\rho_2^{RS}(1-2\chi) \pm \sqrt{\mathcal{A}}\right) \\[10pt] 
       \displaystyle- \frac{w}{\rho \mathcal{C}^{RS}} \left((\varrho_1 - \varrho_2) (u_1- u_2)\alpha_1 \alpha_2 \rho_1^{RS}\rho_2^{RS}(1-2\chi) \pm \sqrt{\mathcal{A}}\right) \\ 1
		\end{pmatrix},
	\end{equation*} 
	where 
    \begin{equation*}
        \mathcal{C}^{RS} = \frac{\alpha \rho_1^{RS}}{\rho} (c_1^{RS})^2  + \frac{\alpha_2 \rho_2^{RS}}{\rho} (c_2^{RS})^2 . 
    \end{equation*}
	
	From the eigenstructure of \eqref{eq:mixture_terms} we can deduce that the relative velocity equation \eqref{eq:Mixture_terms_relVel} is coupled with the partial density equation \eqref{eq:Mixture_terms_alphaRho1} via the enthalpy terms $\hat{h}_1, \hat{h}_2$ resulting in the Mach number dependent eigenvalues $\lambda^\pm$. 
	Therefore, the enthalpy terms must be treated implicitly to obtain an overall Mach number independent CFL restriction. 
	The momentum equation \eqref{eq:Mixture_terms_Mom} is decoupled from \eqref{eq:Mixture_terms_alphaRho1} and \eqref{eq:Mixture_terms_relVel} and can be updated directly once the partial density $\alpha_1 \rho_1$ and relative velocity $u_1 - u_2$ are obtained. 
    This does not yield a restriction on the time step since $\rho u$ is associated with a zero eigenvalue.
	Therefore we consider the subsystem consisting of \eqref{eq:Mixture_terms_alphaRho1} and \eqref{eq:Mixture_terms_relVel} implicitly and obtain with the backward Euler scheme, starting from stage $W^{\star\star}$, the following time discretization 
	\begin{align}
	(\alpha \rho_1)^{\star\star\star} &=  (\alpha \rho_1)^{\star\star}  - \Delta t \partial_x \left(\varrho_1 (\alpha \rho_1)^{\star\star\star} \left(1 - \varrho_1\frac{(\alpha \rho_1)^{\star\star\star}}{\rho^{\star\star}}\right)(u_1 - u_2)^{\star\star\star}\right), \label{eq:alphaRho1_mix}\\
	(u_1 - u_2)^{\star\star\star} &= (u_1 - u_2)^{\star\star} - \Delta t \partial_x \left( \frac{(c_1^{RS})^2}{\rho_1^{RS}  M_1^2} \frac{(\alpha \rho_1)^{\star\star\star}}{\alpha^{\star\star}}  - \frac{(c_2^{RS})^2}{\rho_2^{RS} M_2^2} \frac{\rho^{\star\star} - \varrho_1 (\alpha \rho_1)^{\star\star\star}}{1 - \alpha^{\star\star}} \right). \label{eq:relVeloc_mix_update}
	\end{align}
	In \eqref{eq:relVeloc_mix_update} we have already inserted the definitions of $\hat h_1, \hat h_2$ given in \eqref{eq:Def_h_hat}.
	Substituting $(u_1 - u_2)^{\star\star\star}$ in \eqref{eq:alphaRho1_mix} by the relation given in \eqref{eq:relVeloc_mix_update}, we obtain
	\begin{align}
	\begin{split}
	\label{eq:Update_AlphaRho1}
		(\alpha \rho_1)^{\star\star\star} - (\alpha \rho_1)^{\star\star}  &+ \Delta t \partial_x \left(\varrho_1 (\alpha \rho_1)^{\star\star\star} \left(1 - \frac{\varrho_1(\alpha \rho_1)^{\star\star\star}}{\rho^{\star\star}}\right)(u_1 - u_2)^{\star\star}\right) 
\\
		 &- \Delta t^2 \partial_x \left\lbrace \varrho_1 (\alpha \rho_1)^{\star\star\star} \left(1 -  \frac{\varrho_1 (\alpha \rho_1)^{\star\star\star}}{\rho^\star}\right) \times \right.\\
		 &\left.\partial_x \left( \left(\frac{(c_1^{RS})^2}{\rho_1^{RS}  M_1^2} \frac{1}{\alpha_1^{\star\star}} + \frac{(c_2^{RS})^2}{\rho_2^{RS} M_2^2} \frac{1}{\alpha_2^{\star\star}} \right) (\alpha \rho_1)^{\star\star\star} - \frac{(c_2^{RS})^2}{\rho_2^{RS} M_2^2} \frac{\rho^{\star\star}}{\alpha_2^{\star\star}}\right)\right\rbrace = 0.
	\end{split}
	\end{align}
	It is a nonlinear elliptic problem in $(\alpha \rho_1)^{\star\star\star}$ that can be solved with a linearised iterative scheme as follows
	\begin{align}
	\begin{split}
\label{eq:Iterative_alphaRho1}
(\alpha \rho_1)^{l+1} - (\alpha \rho_1)^{\star\star}  &+ \Delta t \partial_x \left(\varrho_1 (\alpha \rho_1)^{l} \left(1 - \frac{\varrho_1(\alpha \rho_1)^{l}}{\rho^{\star\star}}\right)(u_1 - u_2)^{\star\star}\right) 
\\
&- \Delta t^2 \partial_x \left\lbrace \varrho_1 (\alpha \rho_1)^{l} \left(1 -  \frac{\varrho_1 (\alpha \rho_1)^{l}}{\rho^\star}\right) \right. \times \\
&\left.\partial_x \left( \left(\frac{(c_1^{RS})^2}{\rho_1^{RS}  M_1^2} \frac{1}{\alpha_1^{\star\star}} + \frac{(c_2^{RS})^2}{\rho_2^{RS} M_2^2} \frac{1}{\alpha_2^{\star\star}} \right) (\alpha \rho_1)^{l+1} - \frac{(c_2^{RS})^2}{\rho_2^{RS} M_2^2} \frac{\rho^{\star\star}}{\alpha_2^{\star\star}}\right)\right\rbrace = 0,
\end{split}
\end{align}
where  $l \in \mathbb{N}_0$ denotes the number of iteration and we set the start value as $(\alpha \rho_1)^0 = (\alpha\rho_1)^{\star\star}$.
When the stopping criterion given by the relative $L^1$ error
\begin{equation*}
    \frac{\|{\alpha\rho_1}^{l+1} - (\alpha\rho_1)^l\|_{L^1}}{\| (\alpha \rho_1)^l\|_{L^1}} < \delta,
\end{equation*}
for a given $\delta > 0$ is fulfilled, then we set $(\alpha\rho_1)^{\star\star\star} = (\alpha\rho_1)^{l+1}$. 
We note that after applying a suitable space discretization, the coefficient matrix for solving $(\alpha\rho_1)^{l+1}$ is symmetric and positive definite, since the terms 
\begin{equation}
\label{eq:Coeff_mixture_update}
   s(W) = \varrho_1 (\alpha \rho_1)^{l} \left(1 -  \frac{\varrho_1 (\alpha \rho_1)^{l}}{\rho^\star}\right) \text{ and } \tilde{s}(W) =  \frac{(c_1^{RS})^2}{\rho_1^{RS}  M_1^2} \frac{1}{\alpha_1^{\star\star}} + \frac{(c_2^{RS})^2}{\rho_2^{RS} M_2^2} \frac{1}{\alpha_2^{\star\star}}
\end{equation}
are positive. 
The linear system is therefore well-posed and can be solved efficiently with standard linear solvers. 
The relative velocity is then updated explicitly according to \eqref{eq:relVeloc_mix_update} and the momentum by 
\begin{equation*}
    (\rho u)^{\star\star\star} = (\rho u)^{\star\star} - \Delta t \partial_x (\rho^{\star\star} \chi_1^{\star\star\star} \chi_2^{\star\star\star}((u_1 - u_2)^{\star\star\star})^2).
\end{equation*}
The updated state vector after the third step of the time semi-discrete scheme is then given by $W^{\star\star\star} = (\rho^{\star\star}, (\alpha\rho)^{\star\star},(\alpha_1\rho_1)^{\star\star\star}, (\rho u)^{\star\star\star}, (u_1 - u_2)^{\star\star\star})^T$.

\subsubsection{Treatment of higher order pressure and enthalpy terms}
Till now we have considered the truncated pressure and enthalpy terms \eqref{eq:Def_p_hat}, \eqref{eq:Def_h_hat}. 
This is sufficient to construct a scheme that is applicable on low Mach number flows only. 
To obtain accurate results also for compressible flow regimes, we have to take into account the higher order terms in the Taylor expansions \eqref{eq:p_bar}, \eqref{eq:h_bar}.    
In those regimes, $\bar p_k, \bar h_k$ are of order $\mathcal{O}(1)$ and neglecting them results in large errors and spurious results. 
As these higher order terms vanish for well-prepared initial data in the low Mach number limit, they can be treated explicitly without generating a dependence of the Mach numbers on the CFL condition.
Moreover $\bar{p}$ and $\bar{h}$ are nonlinear and treating them explicitly does not add to the complexity of the scheme. 
Nevertheless, they have to be added carefully to the existing low Mach number scheme to ensure the hyperbolic structure of the subsystems \eqref{sys:Fast_ac}, \eqref{sys:Transport} and \eqref{eq:mixture_terms}. 
We add them in the third step of the time semi-discrete scheme, which is then given by 
%	Including the $\bar{h}$ terms in the $(u_1-u_2)$ equation adds the term 
%	\begin{equation}
%	\label{eq:Hbar_eigenvalue}
%	\frac{\partial_{\alpha \rho_1} \bar{h}_1}{M_1^2} - \frac{\partial_{\alpha \rho_1} \bar{h}_2}{M_2^2} = \frac{1}{\alpha M_1^2}\left(\frac{c_1^2}{\rho_1} - \frac{(c_1^{RS})^2}{\rho_1^{RS}}\right) + \frac{1}{(1 -\alpha) M_2^2}\left(\frac{c_2^2}{\rho_2} - \frac{(c_2^{RS})^2}{\rho_2^{RS}}\right),
%	\end{equation}
%	under the square root of $\lambda^\pm$ in \eqref{eq:EW_mixture}. 
%	Using a Taylor expansion of $c_{1}^2$ in $\rho_{1} - \rho_{1}^{RS}$ (analogously for $c_2^2$) and assuming well-prepared data, we can simplify \eqref{eq:Hbar_eigenvalue}
%	\begin{equation}
%		\frac{\partial_{\alpha \rho_1} \bar{h}_1}{M_1^2} - \frac{\partial_{\alpha \rho_1} \bar{h}_2}{M_2^2} = \frac{(\gamma_1 -2) (c_1^{RS})^2}{\rho_1^{RS} \rho_1}\rho_1^{(2)} + \frac{(\gamma_2 -2) (c_2^{RS})^2}{\rho_2^{RS} \rho_2}\rho_2^{(2)} + \mathcal{O}(M_1^2,M_2^2).
%	\end{equation}
%	If $M_1, M_2 \approx 1$, then we can set the reference values such that \eqref{eq:Hbar_eigenvalue} is positive, as the advective CFL restriction and the mixture CFL restriction are of the same order. 
%	In addition, in the compressible case, one might actually be interested in the resolution of the acoustic waves thus choosing a more restrictive acoustic CFL condition. 
	\begin{enumerate}
		\item Calculate $(\alpha \rho_1)^{\star\star\star}$ by solving \eqref{eq:Iterative_alphaRho1} iteratively.
		\item Update $u_1 - u_2$ using the untruncated enthalpies $h_k = \hat{h}_k + \bar{h}_k ~(k=1,2) $ by
		\begin{equation}
		\label{eq:RelVeloc_update_allspeed}
			(u_1 - u_2)^{\star\star\star} = (u_1 - u_2)^{\star\star}
			- \Delta t \partial_x \left(\frac{h_1(\rho_1^{\star\star\star})}{M_1^2} - \frac{h_2(\rho_2^{\star\star\star})}{M_2^2} \right),
		\end{equation}
		where $\rho_1$ and $\rho_2$ are calculated by relation \eqref{eq:Phase_dens_in_state}. 
		\item Update $\rho u$ including $\bar{p}_k$ terms $(k=1,2)$ by 
		\begin{multline}
		\label{eq:Mom_update_allspeed}
		(\rho u)^{\star\star\star} = (\rho u)^{\star\star} - \Delta t \partial_x \left((\alpha \rho_1)^{\star\star\star} (1 - \frac{(\alpha \rho_1)^{\star\star\star}}{\rho^\star})((u_1 - u_2)^{\star\star\star})^2\right) \\
		- \Delta t\partial_x \left( \alpha^{\star\star}\frac{\bar{p}_1\left(\rho_1^{\star\star\star}\right)}{M_1^2} + \left(1-\alpha^{\star\star}\right) \frac{\bar{p}_2(\rho_2^{\star\star\star})}{M_2^2}\right).
		\end{multline}
	\end{enumerate}
    This results in a correction of the momentum and relative velocity for compressible flow regimes leading to an all-speed scheme for the simulation of two phase flows.  
	
	\subsubsection{Treatment of the relaxation source terms}
	After having established the main scheme for the homogeneous case, we focus now on the numerical approximation of the relaxation source terms acting on the volume fraction and relative velocity. 
    They are given by
	\begin{align}
		\partial_t \rho\alpha &= - \frac{1}{\tau}\left(\frac{\varrho_2}{M_2^2} ~p_2\left(\rho_2\right) - \frac{\varrho_1}{M_1^2}~p_1(\rho_1)\right), \label{eq:Pressure_relax}\\
		\partial_t (u_1 - u_2) &= - \zeta \frac{\varrho_1\alpha\rho_1}{\rho}\left(1 - \frac{\varrho_1\alpha \rho_1}{\rho}\right)(u_1 - u_2) \label{eq:Friction}.
	\end{align}
    We first notice that both equations are decoupled and can be solved simultaneously. 
    Since the remaining variables are not affected by the relaxation processes, the state vector at time $t^{n+1}$ is given by $W^{n+1} = (\rho^{\star\star\star}, \alpha^{n+1}\rho^{\star\star\star}, (\alpha\rho_1)^{\star\star\star}, (\rho u)^{\star\star\star}, (u_1 - u_2)^{n+1})^T$. 
    Since the friction parameter $\zeta > 0$ can be large, equation \eqref{eq:Friction} is integrated implicitly. 
    Due to the linearity of the source term, we can find an analytic update of the relative velocity given by
    \begin{equation}
    \label{eq:Update_Friction}
    \left(1 + \Delta t\zeta \frac{\varrho_1(\alpha\rho_1)^{n+1}}{\rho^{n+1}}\left(1 - \frac{\varrho_1(\alpha \rho_1)^{n+1}}{\rho^{n+1}}\right)\right)(u_1 - u_2)^{n+1} = (u_1 - u_2)^{\star\star\star}.
    \end{equation}
	The pressure relaxation rate is a non-negative parameter $\tau \in  [0, \infty)$, where $\tau = 0$ gives an instantaneous pressure relaxation. 
    The homogeneous equation corresponds to $`` \tau = \infty "$. 
	For fast relaxation rates with $\tau \ll 1$, equation \eqref{eq:Friction} becomes stiff and is discretized implicitly. 
    Since we are free to choose the set of variables, we rewrite \eqref{eq:Friction} in terms of mixture variables $Q$ which is consistent with model \eqref{eq:Comp_model}. 
    We obtain the following implicit equation
	\begin{align}
	\label{eq:Relax_Update_alpha}
		\alpha^{n+1} &= \alpha^{\star\star\star} - \frac{\Delta t}{\tau \rho^{n+1}} \left(\frac{\varrho_2}{M_2^2}\ p_2\left(\frac{1 - \chi^{n+1}}{1-\alpha^{n+1}} \rho^{n+1}\right) -\frac{\varrho_1}{M_1^2}\  p_1\left(\frac{\chi^{n+1}}{\alpha^{n+1}}\rho^{n+1}\right)\right).
	\end{align}
	Since the phase pressures $p_1, p_2$ are nonlinear in $\alpha$ we use the Newton method applied on $g(\alpha) = 0$ to obtain $ \alpha^{n+1}$ with 
	\begin{align*}
		g(\alpha) &= - \alpha - \frac{\Delta t}{\tau \rho^{n+1}} \left(\frac{\varrho_2}{M_2^2} p_2(\alpha) - \frac{\varrho_1}{M_1^2} p_1(\alpha)\right) + \alpha^{\star\star\star}, \\
		\frac{\partial g}{\partial \alpha} &= - 1 - \frac{\Delta t}{\tau} \left(\frac{\varrho_2}{M_2^2} c_2(\alpha)^2 \frac{1-\chi^{n+1}}{1 - \alpha} + \frac{\varrho_1}{M_1^2} c_1(\alpha)^2 \frac{\chi^{n+1}}{\alpha}\right),
	\end{align*}
	where the sound speeds $c_1, c_2$ also depend on $\alpha$. 
	The focus of this work is on the construction of an asymptotic preserving scheme that is consistent with the singular Mach number limits.
	For a different approach on how to treat the relaxation source terms which was used in the context of the Baer-Nunziato model we refer the interested reader to \cite{ChiMue2020}.
	
	To summarize, we have split the non-dimensional model \eqref{eq:NonDim_model} in the following way
	\begin{equation*}
	\partial_t W + \partial_x f^a\left(W\right) + \partial_x f^b\left(W\right) + \partial_x f^c\left(W\right) =- r\left(W\right),  
	\end{equation*}
	where the fluxes of the subsystems \eqref{sys:Fast_ac}, \eqref{sys:Transport}, \eqref{eq:mixture_terms} are given by 
	\begin{align*}
	f^a(W) &= 
	\begin{pmatrix}
	\rho u \\ 0 \\ 0 \\ \frac{\alpha \hat p_1}{M_1^2} + \frac{(1-\alpha)\hat p_2}{M_2^2} \\ 0
	\end{pmatrix}, \quad
	f^b(W) = 
	\begin{pmatrix}
	0 \\ \alpha \rho u \\ \alpha \rho_1 u \\ \rho u^2 \\ \frac{1}{2}(u_1^2 - u_2^2)
	\end{pmatrix}, \\ 
	f^c(W) &= 
	\begin{pmatrix}
	0 \\ 0 \\ \alpha \rho_1 \left(1 - \frac{\alpha \rho_1}{\rho}\right) (u_1 - u_2) \\
	\alpha \rho_1 \left(1 - \frac{\alpha \rho_1}{\rho}\right) (u_1 - u_2)^2 + \frac{\alpha \bar p_1}{M_1^2} + \frac{(1-\alpha)\bar p_2}{M_2^2} \\
	\frac{h_1}{M_1^2} - \frac{h_2}{M_2^2}
	\end{pmatrix}
	\end{align*}
	and the relaxation source term reads
	\begin{equation*}
	r(W) = 
	\begin{pmatrix}
	0 \\ \frac{1}{\tau} \left(\frac{p_2}{M_2^2} - \frac{p_1}{M_1^2}\right) \\ 0\\ 0 \\ \zeta \chi \left(1-\chi\right)(u_1 - u_2)
	\end{pmatrix}.
	\end{equation*}
	The time semi-discrete scheme is then given by the following operator splitting
	\begin{subequations}
		\label{eq:Time_Semi_scheme}
		\begin{align}
		W^\star &= W^n - \Delta t ~\partial_x f^a\left(W^\star\right), \label{eq:Time_Semi_Acc}\\
		W^{\star\star} &= W^\star - \Delta t ~\partial_x f^b\left(W^\star\right), \label{eq:Time_Semi_Adv}\\ 
		W^{\star\star\star} &= W^{\star\star} - \Delta t ~\partial_x f^c\left(W^{\star\star\star}\right), \label{eq:Time_Semi_Mix}\\
		W^{n+1} &= W^{\star\star\star} - \Delta t ~r\left(W^{n+1}\right). \label{eq:Time_Semi_Relax}
		\end{align}
	\end{subequations}

	\subsection{Fully discrete scheme}
    \label{sec:FullDisc}
    
	In this section, we describe the fully discrete numerical scheme associated with the time semi-discrete stages derived in the previous section. 
	In time, we set as before $t^{n+1} = t^n + \Delta t$, where $\Delta t$ obeys a time step restriction given by the CFL condition \eqref{eq:CFL_mat_transport}.
	In space, we consider a computational domain $\Omega$ divided into cells $C_i = (x_{i-1/2}, x_{i+1/2})$ of uniform step size $\Delta x$ with the cell center $x_i = i \Delta x$ for $i=1,\dots,N$.
	We use a finite volume framework, where the solution on cell $C_i$ at time $t^n$ is approximated by the average given by 
	\begin{equation*}
	    W_i^n \approx \frac{1}{\Delta x} \int_{\Omega_i} w(x,t^n) dx.
	\end{equation*}
	For the explicit advective step \eqref{eq:Time_Semi_Adv}, we apply a standard finite volume scheme using the Rusanov flux.
	The update on cell $C_i$ is given by 
	\begin{equation*}
	W_i^{\star\star} = W_i^\star - \frac{\Delta t}{\Delta x}\left( F^a(W_i^\star, W_{i+1}^\star) - F^a(W_{i-1}^\star, W_i^\star)\right),
	\end{equation*}	 
	where the numerical flux function $F^a(W_i,W_{i+1})$ is given by 
	\begin{equation*}
	F^a(W_i,W_{i+1}) = \frac{1}{2} (f^a(W_i) + f^a(W_{i+1}) - \frac{1}{2} a_{i+1/2} \left(W_{i+1} - W_i\right)
	\end{equation*}
	with $a_{i+1/2} = \max_{k=1,\dots,5}(|\lambda_k(W_i,W_{i+1})|)$.
	The eigenvalues $\lambda_k$ are given in \eqref{eq:EW_transport} by the advection step and are of the order of the material wave. 
	Since this step is the only explicit part of the numerical scheme, the CFL condition is given by \eqref{eq:CFL_mat_transport}.
	For the implicit steps \eqref{eq:Time_Semi_Acc} and \eqref{eq:Time_Semi_Mix} we apply centered differences for the space derivatives.        
	Thereby, the mixed derivatives \eqref{eq:Coeff_mixture_update} in \eqref{eq:Iterative_alphaRho1} are discretized by 
	\begin{align*}
	\partial_x \left( s(W)\ \partial_x \left(\tilde{s}(W) W \right) \right) \cong &~ \frac{1}{\Delta x^2}\left( s_{i+1/2}\tilde{s}(W_{i+1}) W_{i+1} - (s_{i+1/2} + s_{i-1/2})\tilde{s}(W_i) W_i \right.\\
	& \left. \qquad \quad+ ~s_{i-1/2}\tilde{s}(W_{i-1}) W_{i-1}\right), \\
	s_{i+1/2} \cong &~ \frac{1}{2}(s(W_i) + s(W_{i+1})).
	\end{align*}
	Note that by construction, the diffusion of the all-speed scheme is independent of the Mach number due to the use of centred differences on the stiff pressure and enthalpy terms and Mach number independent eigenvalues in the transport step.
	The scheme is therefore well suited to simulate flows in the low Mach number regime.

%	For the (semi) explicit updates in the acoustic and mixture parts \eqref{eq:Mom_update_acc}, \eqref{eq:Mom_update_allspeed}, \eqref{eq:RelVeloc_update_allspeed} centred fluxes are applied as follows
%	\begin{align}
%	W_i^{\star} &= W_i^n - \frac{\Delta t}{2\Delta x} \left(f^a(W_{i+1}^{n};\rho_{i+1}^\star) - f^a(W_{i-1}^n;\rho_{i-1}^{\star})\right), \\
%	W_i^{\star\star\star} &= W_i^{\star\star} - \frac{\Delta t}{2\Delta x} \left(f^c(W_{i+1}^{\star\star};(\alpha\rho_1)_{i+1}^{\star\star\star}) - f^c(W_{i-1}^{\star\star};(\alpha\rho_1)_{i-1}^{\star\star\star})\right).
%    \end{align}	 
%    The values $\rho_i^\star$ and $(\alpha\rho_1)_i^{n+1}$ are obtained with via the implicit linear iterative elliptic equations \eqref{eq:AcIteration} and \eqref{eq:Iterative_alphaRho1} with the second order centred differences in space given above.
%    In case extra diffusion is needed to make the numerical simulation more robust, especially when the eigenvalues $\lambda_k$ of the advection step are small, extra diffusion can be added in the (semi) explicit update in the acoustic and mixture parts.
	
	\section{AP property}
	As we have seen in Section \ref{sec:WP}, for well-prepared initial data, the continuous compressible model converges formally towards incompressible equations when the Mach number tends to zero. 
	For the all-speed RS-IMEX scheme we show the discrete analogue, i.e. the numerical scheme is a consistent discretization of the limit model in the low Mach number limit. 
	This property is called asymptotic preserving (AP) and it is essential for ensuring the correct behaviour of the numerical solution in low Mach number regimes. 
	
	We will show this property for the time semi-discrete scheme. Indeed, to obtain a consistent discretization with the limit equations an appropriate time discretization is essential.
	Thereby we use techniques that are used in the context of proving the AP property of IMEX schemes for (isentropic) Euler equations, see e.g. \cite{CorDegKum2012,DegTan2011}.
	First, we consider the case of a compressible and a weakly compressible phase (Case 1 in Section \ref{sec:WP}), then the case of two weakly compressible phases in the same low Mach number regime, i.e. $M_1 = M_2 = M \ll 1$ (Case 2 in Section \ref{sec:WP}).
	As done in the derivation of the well-prepared data, we neglect the scaling parameters $\varrho_1,\varrho_2$.
	
	\paragraph{Case 1}
	Let the data at time $t^n$ be well-prepared, i.e. $W^n \in \Omega^{wp}_1$ given in \eqref{eq:wp_1} with constant reference state $\rho_2^{(0)} = \rho_2^{RS}$ and $\alpha^{(0)}$ constant. 
	The Mach number expansions for the variables of the second phase is then given by
	\begin{align}
	\begin{split}
	\label{eq:Mach_num_exp_Case1}
	\alpha^n &= \alpha^{(0)} + \mathcal{O}(M), \quad \alpha^{(0)} \text{ const. }\\
	\rho_2^n &= \rho_2^{RS} + M^2 \rho_2^{(2),n} + \mathcal{O}(M^3), \\
	u_2^n &= u_2^{(0),n} + \mathcal{O}(M) \text{ with } \partial_x u_2^{(0),n} = 0 \\
    p_2^{RS} &= 0, ~p_2^{(2)} = p_1.
    \end{split}
	\end{align} 
	In particular, we find 
    \begin{equation*}
        \chi^n = \frac{\alpha^{(0)}\rho_1^{RS}}{\rho^{(0)}} + M \frac{\alpha^{(1)}\rho_1^{RS}\rho^{(0)} - \alpha^{(0)}\rho_1^{RS}\rho^{(1)}}{{\rho^{(0)}}^2} + \mathcal{O}(M^2) = \chi^{(0),n} + M \chi^{(1),n} + \mathcal{O}(M^2).
    \end{equation*} 
	We assume for all sub-steps of the numerical scheme that we have the following Mach number expansions of phase two
	\begin{align*}
	\begin{split}
	\alpha^\star &= \alpha^{(0),\star} + \mathcal{O}(M), \\
	\rho_2^\star &= \rho_2^{(0),\star} + M \rho_2^{(1),\star} + M^2 \rho_2^{(2),\star} + \mathcal{O}(M^3), \\
	u_2^\star &= u_2^{(0),\star} + \mathcal{O}(M).
	\end{split}
	\end{align*}
	An analogue notation holds for the sub-steps $(\cdot)^{\star\star}, (\cdot)^{\star\star\star}$ and the final update $(\cdot)^{n+1}$. 
	Considering the update of the density $\rho^\star$ given in \eqref{eq:Acc_density}, we can rewrite it with $p_2^{RS} = 0$ as  
	\begin{equation}
	\label{eq:Update_dens_Case1}
		\rho^\star - \Delta t^2 \partial_x^2 \left(\frac{(c_2^{RS})^2}{M^2}\chi_2^n \rho^\star\right) = \rho^n - \Delta t \partial_x (\rho u)^n + \Delta t^2 \partial_x^2 \left(\alpha^n \hat p_1^\star -  \frac{(c_2^{RS})^2}{M^2} \alpha_2^n \rho_2^{RS}  \right).
	\end{equation}
	Inserting the Mach number expansions \eqref{eq:Mach_num_exp_Case1} and sorting by order of the Mach numbers, we obtain for the $\mathcal{O}(M^{-2})$ terms
	\begin{equation*}
		\partial_x^2 (\chi_2^{(0),n} \rho^{(0),\star}) = 0 \Leftrightarrow \partial_x^2 \left(\frac{\chi_2^{(0),n}}{\alpha^{(0),n}} \rho^{(0),\star} \right) = 0 \Leftrightarrow \partial_x^2 \rho_2^{(0),\star} = 0.
	\end{equation*}
	The last equivalence holds since we solve \eqref{eq:Update_dens_Case1} in mixture variables $\rho, \alpha, \chi$ and $\alpha^\star = \alpha^n, ~\chi^\star = \chi^n$. 
	For periodic or no-flux boundary conditions we formally obtain that $\rho_2^{(0),\star}$ is constant.
	For the $\mathcal{O}(M^{-1})$ terms, we find 
	\begin{equation*}
		\partial_x^2 \left(\frac{\rho_2^{RS}}{\rho^{(0),n}} \rho^{(1),\star} - \frac{ \rho_2^{RS}\rho^{(1),n}}{{\rho^{(0),n}}^2}\rho^{(0),\star}\right) = 0 \Leftrightarrow \partial_x^2 \rho_2^{(1),\star} = 0.
	\end{equation*}
	This implies formally $\rho_2^{(1),\star}$ is constant for the above given boundary conditions.
	Integrating the zero order terms of the density update 
	\begin{equation}
	\label{eq:AP_Acc_dens}
	\rho^{(0),\star} - \rho^{(0),n} + \Delta t \partial_x (\rho^{(0),\star} u^{(0),\star}) = 0
	\end{equation}
	on the computational domain, we obtain $\rho^{(0),\star} = \rho^{(0),n}$ and analogously for the first order terms $\rho^{(1),\star} = \rho^{(1),n}$. 
	Since $\alpha$ and $\chi$ do not change in the acoustic step, we obtain $\rho_2^{(0),\star} = \rho_2^{RS}$.
	Further we obtain with $(\alpha\rho_1)^{(1),\star} = \alpha^{(1),n} \rho_1^n$ that $\rho_2^{(1)\star} = 0$.  
	Regarding the velocity $u_2^{(0),\star}$, we obtain
	\begin{equation*}
		\partial_x u_2^{(0),\star} = \partial_x \left(u^{(0),\star} - \chi^{(0),n}(u_1^n - u_2^{(0),n})\right) = \mathcal{O}(\Delta t).
	\end{equation*}
	
	In the transport step, we find after some simple reformulations 
	\begin{align*}
		\alpha^{(0),\star\star} &= \alpha^{(0),n} + \Delta t \alpha^{(0),n}\partial_x( (\rho^{(0),\star} - \rho^{(0),n}) u^{(0),\star}) = \alpha^{(0),n},\\
		\alpha_2^{(0),\star\star} \rho_2^{(0),\star\star} &= \rho^{(0),\star} - (\alpha\rho_1)^{(0),\star\star} = \alpha_2^{(0),n}\rho_2^{RS},		
	\end{align*}
	from which follows $\rho_2^{(0),\star\star} = \rho_2^{RS}$. 
	Analogously to $u_2^{(0),\star}$, we also get $\partial_x u_2^{(0),\star\star} = \mathcal{O}(\Delta t)$.
	
	In the mixture step, the $\mathcal{O}(M^{-2})$ terms in \eqref{eq:Update_AlphaRho1} yield
	\begin{equation*}
		\partial_x\left( (\alpha\rho_1)^{\star\star\star} \left(1 - \frac{(\alpha\rho_1)^{\star\star\star}}{\rho^{(0),\star\star}}\right) \partial_x \rho_2^{(0),\star\star\star}\right) = 0.
	\end{equation*}
	Together with $\partial_x \rho_2^{(0),\star\star\star} = 0$ from the zero order terms of the relative velocity equation \eqref{eq:relVeloc_mix_update} it follows that $\rho_2^{(0),\star\star\star}$ is constant. 
	In addition holds 
	\begin{equation*}
		\alpha_2^{(0),n}\rho_2^{(0),\star\star\star} = \rho^{(0),\star} - \alpha^{(0),n}\rho_1^{\star\star\star} = \alpha^{(0),n} (\rho_1^\star - \rho_1^{\star\star\star}) + \alpha_2^{(0),n}\rho_2^{RS}.
	\end{equation*}
	Therefore we have $\rho_2^{(0),\star\star\star} = \rho_2^{RS} + \mathcal{O}(\Delta t)$.
	Analogously, we find from the $\mathcal{O}(M^{-1})$ terms in \eqref{eq:Update_AlphaRho1} and \eqref{eq:relVeloc_mix_update} that $\rho_2^{(1),\star\star\star}$ is constant. 
	Integration of \eqref{eq:alphaRho1_mix} on the computational domain leads to
	\begin{equation*}
			(\alpha_2 \rho_2)^{(1),\star\star\star}\alpha_2^{(0),n}\rho_2^{(1),\star\star\star} - \alpha_2^{(1),\star\star\star}\rho_2^{(0),\star\star\star} = \rho^{(1),\star} - \alpha^{(1),\star\star}\rho_1^{\star\star\star}
	\end{equation*} 
	and thus $\rho_2^{(1),\star\star\star} = \mathcal{O}(\Delta t)$.
	Regarding the pressure relaxation source term, we obtain immediately $p_2^{(0)} = 0 = p_2^{RS}$ for the $\mathcal{O}(M^{-2})$ terms and hence $\alpha^{(0),n+1} = \alpha^{(0),n}$ constant. 
	For the $\mathcal{O}(M^{-1})$ terms we obtain using the pressure expansion \eqref{eq:Exp_p} that $\rho_2^{(1),n+1} = 0$.
	Consequently, $\alpha^{(1),n+1} = \alpha^{(1),\star\star}$. 
	Since $\alpha^{(0),n+1}$ is constant, we find from the zero order terms $p_2^{(2)} = p_1$ which gives the pressure constraint in the well-prepared data. 
	
	From the density update and the final update of $\alpha \rho_1$, we obtain the final update of $\rho_2 \alpha_2$ and thus $\partial_x u^{(0),n+1} = \mathcal{O}(\Delta t)$. 
	Summarizing, the data $W^{n+1}$ preserves the asymptotics for $\Delta t \to 0$.
	We have proved the following theorem.
	\begin{theorem}[AP property compressible/weakly compressible]
        \label{theo:AP1}
		For well-prepared initial data $W^n \in \Omega_1$ and periodic or no-flux boundary conditions, the time semi discrete scheme described in Section \eqref{sec:TimeSemi} is asymptotic preserving and a consistent time discretization of the limit equations \eqref{eq:Limit1_cons} when $M \to 0$. 
	\end{theorem}
    
    \paragraph{Case 2} 
    We consider the case of two weakly compressible phases in the same Mach number regime with constant reference states $\rho_1^{(0),n} = \rho_1^{RS}$ and $\rho_2^{(0)} = \rho_2^{RS}$ and $\alpha^{(0),n}$ constant. 
    According to \eqref{eq:wp_M}, the well-prepared data at time $t^n$ is then given by
    \begin{align}
    \label{eq:AP_WP_Case2}
    \begin{split}
        \alpha^n(x)&= \alpha^{(0)} + \mathcal{O}(M), \quad \alpha^{(0)} \text{ const.},\\
        \rho_1^n(x) &= \rho_1^{RS} + M^2 \rho_1^{(2)} + \mathcal{O}(M^3), \quad \rho_2^n(x) = \rho_2^{RS} + M^2 \rho_2^{(2)} + \mathcal{O}(M^3), \\
        u_1^n(x) &= u_1^{(0)} + \mathcal{O}(M), \quad u_2^n(x) = u_2^{(0)} + \mathcal{O}(M),\quad \partial_x u_1^{(0)} = 0,  \quad \partial_x u_2^{(0)} = 0, \\
        p_1^{RS} &= p_2^{RS}, \quad p^{(2)} = p_1^{(2)} = p_2^{(2)}.
    \end{split}
    \end{align}
    Analogously to Case 1, we consider a general Mach number expansion for the state variables of the respective steps of the numerical scheme. 
    For the first step, we define
    \begin{equation*}
        W^\star(x) = W^{(0),\star} + M W^{(1),\star} + \mathcal{O}(M^2).
    \end{equation*}
    Using well-prepared data \eqref{eq:AP_WP_Case2}, we can rewrite the update for $\rho^{\star}$ as 
    \begin{equation*}
        \rho^\star - \Delta t^2 \partial_x^2 \left(\left(c^{RS}\left(\frac{\alpha^n \rho_1^n}{\rho^n}\right)\right)^2 \rho^\star\right) = \rho^n - \Delta t \partial_x (\rho u) - \Delta t^2 \partial_x^2 \left(\left(c^{RS}\left(\frac{\alpha^n \rho_1^{RS}}{\rho^n}\right)\right)^2 \rho^n\right),
    \end{equation*}
    where $c^{RS}$ is defined by \eqref{eq:Acc_cRS} and depends on $1/M$. 
    Inserting the Mach number expansions and sorting by terms of the Mach number, we find for the $\mathcal{O}(M^{-2})$ terms
    \begin{equation*}
        \partial_x^2 \rho^{(0),\star} = \partial_x^2 \rho^{(0),n} = 0.
    \end{equation*}
    This implies formally that $\rho^{(0),\star}$ is constant for periodic or no-flux boundary conditions. 
    Integrating the zero order density terms \eqref{eq:AP_Acc_dens} over the computational domain, we get $\rho^{(0),\star} = \rho^{(0),n}$. 
    Further, since we solve the acoustic step in the variables $(\rho, \alpha, \chi, \rho u, u_1 - u_2)$, we have $\alpha^\star = \alpha^n$ and therefore $\alpha^{(0),\star} = \alpha^{(0),n}$ constant, $\alpha^{(1),\star} = \alpha^{(1),n}$. 
    Since $\alpha^{(0)}$ constant and positive, it follows $\rho_1^{(0),\star} = \rho_1^{RS}$ and $\rho_2^{(0),\star} = \rho_2^{RS}$ which implies $\partial_x u^{(0),\star} = 0$.
    Since $\chi^{(0),\star}$ is constant it follows $\partial_x u_1^{(0),\star}, \partial_x u_2^{(0),\star} = 0$. 
    For the $\mathcal{O}(M^{-1})$ terms we find 
    \begin{equation*}
        \partial_x^2 \rho^{(1),\star} = \partial_x^2 \rho^{(1),n} \Leftrightarrow \partial_x^2 (\alpha^{(0),n}\rho_1^{(1),\star} + (1-\alpha^{(0),n})\rho_2^{(1),\star}) = 0.
    \end{equation*}
    Thus, formally $\rho_1^{(1),\star},\rho_2^{(1),\star}$ are constant.  
    Integrating the analogue equation to \eqref{eq:AP_Acc_dens} for the first order density perturbation $\rho^{(1),\star}$, we find $\alpha^{(0),n} \rho_1^{(1),\star} + \alpha_2^{(0),n} \rho_2^{(1),\star} = 0$ hence $\rho_1^{(1),\star}, \rho_2^{(1),\star} = 0$. 
    From the momentum update, we obtain then 
    \begin{equation*}
        u^{(0),\star} = u^{(0),n} -\Delta t \frac{ \partial_x (\alpha^{(0),n}\hat p_1^{(2),\star} + \alpha_2^{(0),n} \hat p_2^{(2),\star})}{\rho^{(0),n}}.
    \end{equation*}
    Hence $W^\star$ preserves the asymptotics for $\Delta t \to 0$. 
    We assume that the update after the transport step obeys $W^{\star\star} = W^{(0),\star\star} + M W^{(1),\star\star} + \mathcal{O}(M^2)$. 
    Then step \eqref{eq:Time_semi_transport} with $\partial_x u^{(0),\star} = 0$ for the $W^{(0),\star\star}$ terms leads to the following relations
    \begin{align*}
        \rho^{(0),\star\star} = \rho^{(0),\star} = \rho^{(0),n}, \\
        \alpha^{(0),\star\star} = \alpha^{(0),n}, \\
        \rho_1^{(0),\star\star} = \rho_1^{RS}, \\
        u^{(0),\star\star} = u^{(0),\star}, \\
        (u_1 - u_2)^{(0),\star\star} = (u_1 - u_2)^{(0),n}.
    \end{align*}
    In particular, this yields $\partial_x u_1^{(0),\star\star} = 0 = \partial_x u_2^{(0),\star\star}$. 
	Considering the $\mathcal{O}(M)$ terms we obtain 
	\begin{equation*}
		\alpha^{(1),\star\star} = \alpha^{(1),n} - \Delta t u^{(0),\star} \partial_x \alpha^{(1),n}
	\end{equation*}
	which is consistent with the continuous model. 
	Using this relation we obtain $$\rho_1^{(1),\star\star} = -\Delta t \rho_1^{RS} \partial_x \alpha^{(1),n} = \mathcal{O}(\Delta t).$$ 
	Therefore, $W^{\star\star}$ after the transport step preserves the asymptotics for $\Delta t \to 0$. 
	Assume that the update after the mixture step obeys $W^{\star\star\star} = W^{(0),\star\star\star} + M W^{(0),\star\star\star} + \mathcal{O}(M^2)$. 
	Inserting the Mach number expansions of $W^{\star\star\star}$ and $W^{\star\star}$ into the update \eqref{eq:Update_AlphaRho1}, we derive for the $\mathcal{O}(M^{-2})$ terms the following relation
	\begin{equation*}
		\partial_x \left(\rho_1^{(0),\star\star\star}(\rho^{(0),n} - \alpha^{(0),n}\rho_1^{(0),\star\star\star})\partial_x \rho_1^{(0),\star\star\star}\right) = 0.
	\end{equation*}
	Further, from the relative velocity update \eqref{eq:relVeloc_mix_update} it follows $\partial_x \rho_1^{(0),\star\star\star} = 0$. 
	Consequently, $\rho_1^{(0),\star\star\star}$ is constant for periodic or non-flux boundary conditions. 
	Integrating \eqref{eq:alphaRho1_mix} on the computational domain yields $(\alpha\rho_1)^{(0),\star\star\star} = \alpha^{(0),n}\rho_1^{RS}$ and thus $\rho_1^{(0),\star\star\star} = \rho_1^{RS}$. 
	Note that $\alpha$ does not change in this step and we have $\alpha^{(0),\star\star\star} = \alpha^{(0),n}$.
	An immediate consequence is $\rho_2^{(0),\star\star\star} = \rho_2^{RS}$ and $\partial_x (u_1 - u_2)^{(0),\star\star\star} = 0$. 
	For the $\mathcal{O}(M^{-1})$ terms we obtain then
	\begin{equation*}
		\partial_x^2 \left(\frac{(c_1^{RS})^2}{\rho_1^{RS}}\rho_1^{(1),\star\star\star} + \frac{(c_2^{RS})^2}{\rho_2^{RS} }\rho_2^{(1),\star\star\star}\right) = 0,
	\end{equation*}
	for which we formally obtain constant $\rho_1^{(1),\star\star\star}, \rho_2^{(1),\star\star\star}$.
	Considering the first order terms of \eqref{eq:alphaRho1_mix} we obtain that they are at least of order $\mathcal{O}(\Delta t)$. 
	For the higher order terms $\bar{p}$ we find that they are of the order $\bar{p} = \mathcal{O}(\Delta t M^2)$. 
	Therefore \eqref{eq:Mom_update_allspeed} yields for the velocity $u^{(0),\star\star\star} = u^{(0),\star\star} + \mathcal{O}(\Delta t^2)$ and hence its space derivative is zero. 
	Since $\rho_1^{(1),\star\star\star}, \rho_2^{(1),\star\star\star}$ are constant, the following relation at order $\mathcal{O}(1)$ is obtained from \eqref{eq:RelVeloc_update_allspeed}
	\begin{equation*}
		(u_1 - u_2)^{(0),\star\star\star} = (u_1 - u_2)^{(0),n} - \Delta t \partial_x (h_1^{(2),\star\star\star} - h_2^{(2),\star\star\star}).
	\end{equation*}
	Summarizing, $W^{\star\star\star}$ preserves the asymptotics for $M \to 0$. 
	Considering also the friction relaxation terms in the relative velocity \eqref{eq:Update_Friction}, we find the following time semi-discrete scheme for the limit equations 
	\begin{align}
		\label{eq:Time-semi-Mach-limit-eq}
	\begin{split}
		u^{(0),n+1} &= u^{(0),n} -\Delta t \frac{ \partial_x (\alpha^{(0),n}\hat p_1^{(2),\star} + \alpha_2^{(0),n} \hat p_2^{(2),\star})}{\rho^{(0),n}} \\
		(u_1 - u_2)^{(0),n+1} &= (u_1 - u_2)^{(0),n} - \Delta t \partial_x (h_1^{(2),\star\star\star} - h_2^{(2),\star\star\star}) - \zeta \chi^{(0),n+1}\chi_2^{(0),n+1} (u_1 - u_2)^{(0),n+1}
	\end{split}
	\end{align}
	which is consistent with the formulation of the limit equations \eqref{eq:Limit_M_same} in terms of $u_1, u_2$. 
	Note that $\partial_x p_k^{(2)} = \partial_x \hat{p}_k^{(2)}$ for $k=1,2$ and that adding the friction source term does not influence the zero order velocity derivatives $\partial_x u_1^{(0),n+1} = 0 = \partial_x u_2^{(0),n+1}$. 
	
	In the pressure relaxation \eqref{eq:Relax_Update_alpha}, we find for the $\mathcal{O}(M^{-2})$ terms 
	\begin{equation*}
		p_2(\rho_2^{(0),n+1}) = p_1(\rho_2^{(0),n+1}).
	\end{equation*}
	Hence $\alpha^{(0),n+1} = \alpha^{(0),n}$ since $\rho^{(0),n+1}$ and $\chi^{(0),n+1}$ are fulfilling the well-prepared constraint on the pressure $p_1^{(0)} = p_2^{(0)}$. 
	For the $\mathcal{O}(M^{-1})$ terms we find
	\begin{equation*}
	p_2(\rho_2^{(1),n+1}) = p_1(\rho_2^{(1),n+1}) \Leftrightarrow c_1^{(0)}\rho_1^{(1),n+1} = c_2^{(0)} \rho_2^{(1),n+1}
	\end{equation*}
	where $c_1^{(0),n+1}$ and $c_2^{(0),n+1}$ are known from the zero order terms. 
	Using the Mach number expansions for $\alpha^{(1),n+1} = \alpha^{(1),\star\star\star}$ leads to $\rho_1^{(1),n+1} = \rho_1^{(0),\star\star\star} = \mathcal{O}(\Delta t)$. 
	Analogously, one obtains $\rho_2^{(1),n+1} = \mathcal{O}(\Delta t)$. 
	Since $\alpha^{(0),n}$ and $\alpha^{(0),n+1}$ coincide, we obtain for the $\mathcal{O}(1)$ terms that 
	%$\alpha^{(2),n+1}$ is such that 
	\begin{equation*}
		p_2\left(\rho_2^{(2),n+1}\right) = p_1\left(\rho_2^{(2),n+1}\right) \Leftrightarrow 	\hat p_2^{(2),n+1} = \hat p_1^{(2),n+1}
	\end{equation*} 
	which means at time $t^{n+1}$ the pressure constraint on the well-prepared data is fulfilled. 
	Summarizing the above discussion, $W^{n+1}$ preserves the asymptotics for $\Delta t \to 0$ with the consistent time semi-discrete limit equations \eqref{eq:Time-semi-Mach-limit-eq} and we have proved the following theorem.
	\begin{theorem}[AP property for two weakly compressible phases]
    \label{theo:AP2}
		For well-prepared initial data $W^n \in \Omega_M^{wp}$ and periodic or no-flux boundary conditions, the time semi-discrete scheme described in Section \eqref{sec:TimeSemi} is asymptotic preserving and a consistent time discretization of the limit equations \eqref{eq:Limit2_const_dens} for $M \to 0$. 
	\end{theorem}
	The AP property for the case of two different low Mach number regimes with well-prepared data $W^n \in \Omega_2^{wp}$ can be obtained following the lines of the proof of the AP property discussed for Case 1 and 2.

	\section{Numerical results}
	\label{sec:NumRes}
	In this section, we illustrate by numerical experiments theoretical properties of the proposed RS-IMEX scheme \eqref{eq:Time_Semi_scheme}. 
	The first set of test cases concerns the homogeneous equations without pressure relaxation and friction terms. 
	The initial conditions, if not mentioned otherwise, are given in non-dimensional form. 
	They can be straightforwardly transformed into dimensional variables following the scaling procedure in Section \ref{sec:NonDim}.
	We compare the numerical results obtained with our scheme \eqref{eq:Time-semi-Mach-limit-eq} with a reference solution computed by a first order explicit Rusanov scheme. 
	The latter requires an acoustic time stepping given by the following CFL condition 
	\begin{equation*}
	\Delta t \leq \nu_{ac} \frac{\Delta x}{\max(| u_1^n \pm c_1^n/M_1|, |u_2^n \pm c_2^n/M_2|)}.
	\end{equation*}
	
	\subsection{Isentropic Euler equations}
	\label{sec:NumRes:Isentropic}
\begin{figure}[t!]
	 	\begin{subfigure}[c]{\textwidth}
	 	\includegraphics[scale=0.45]{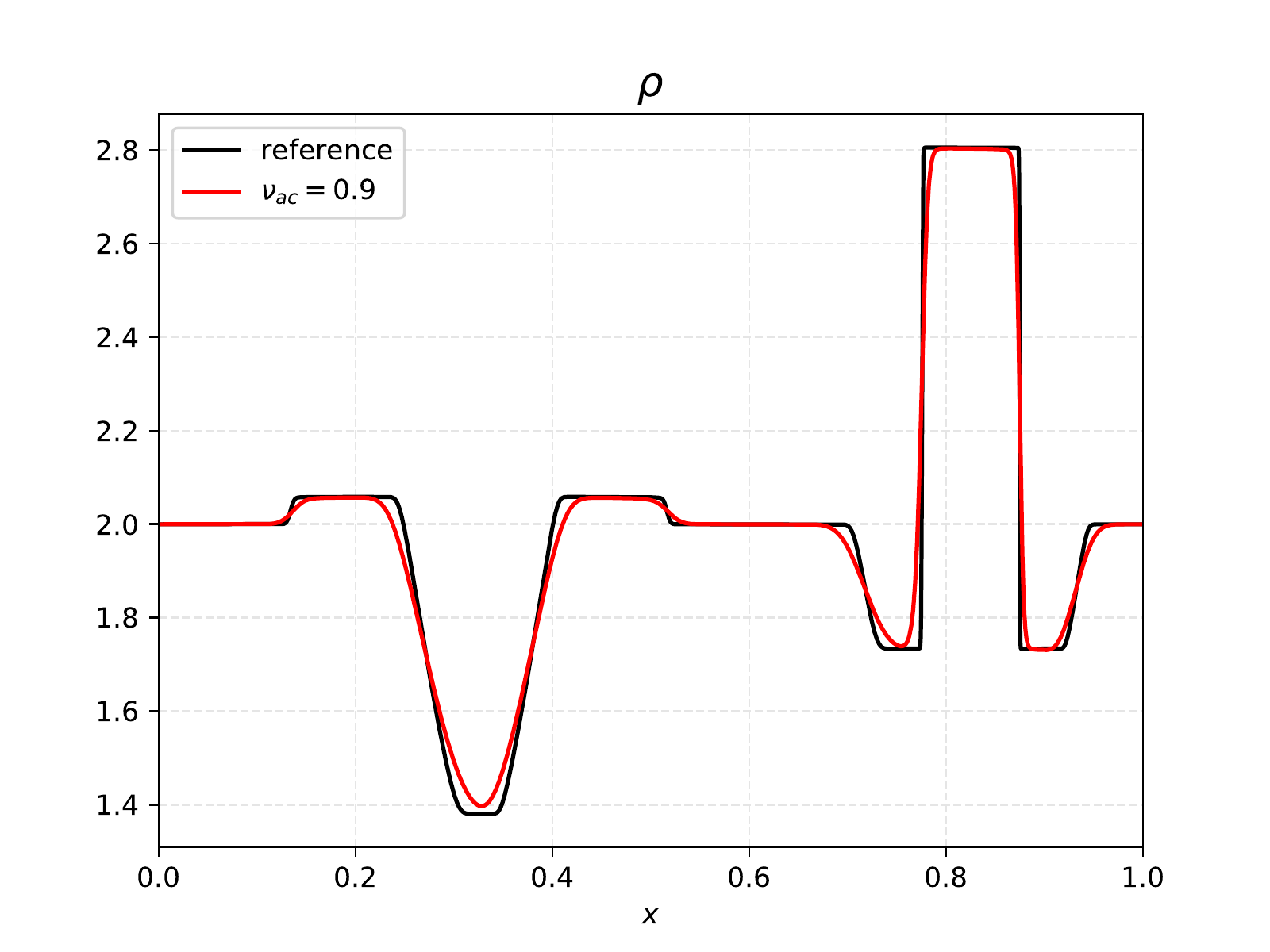}
	 	\includegraphics[scale=0.45]{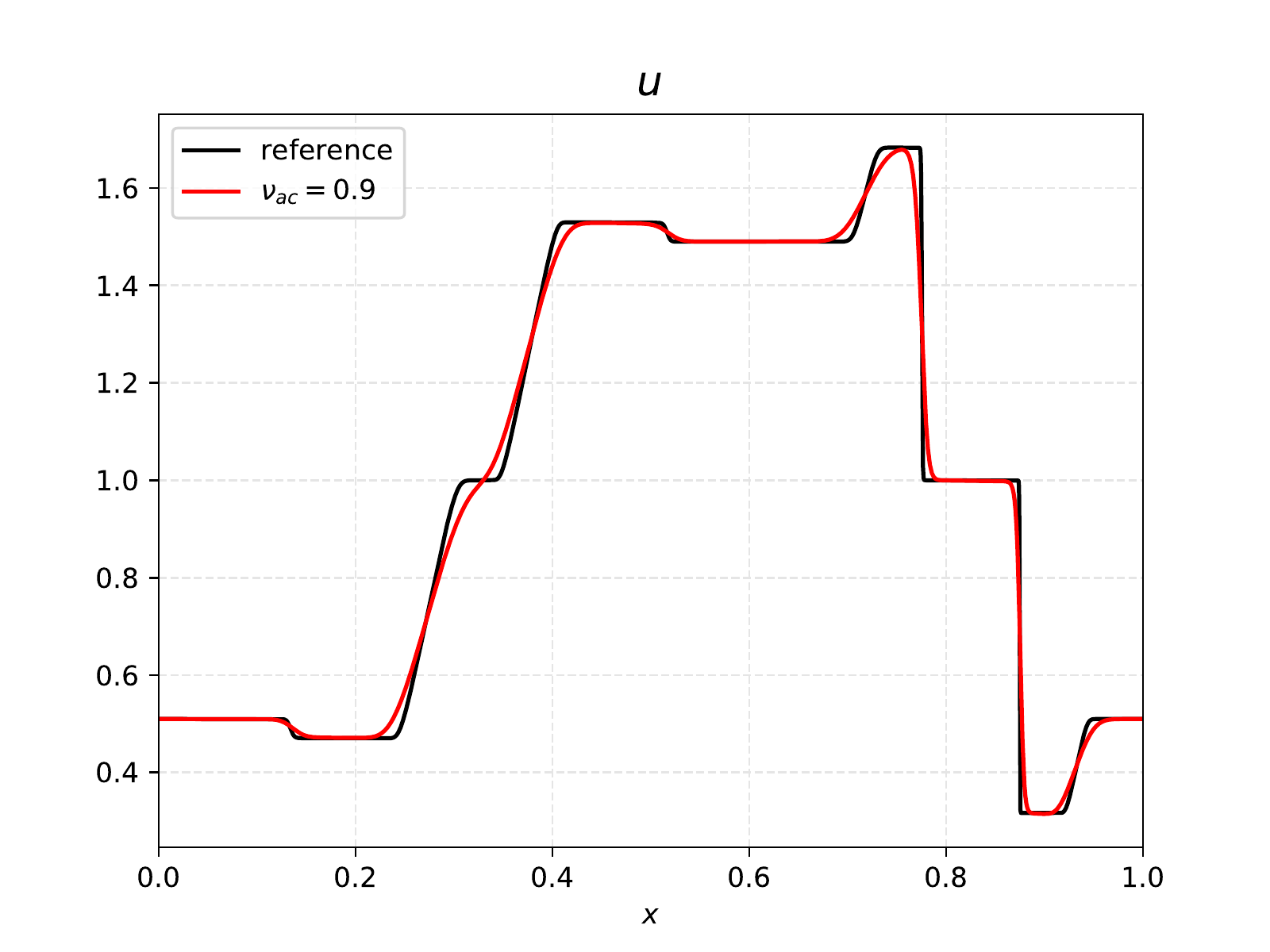}
	 	\subcaption{$M=0.99$}
	 	\label{fig:EulerDTComp}
	 	\end{subfigure}
 		\begin{subfigure}[c]{\textwidth}
	 	\includegraphics[scale=0.45]{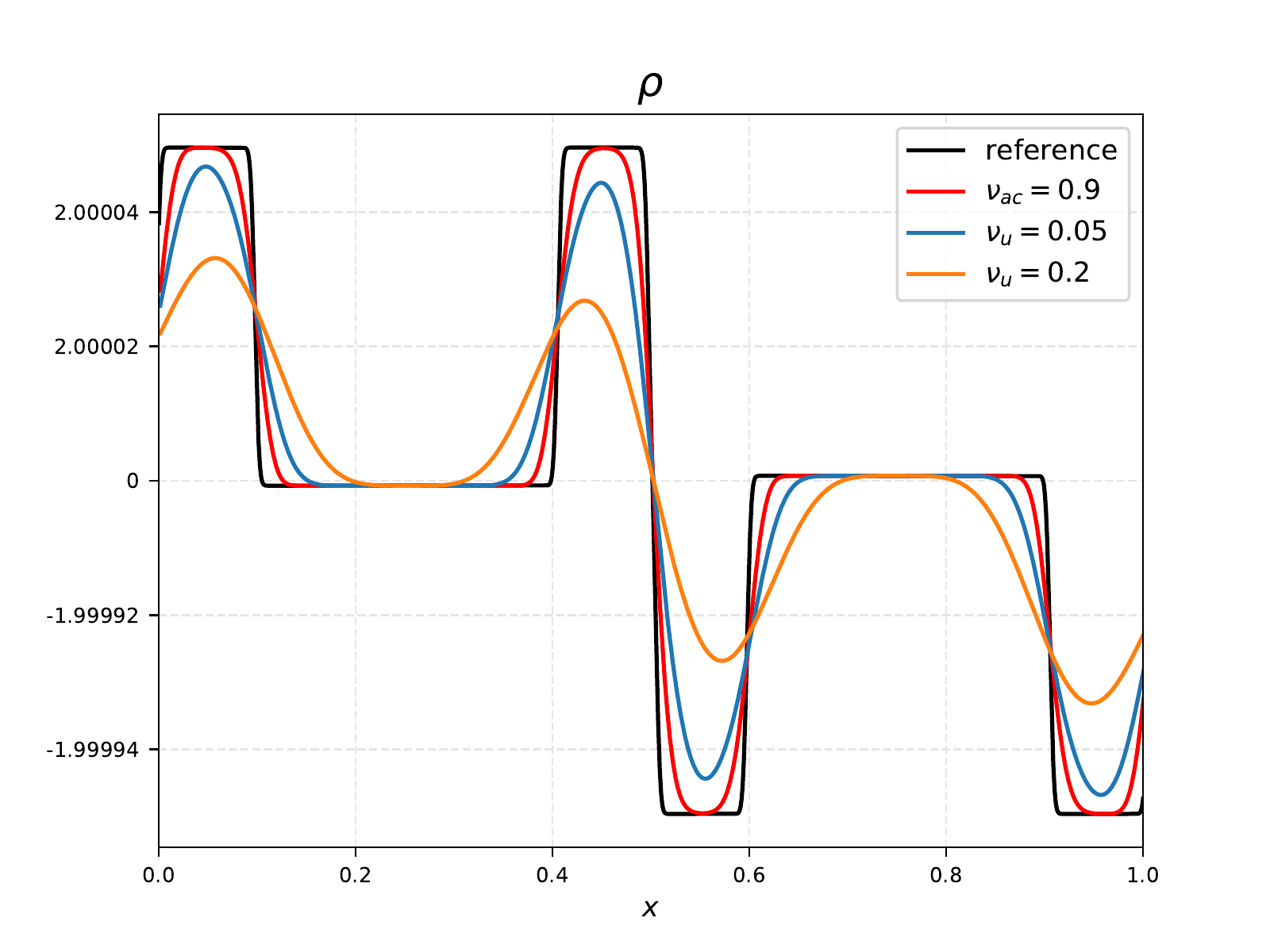}
	 	\includegraphics[scale=0.45]{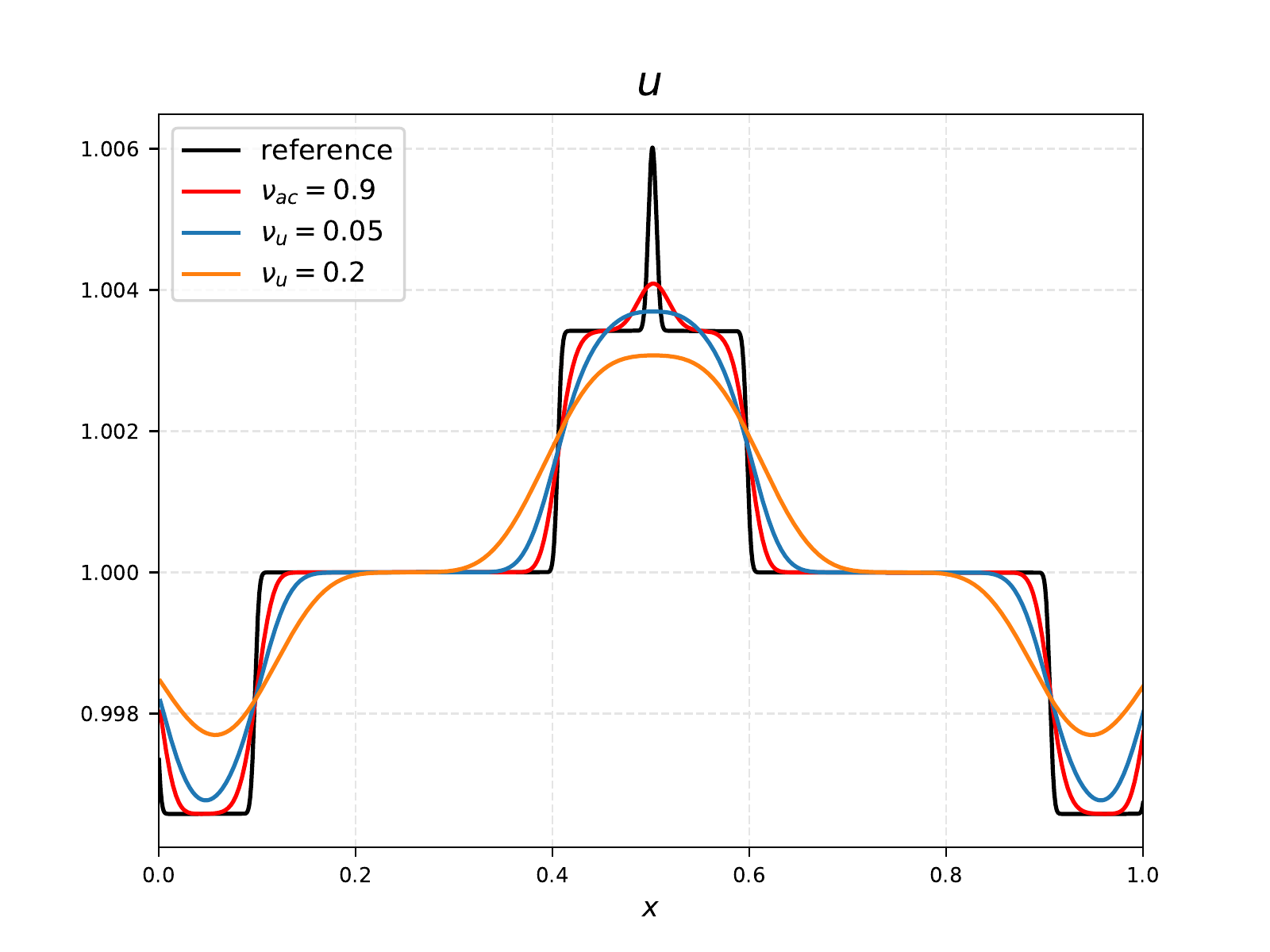}
	 	\subcaption{$M=10^{-2}$}
	 	\label{fig:EulerDTWeak}
	 	\end{subfigure}
	 	\caption{Riemann Problem from Section \ref{sec:NumRes:Isentropic}: Numerical results for density $\rho$ and velocity $u$ for $M=0.99$ (top panel) and $M = 10^{-2}$ (bottom panel) at final time $T_f = 0.075, 0.0015$, respectively, and grid size $\Delta x = 10^{-3}$.}
 		\label{fig:EulerDT}
	 \end{figure}

	With the following test we validate the consistency of the RS-IMEX all-speed scheme \eqref{eq:Time_Semi_scheme} with the isentropic Euler equations.
	Note that the model \eqref{eq:Comp_model} reduces to the isentropic Euler equations for single phase flow.
	We consider the low Mach number test case by Degond \& Tang described in \cite{DegTan2011,DimLouMicVig2018}. 
	We assign to each phase the ideal gas law with $\gamma = 1.4, \kappa = 1$ and set $\alpha = 0.5$. 
	The initial data are well-prepared and given by
	\begin{align}
	\label{eq:EulerDT}
	\begin{split}
		\rho(x,0) &= 
		\begin{cases}
			2  &\text{ for } x < 0.2, x \in (0.3,0.7], x > 0.8,\\
			2 + M^2 & \text{ for } x \in (0.2, 0.3],\\
			2 - M^2 & \text{ for } x \in (0.7,0.8],\\
		\end{cases} \\
		u(x,0) &= 
		\begin{cases}
			1 - M^2/2 & \text{ for } x < 0.2, x > 0.8,\\
			1, & \text{ for } x \in (0.2,0.3] \cup (0.7,0.8],\\
			1 + M^2/2 &\text{ for } x \in (0.3,0.7].
		\end{cases}
	\end{split}
	\end{align}  
	The solution consists of several Riemann problems leading to shocks and contact discontinuities that are stronger the larger the Mach number is chosen. 
	The computational domain is given by $[0,1]$ and is discretized using $\Delta x = 10^{-3}$.
	We first consider the compressible regime with $M = 0.99$ with $\nu_{ac} = 0.9$ leading to $\Delta t = 1.6\cdot 10^{-4}$.
	The results are given in Figure \ref{fig:EulerDTComp} where we see that our RS-IMEX all-speed scheme captures all shock positions correctly. 
	
	Next, we consider a weakly compressible regime with $M=10^{-2}$. 
	For the RS-IMEX scheme, we consider different time step sizes given by $\nu_{ac} = 0.9$ and $\nu_u = 0.05, 0.2$ which corresponds to $\Delta t = 3.6\cdot 10^{-6}$ (417 time steps), $\Delta t = ~2.5\cdot 10^{-5}$ (60 time steps) and $\Delta t = 10^{-4}$ (15 time steps), respectively. 
	The reference solution computed with the explicit Rusanov scheme thus $\nu_{ac} = 0.9$ and $\Delta t = 3.6\cdot 10^{-6}$ (417 time steps).
	The results are given in Figure \ref{fig:EulerDTWeak}. 
	We can clearly observe that the scheme is able to capture accurately all shocks and rarefactions for an acoustic time step, whereas, as expected, the acoustic waves are more diffused the larger the time step is chosen. 
	Note, that the wave fan of the isentropic Euler equations consists only of acoustic waves which are smeared for large time steps in our RS-IMEX scheme.  
	
	\subsection{Accuracy}
	\label{sec:NumRes:Acc}
	The following numerical test concerns the accuracy of the RS-IMEX scheme.
	We recall that by construction, the scheme is formally first order. 
	However due to the operator splitting in several sub-systems, we want to validate the experimental order of convergence (EOC). 
	We consider a double rarefaction test case with smooth initial data given by
	\begin{align*}
	\rho_1(x,0) &= 1, \quad \rho_2(x,0) = 1, \quad \alpha = 0.99 \\
	u_1(x,0) &= 
	\begin{cases}
	-2 & \text{ for } x < 0.4 \\
	30 -16 x + 200 x^2 & \text{ for } x \in [0.4, 0.5) \\
	-70 + 240 x - 200 x^2 & \text{ for } x \in [0.5, 0.6) \\
	2 &\text{ for } x \geq 0.6 \\
	\end{cases}, \\
	u_2(x,0) &= u_1(x,0).\\
	\end{align*}
	We assign to phase one the ideal gas law with $\gamma_1 = 1.4$ and $\kappa_1 = 1$ and to phase two the stiffened gas equation with $\gamma_2 = 2.8, ~p_\infty = 1, ~\kappa_2 = 2$. 
	The reference solution was computed with the Rusanov scheme on a fine grid with $N = 2^{15}$ grid cells on the domain $[0,1]$ up to a final time $T_f = 0.1 M_2$ using $\nu_{ac} = 0.9$. 
	In Table \ref{tab:ErrorConv} the $L^1$ error and EOC with respect to that reference solution is displayed.
	The EOC is of first order for all considered Mach number regimes. 
	In addition, the magnitude of the error in the density $\rho$ as well as in $\alpha\rho_1$ is of order of $\max(M_1,M_2)^2$. 
	This confirms the asymptotic preserving property of the RS-IMEX scheme \eqref{eq:Time_Semi_scheme}.
	The initial jump in the velocities $u_1, u_2$ triggers a perturbation of the densities $\rho_1, \rho_2$ of order $\mathcal{O}(M_1^2), \mathcal{O}(M_2^2)$ respectively. 
	Analogous results have been obtained for different initial values for $\alpha$ which we do not present here. 
	
	\begin{table}[h!]
	\renewcommand{\arraystretch}{1.25}
		\begin{tabular}{c c c c c c c c c c}
		$M_1, M_2$ & $N$& $\rho$& & $\alpha_1 \rho_1$ & & $\rho u$ & & $u_1 - u_2$ & \\\hline\hline
		\multirow{7}{*}{$10^{-1}, 10^{-1}$}&$2^6$ & 4.235e-03 &---& 3.650e-03 &---& 5.500e-02 &---&
		1.008e-01&---\\ 
		&$2^7$ & 2.089e-03 &1.02&1.844e-03 &0.98& 2.708e-02 &1.02&
		5.832e-02 &0.79\\ 
		&$2^8$ &1.057e-03 & 0.98& 9.541e-04 &0.95 & 1.401e-02 &0.95&
		3.319e-02 &0.81\\ 
		&$2^9$ &5.476e-04& 0.94& 4.999e-04 &0.93&7.389e-03 &0.92&
		1.820e-02 &0.86\\ 
		&$2^{10}$ & 2.514e-04 &1.12& 2.308e-04 &1.11&3.379e-03 &1.12&
		9.320e-03&0.96\\ 
		&$2^{11}$ &1.209e-04  &1.05& 1.111e-04 &1.06& 1.647e-03 &1.03&
		4.629e-03 &1.01\\ 
		&$2^{12}$ & 5.921e-05 &1.03& 5.456e-05 &1.02& 8.156e-04 &1.02&
		2.176e-03&1.09\\\hline
		\multirow{7}{*}{$10^{-3},10^{-3}$}&$2^6$ &3.480e-05 &---& 3.502e-05 &---& 4.182e-02 &---&
		1.324e-01&---\\ 
		&$2^7$ & 1.672e-05 &1.06&1.673e-05 &1.06& 1.983e-02 &1.07&
		7.081e-02 &0.90\\ 
		&$2^8$ &8.512e-06 & 0.97&8.523e-06 &0.97 & 1.007e-02 &0.97&
		3.843e-02 &0.88\\ 
		&$2^9$ &4.465e-06& 0.93& 4.468e-06 &0.93&5.298e-03 &0.92&
		2.056e-02 &0.90\\
		&$2^{10}$ & 2.000e-06 &1.16& 2.005e-06 &1.15&2.367e-03 &1.16&
		9.961e-03&1.05\\ 
		&$2^{11}$ &9.512e-07  &1.07& 9.537e-07 &1.07& 1.128e-03 &1.07&
		4.936e-03 &1.01\\ 
		&$2^{12}$ & 4.607e-07 &1.05& 4.610e-07 &1.05& 5.513e-04 &1.03&
		2.401e-03&1.03\\ \hline
		\multirow{7}{*}{$10^{-1},10^{-2}$}&$2^6$ & 1.054e-03 &---& 1.049e-03 &---& 3.167e-02 &---&
		7.785e-02&---\\ 
		&$2^7$ & 5.261e-04 &1.00&5.243e-04 &1.00& 1.581e-02 &1.00&
		3.977e-02 &0.97\\ 
		&$2^8$ &2.611e-04 & 1.01& 2.604e-04 &1.01 & 7.865e-03 &1.01&
		2.053e-02 &0.95\\ 
		&$2^9$ &1.290e-04& 1.02& 1.286e-04 &1.02&3.891e-03 &1.02&
		9.990e-03 &1.03\\ 
		&$2^{10}$ & 6.277e-05 &1.04&6.266e-05 &1.04&1.907e-03 &1.05&
		4.961e-03&1.01\\ 
		&$2^{11}$ &3.072e-05  &1.03& 3.076e-05 &1.03& 9.189e-03 &1.09&
		2.469e-03 &1.01\\ 
		&$2^{12}$ & 1.583e-05 &0.96&1.589e-05 &0.95& 4.314e-03 &1.02&
		1.082e-03&1.19\\\hline 
		\multirow{7}{*}{$10^{-2},10^{-3}$}&$2^6$ &1.083e-04 &---& 1.084e-04 &---& 3.130e-02 &---&
		7.721e-02&---\\ 
		&$2^7$ & 5.466e-05 &0.99&5.463e-05 &0.99& 1.564e-02 &1.00&
		3.905e-02 &0.98\\ 
		&$2^8$ &2.713e-05 & 1.01&2.710e-05 &1.01 & 7.785e-03 &1.01&
		2.007e-02 &0.96\\ 
		&$2^9$ &1.341e-05& 1.02& 1.338e-05 &1.02&3.864e-03 &1.01&
		9.749e-03 &1.04\\
		&$2^{10}$ & 6.662e-06 &1.01& 6.650e-06 &1.01&1.900e-03 &1.02&
		4.833e-03&1.01\\ 
		&$2^{11}$ &3.435e-06  &0.96& 3.431e-06 &0.95& 9.228e-04 &1.04&
		2.404e-03 &1.01\\ 
		&$2^{12}$ & 1.901e-06 &0.85& 1.899e-06 &0.85& 4.407e-04 &1.07&
		1.051e-03&1.19\\\hline\hline 
	\end{tabular}
	\caption{Riemann Problem from Section \ref{sec:NumRes:Acc}: $L^1$ error and convergence rates for the smooth rarefaction test case with $\nu_{ac} = 9$ and final time $T_f = 0.2 M_2$ on the computational domain $[0,1]$.}
	\label{tab:ErrorConv}
	\end{table} 

	\subsection{Riemann Problem}
	\label{sec:NumRes:RP}
	The basis of the next series of numerical results is a Riemann problem consisting of a jump in the phase densities $\rho_1, \rho_2$ analogously to the Riemann problem for the isentropic Euler equations \eqref{eq:EulerDT}. 
	For all tests we choose the ideal gas law with $\gamma_1=1.4, \kappa_1=1$ for phase one and the stiffened gas equation with $\gamma = 2.8, \kappa_2=2, p_\infty=1$ for phase two. 
	The initial data for phase densities and velocities are given by
	\begin{align}
	\label{eq:RP}
	\begin{split}
		\rho_1(x,0) &= 
		\begin{cases}
			1 + M_1^2\rho_1^{(2)} & \text{ for } x < 0.5 \\
			1 &\text{ for } x \geq 0.5 \\
		\end{cases}, \quad 
		\rho_2(x,0) = 
		\begin{cases}
		1 + M_2^2\rho_2^{(2)} & \text{ for } x < 0.5 \\
		1 &\text{ for } x \geq 0.5 \\
		\end{cases} \\
		u_1(x,0) &= 0.25 = u_2(x,0).
	\end{split}
	\end{align}
	The initial values for $\rho_1^{(2)}$ and $\rho_2^{(2)}$ will be specified below.
	
\begin{figure}[t!]
 		\begin{subfigure}[c]{\textwidth}
 			\includegraphics[scale=0.33]{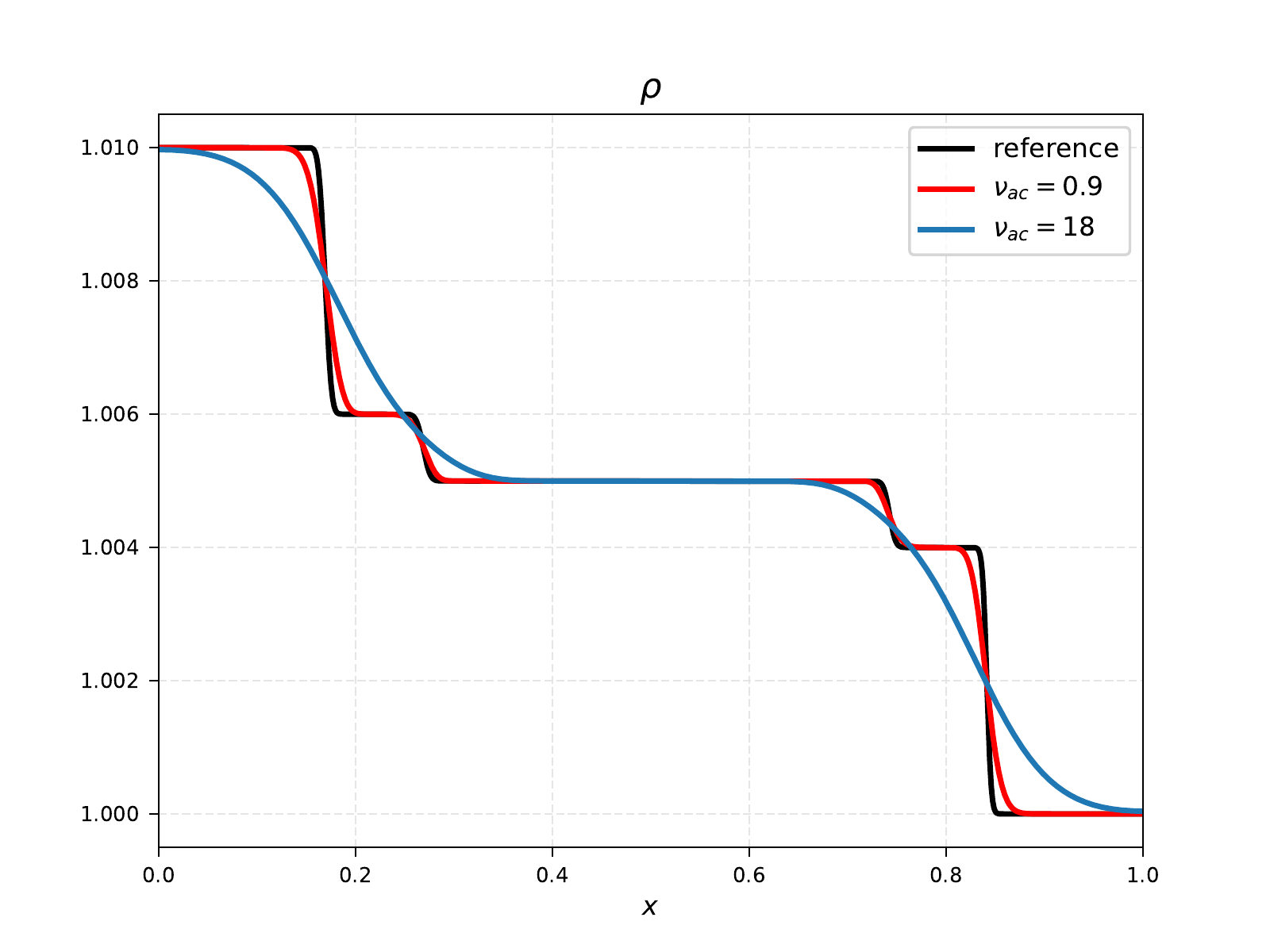}
 			\includegraphics[scale=0.33]{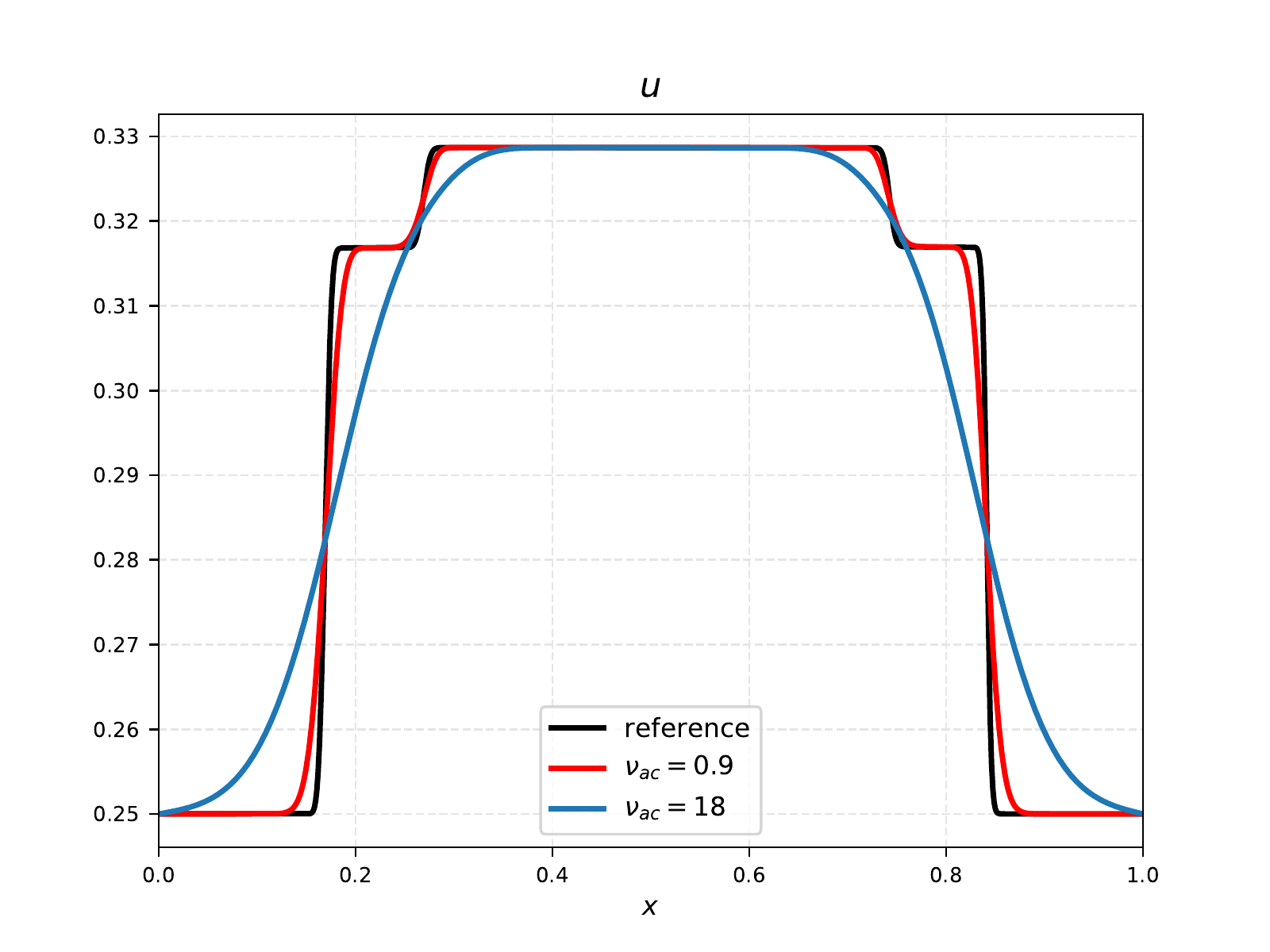}
 			\includegraphics[scale=0.33]{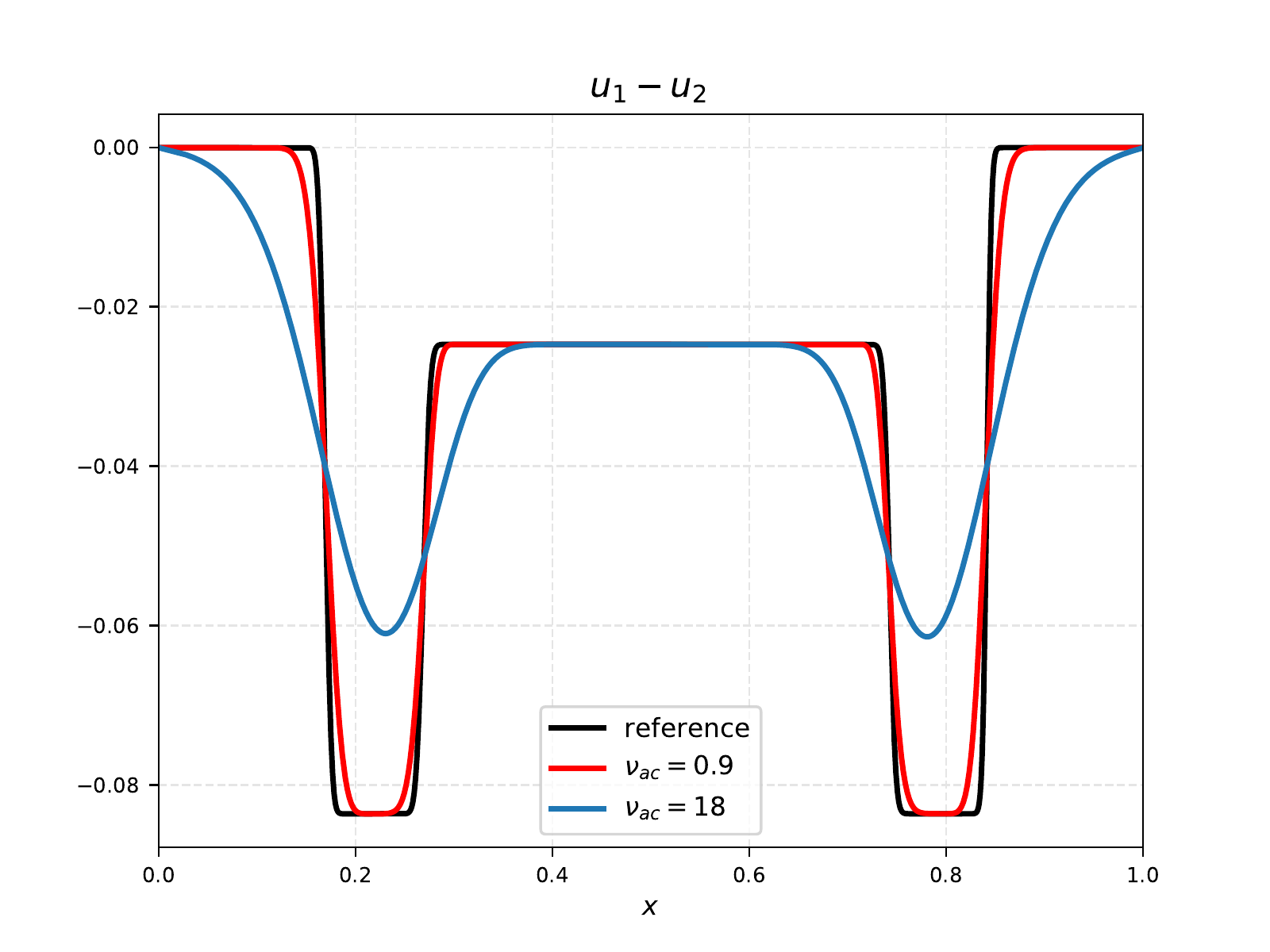}
 			\subcaption{$M=10^{-1}$}
 		\end{subfigure}
 		\begin{subfigure}[c]{\textwidth}
 			\includegraphics[scale=0.33]{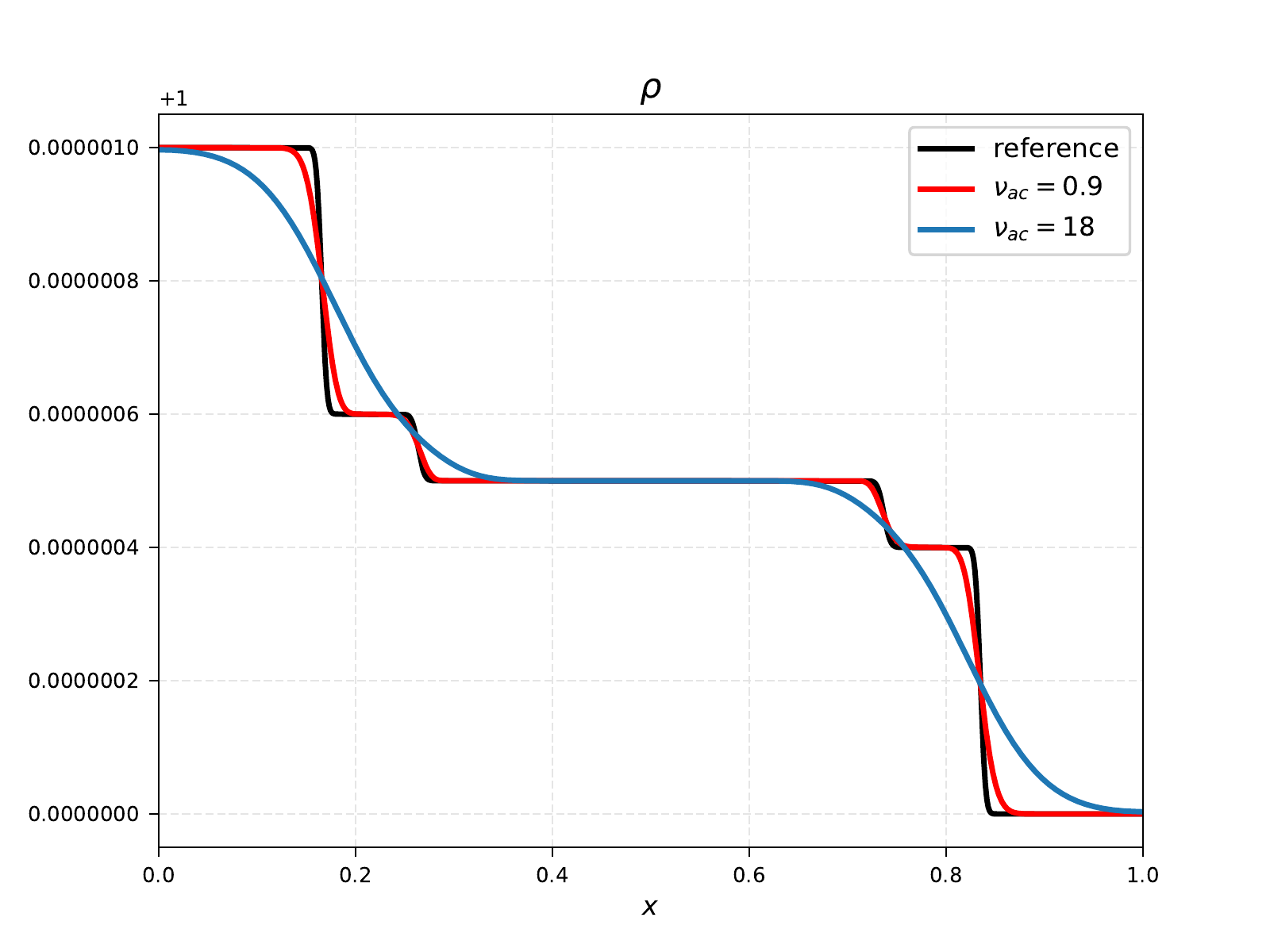}
 			\includegraphics[scale=0.33]{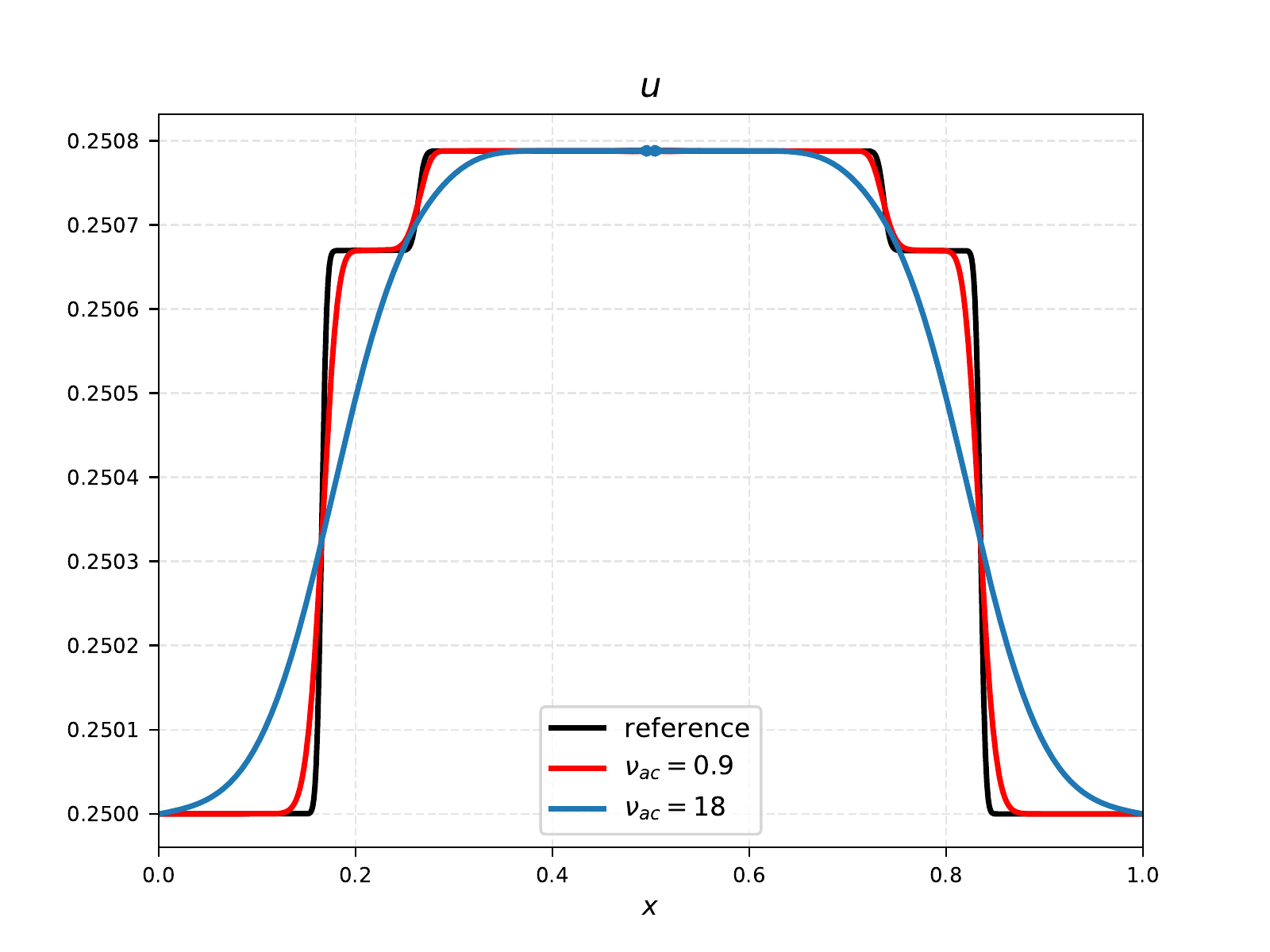}
 			\includegraphics[scale=0.33]{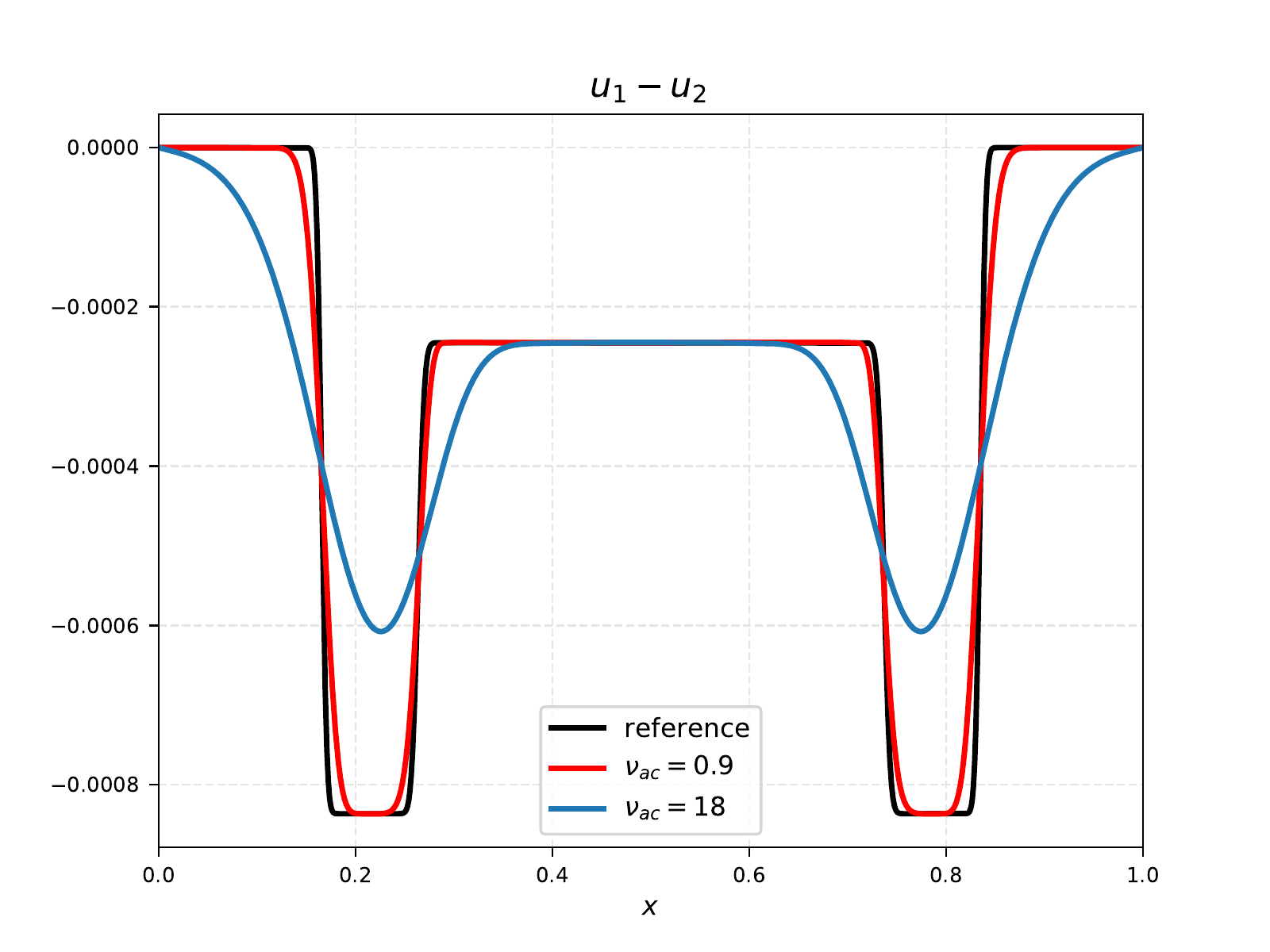}
 			\subcaption{$M=10^{-3}$}
 		\end{subfigure}
 			\caption{Riemann Problem from Section \ref{sec:NumRes:RP}: Numerical results for mixture density $\rho$, mixture velocity $u$ and relative velocity $u_1 - u_2$ for $M=10^{-1}$ (top panel) and $M = 10^{-3}$ (bottom panel) at final time $T_f = 0.2 M$, respectively, and grid size $\Delta x = 10^{-3}$.}
 		\label{fig:RPMSame}
 	\end{figure}

	\textbf{Same Mach number regime.}
	In the first series of numerical tests we consider two phases in the same flow regime with $M_1 = M = M_2$ and $\rho_1^{(2)} = 1 = \rho_2^{(2)}$.
	This leads to well-prepared initial data with $p_1^{RS} = 1 = p_2^{RS}$.  
	We start with a homogenous mixture governed by the homogeneous model with constant volume fraction $\alpha$.
	Afterwards we consider a jump in $\alpha$ and a smooth transition modelling a sharp and a diffusive interface.  
	The computational domain is $[0,1]$ with a uniform mesh size of $\Delta x = 10^{-3}$.
	The final time $T_f = 0.2 M$ is chosen such that all waves are still contained in the computational domain to avoid boundary effects. 
	In Figure \ref{fig:RPMSame} the numerical results are displayed for $M=10^{-1}$ in the top panel and $M=10^{-3}$ in the bottom panel. 
	The computations are done with an acoustic time step with $\nu_{ac} = 0.9$ and a larger time step with $\nu_{ac}=18$ which corresponds to $\Delta t = 2.9\cdot 10^{-5}, ~2.9\cdot 10^{-4}$ for $M=10^{-1}$ and $\Delta t = 2.9\cdot 10^{-7}, ~2.9\cdot 10^{-6}$ for $M=10^{-3}$, respectively. 
	We can see that the acoustic waves are captured accurately by the RS-IMEX scheme with a time step oriented to the acoustic waves, whereas they are faded out for larger time steps. 
	Note that since $\alpha$ is constant the material wave at $x=0.5$ is not visible. 
	
\begin{figure}[t!]
 	 	\begin{center}
 		\includegraphics[scale=0.4]{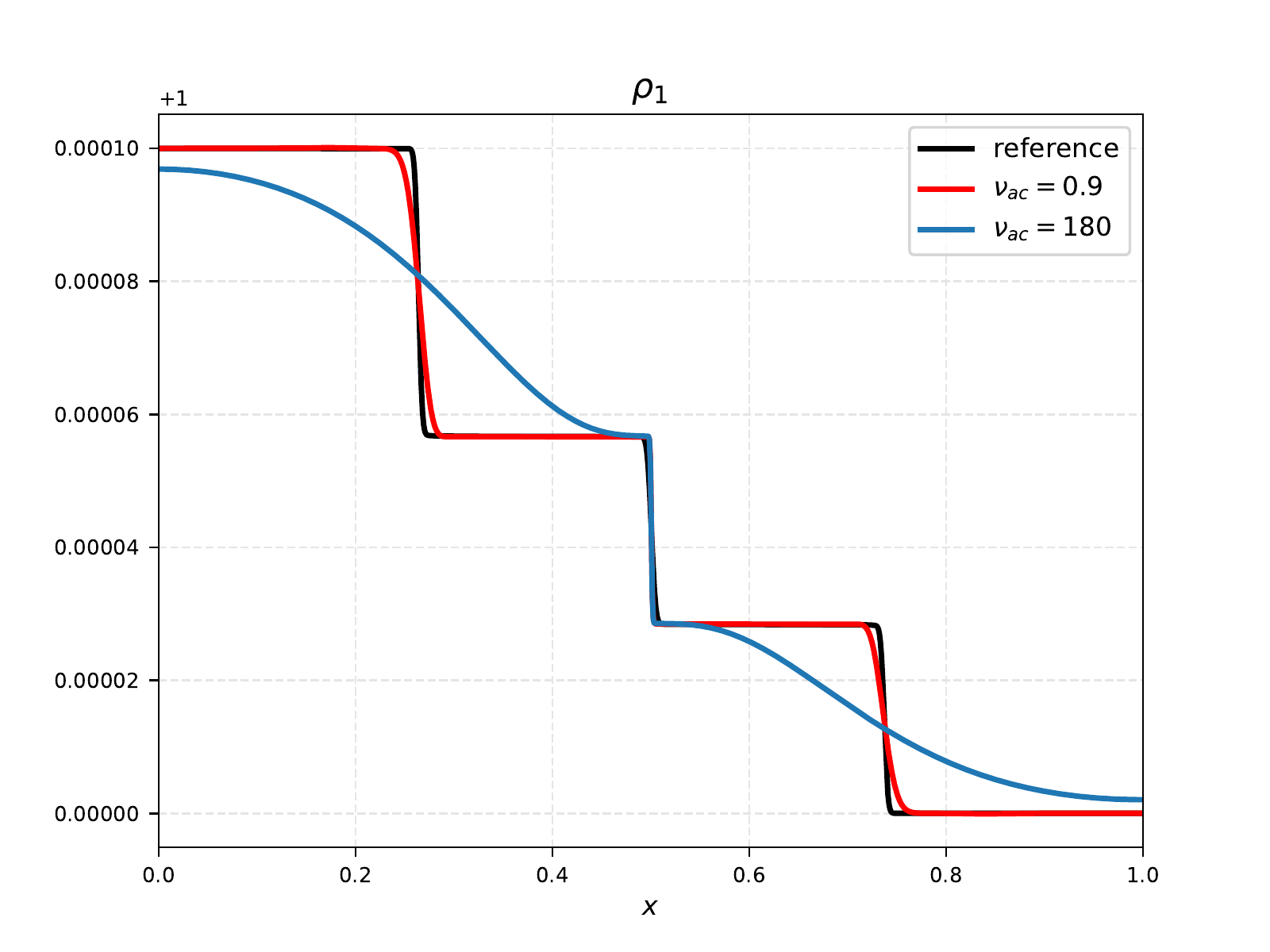}
 		\includegraphics[scale=0.4]{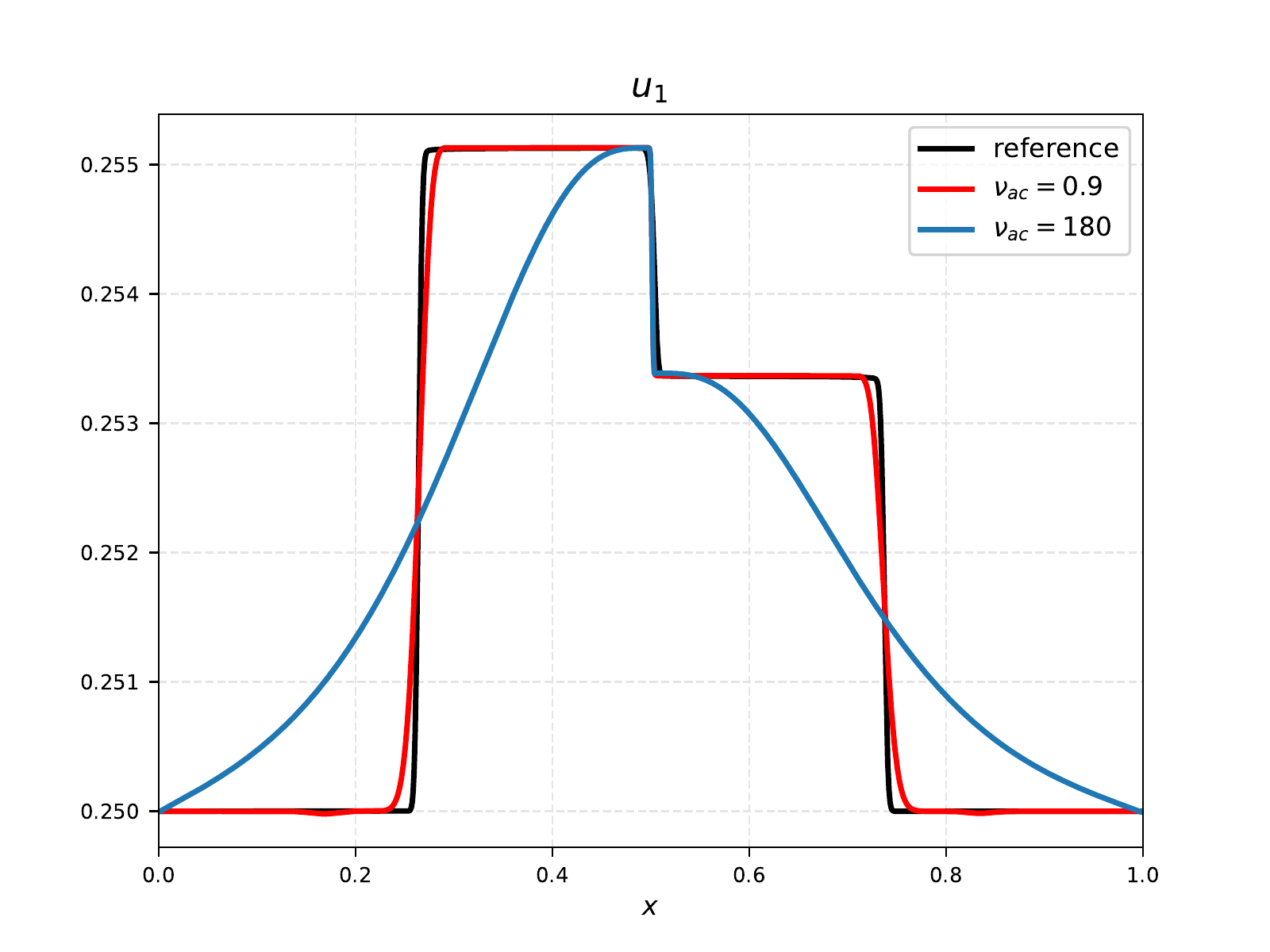}\\
 		\includegraphics[scale=0.4]{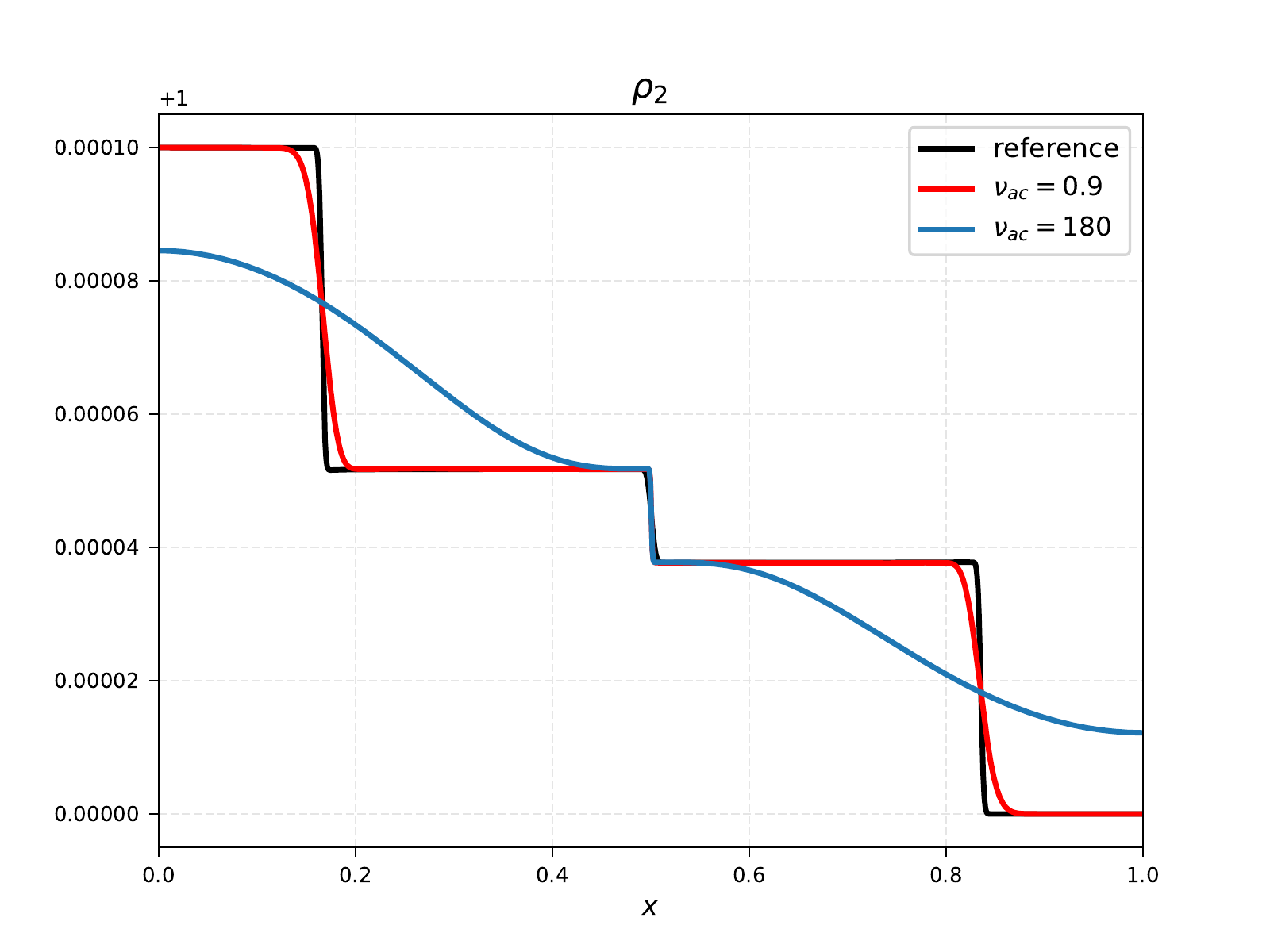}
 		\includegraphics[scale=0.4]{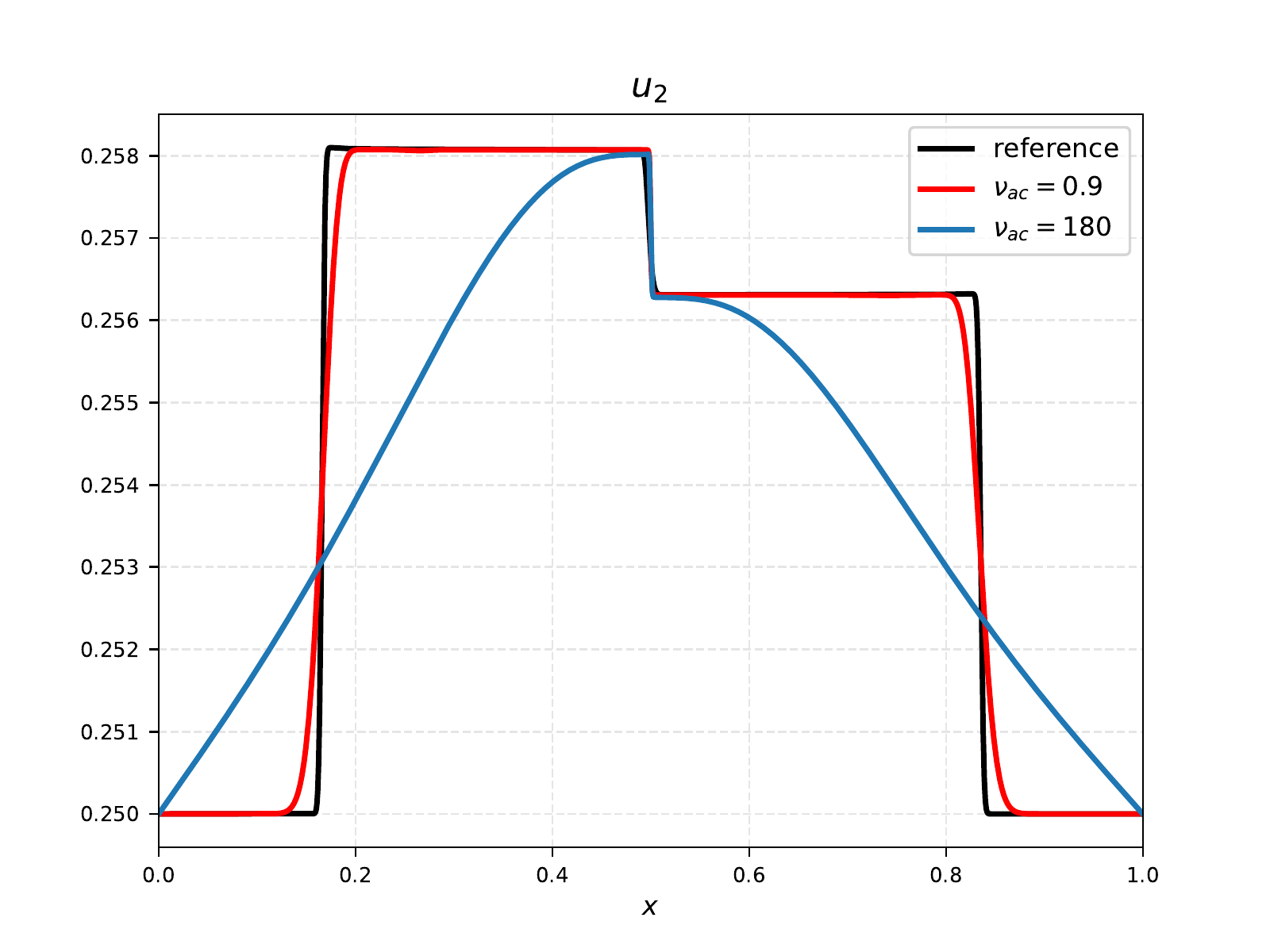}\\
 		\includegraphics[scale=0.4]{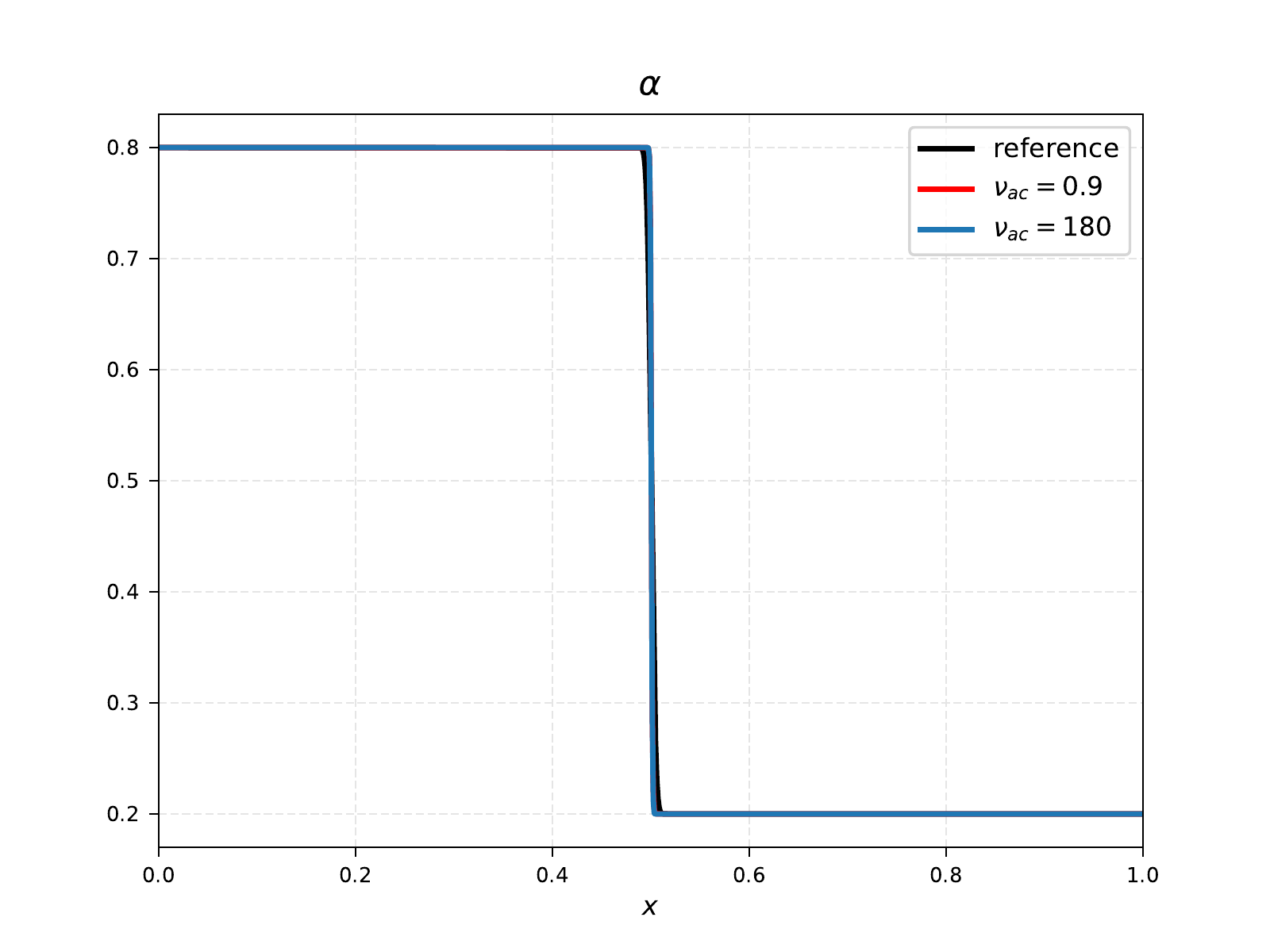}
 		\includegraphics[scale=0.4]{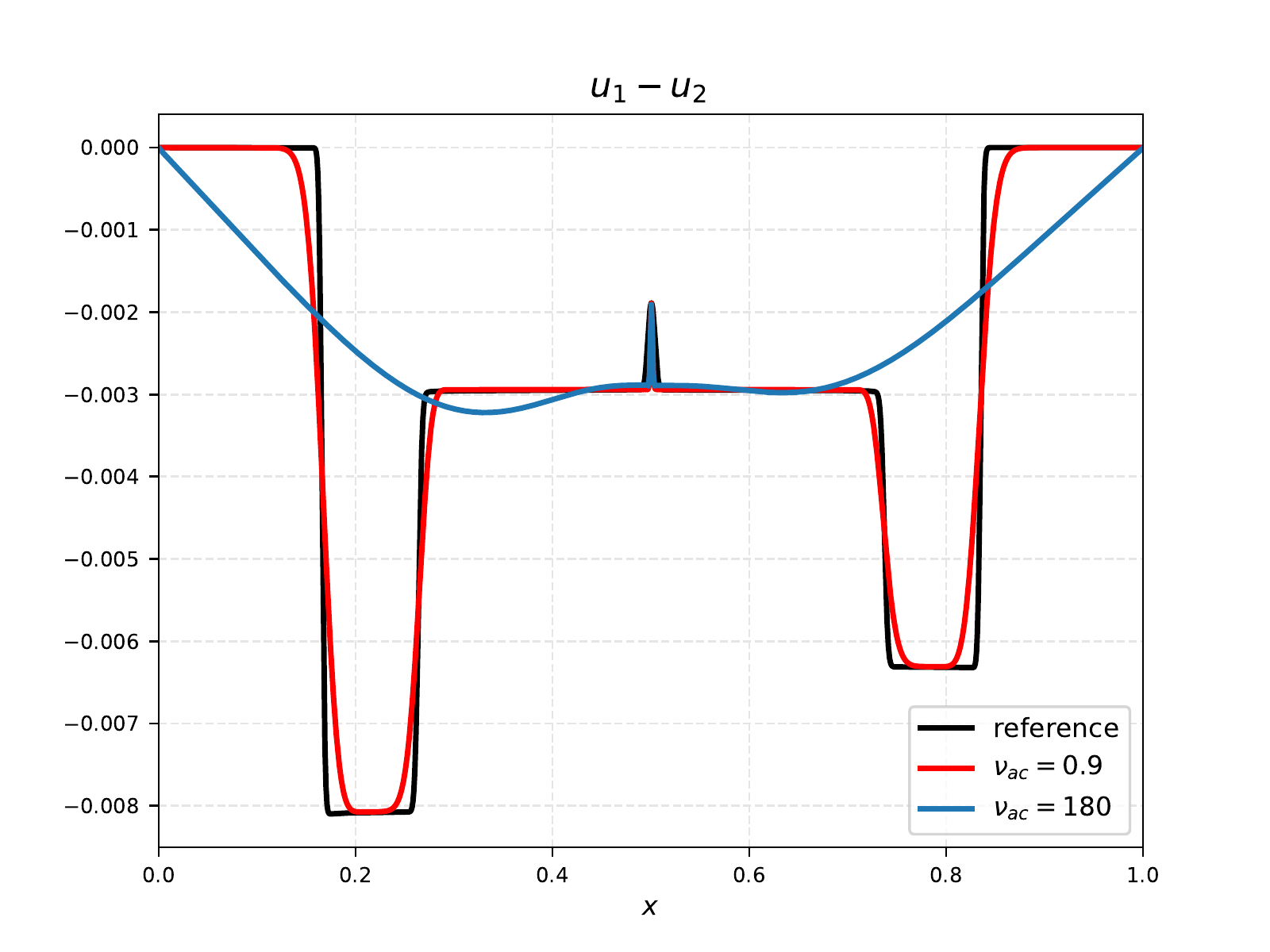}
 		\end{center}
 		\caption{Riemann Problem from Section \ref{sec:NumRes:RP} with an initial jump in $\alpha$: Numerical results for mixture density $\rho$, mixture velocity $u$ and relative velocity $u_1 - u_2$ for $M=10^{-2}$ at final time $T_f = 0.2 M$ with grid size $\Delta x = 10^{-3}$.}
 		\label{fig:RPMSameJumpA}
 	\end{figure}

	To capture the material wave, we consider a jump in the volume fraction at $x = 0.5$ given by $\alpha_L = 0.8$ and $\alpha_R = 0.2$ which is transported with $u$. 
	In Figure \ref{fig:RPMSameJumpA} the numerical results for $M=10^{-2}$ are presented. 
	The simulation is done with an acoustic time step with $\nu_{ac} = 0.9$ and a material time step resulting in $\nu_{ac} = 180$ leading to $\Delta t = 2.9 \cdot 10^{-6}$ and $\Delta t = 5.8 \cdot 10^{-4}$ respectively. 
	We want to stress that the material wave at $x=0.5$ is captured also for large time steps as sharp as the reference solution calculated with $N=30000$ cells although using a time step that is 200 times larger.
	
\begin{figure}[t!]
	 	\begin{center}
	 		\includegraphics[scale=0.4]{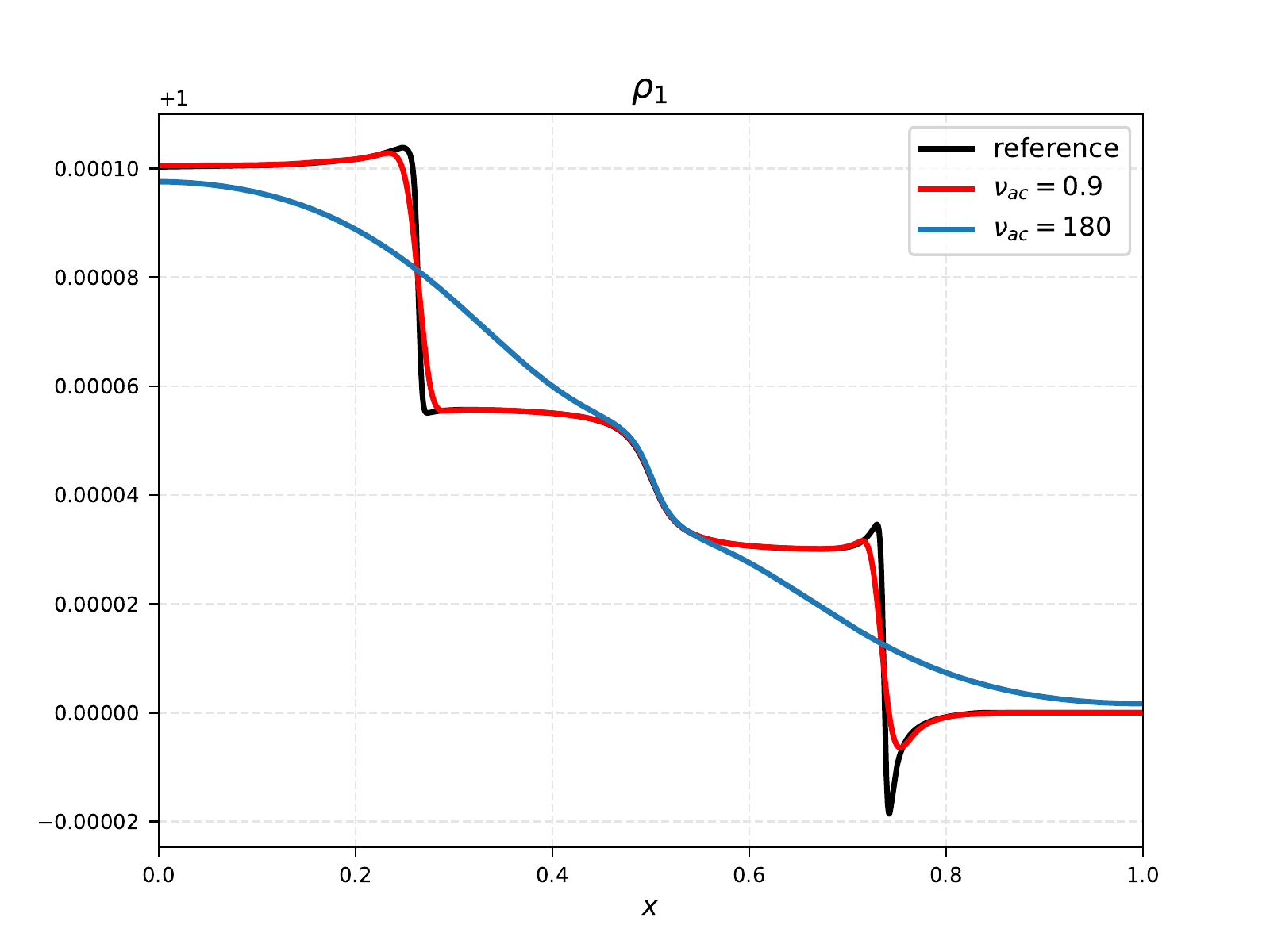}
	 		\includegraphics[scale=0.4]{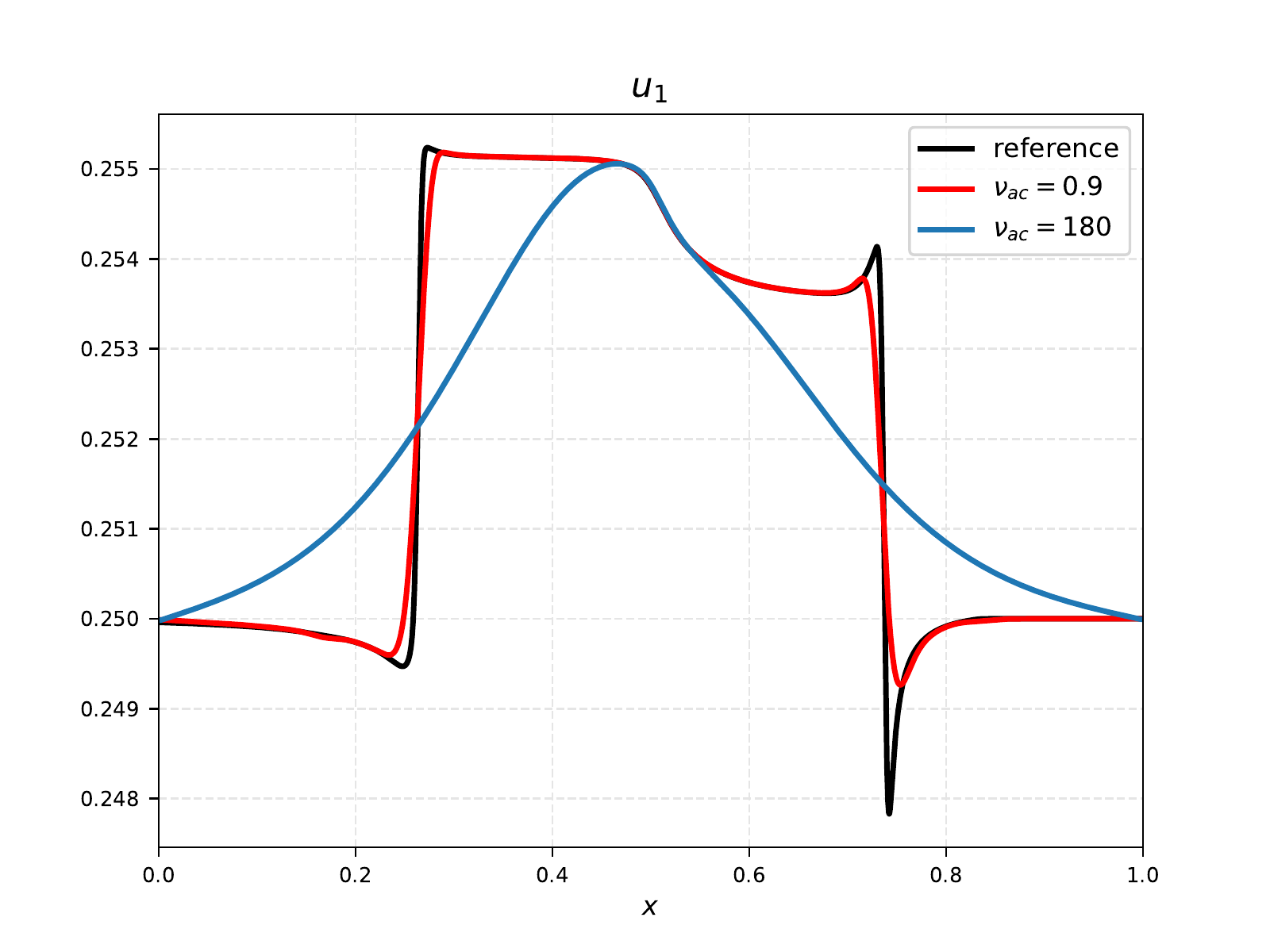}\\
	 		\includegraphics[scale=0.4]{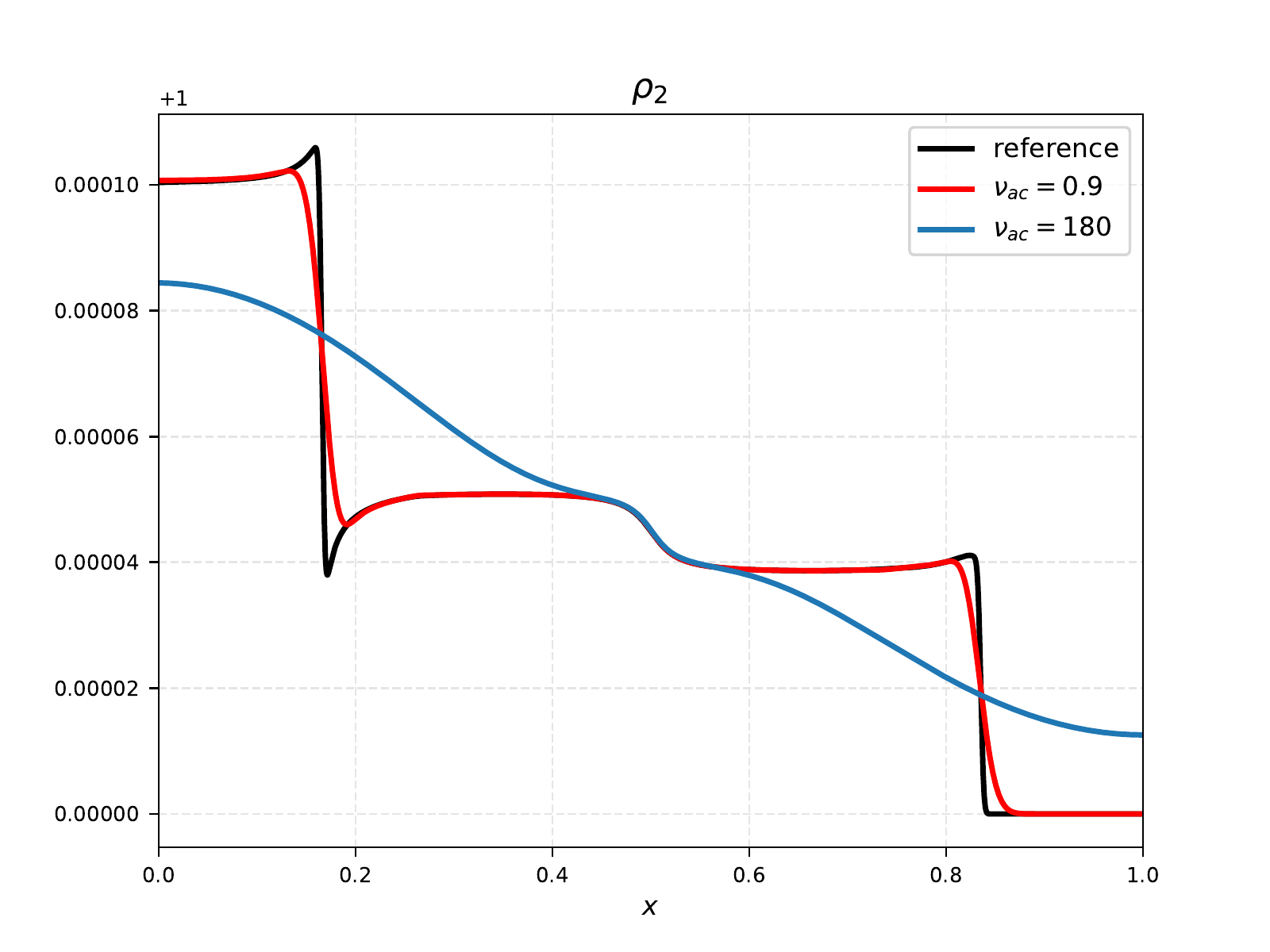}
	 		\includegraphics[scale=0.4]{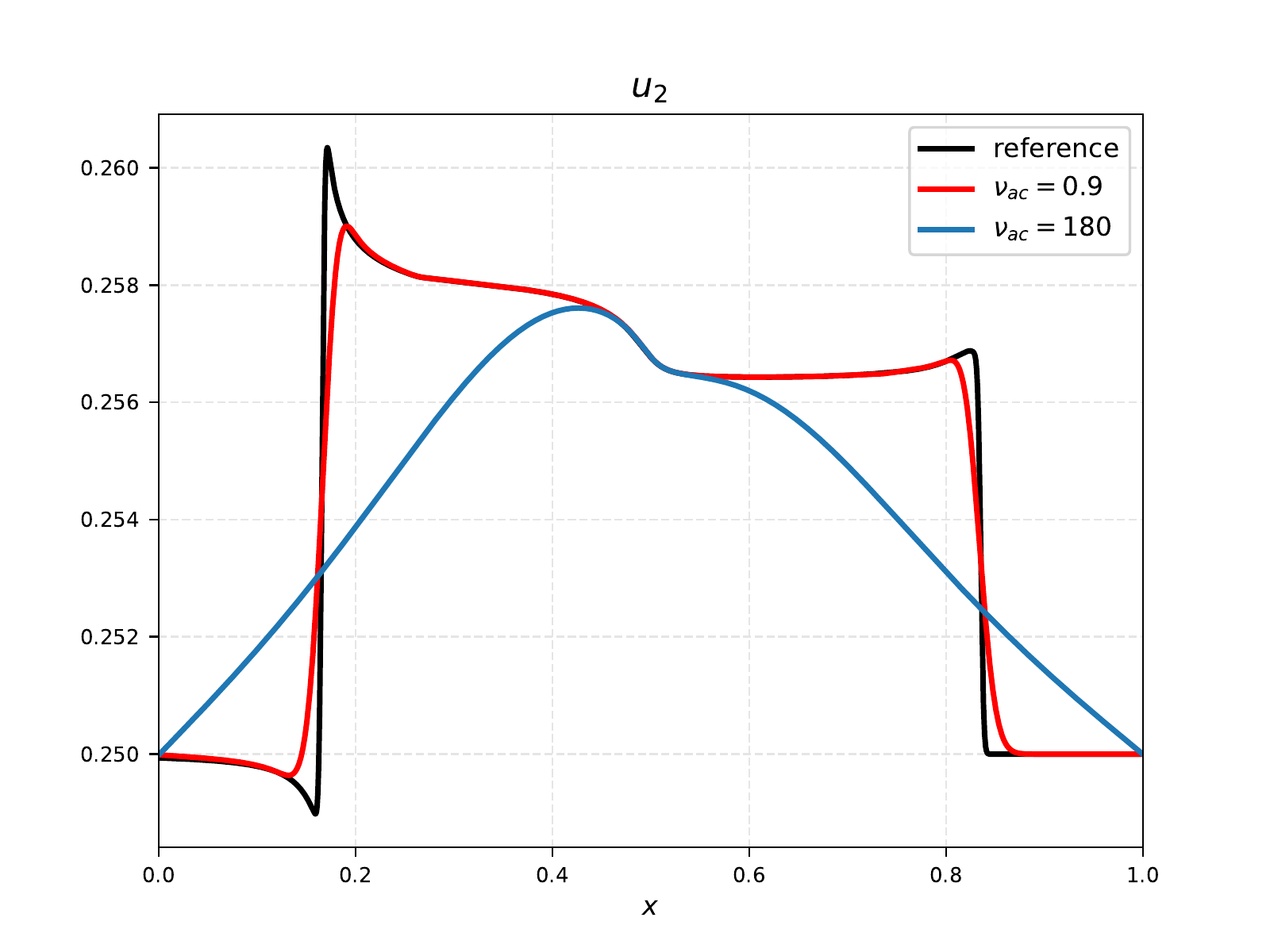}\\
	 		\includegraphics[scale=0.4]{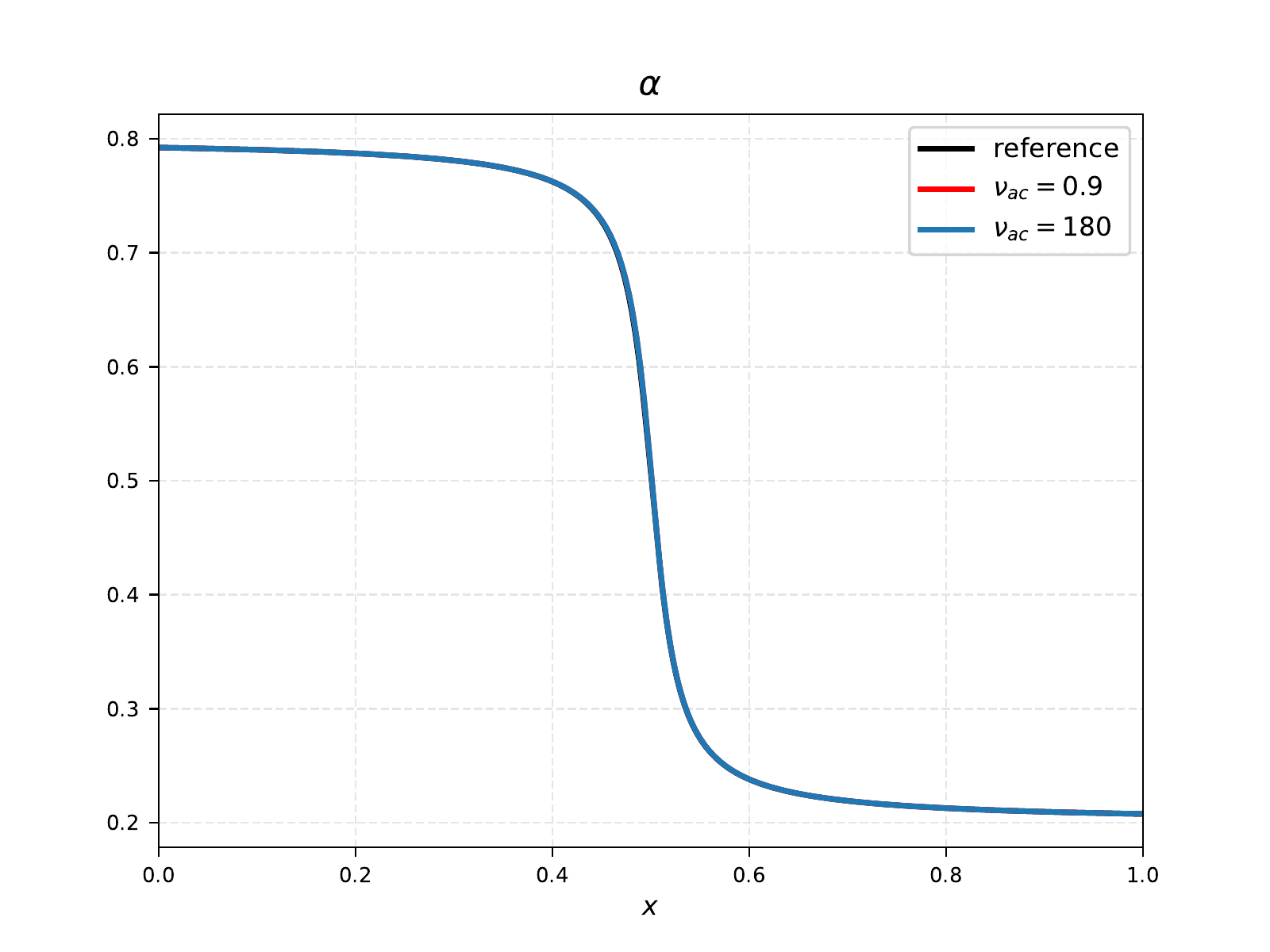}
	 		\includegraphics[scale=0.4]{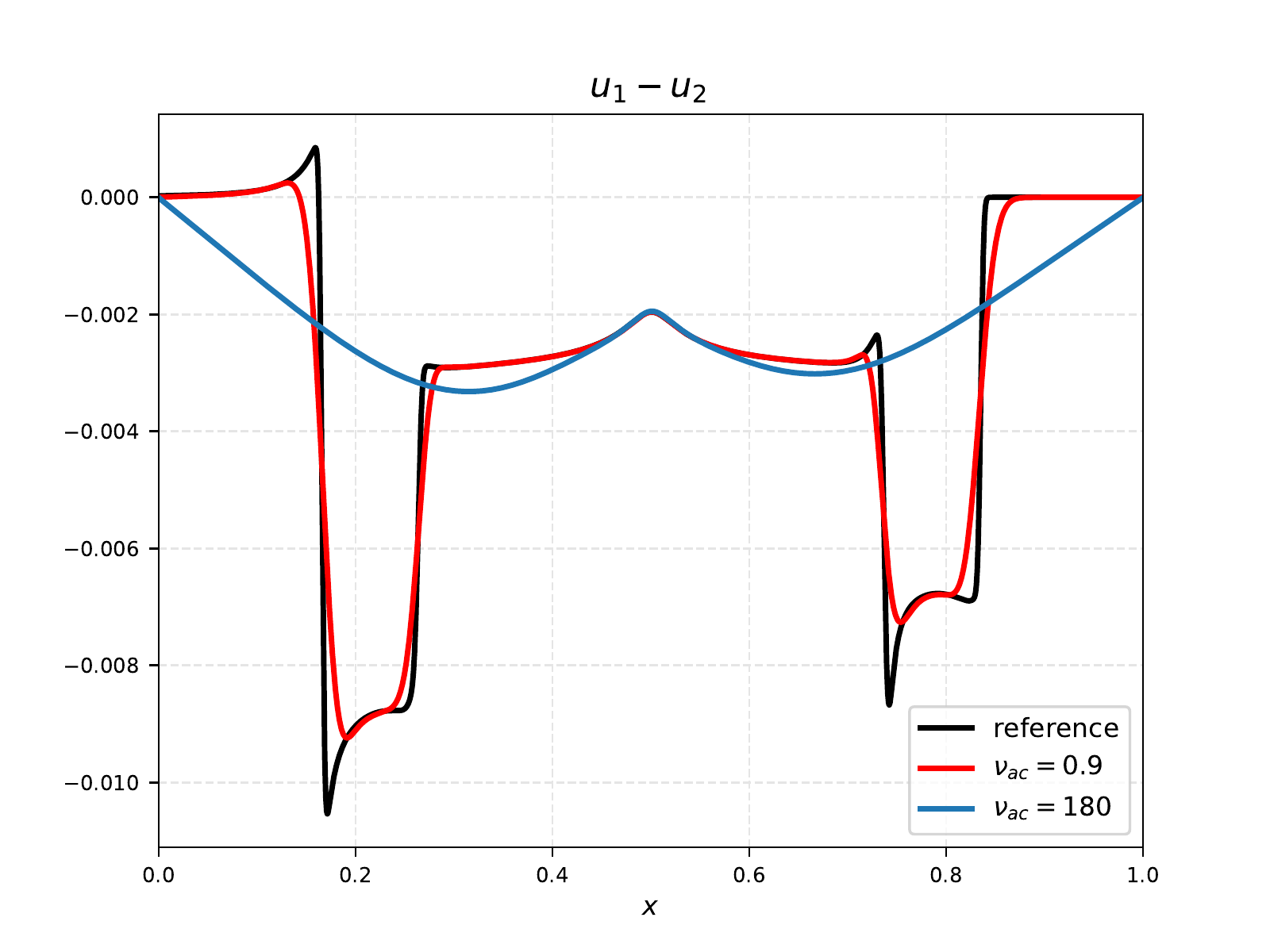}
	 	\end{center}
	 	\caption{Riemann Problem from Section \ref{sec:NumRes:RP} with smooth initial $\alpha$: Numerical results for mixture density $\rho$, mixture velocity $u$ and relative velocity $u_1 - u_2$ for $M=10^{-2}$ at final time $T_f = 0.2 M$ with grid size $\Delta x = 10^{-3}$.}
	 	\label{fig:RPMSameSmoothA}
	\end{figure}

	We repeat this test using an initial smooth tangent transition between $\alpha_L$ and $\alpha_R$. 
	The results are given in Figure \ref{fig:RPMSameSmoothA}. 
	Analogously to the previous case, the RS-IMEX scheme \eqref{eq:Time_Semi_scheme} captures accurately the diffusive interface even for large time steps with $\nu_{ac}=180$. 
 	
	\textbf{Different Mach number regimes.}
	The next numerical test concerns the Riemann problem with initial data \eqref{eq:RP} for different flow regimes of the respective phase given by $M_1 = 10^{-2}, ~M_2 = 10^{-3}$. 
	We consider a well-prepared homogeneous mixture governed by the homogeneous model with constant $\alpha = 0.9$ and $\rho_1^{(2)} = 1$ and $p_2^{(2)} = p_1^{(1)}$. 
	Since the acoustic waves of phase two are significantly faster than the ones of phase one, we consider a larger computational domain given by $[0,50]$ and $\Delta x = 10^{-2}$ and the final time $T_f = 7.5 \cdot 10^{-2}$.  
 
\begin{figure}[t!]
	\begin{subfigure}[c]{\textwidth}
		\begin{center}
		\includegraphics[scale=0.41]{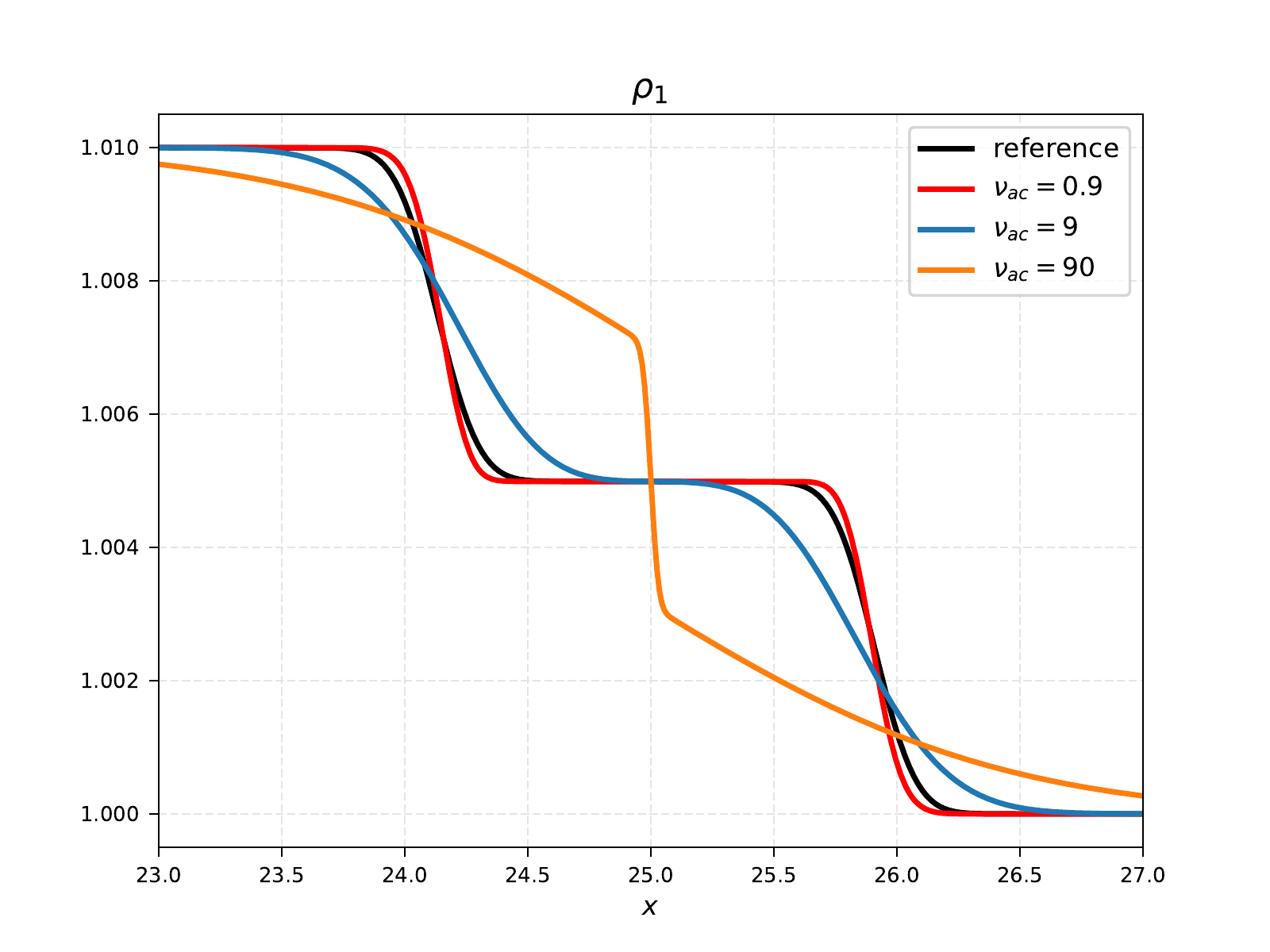}
		\includegraphics[scale=0.41]{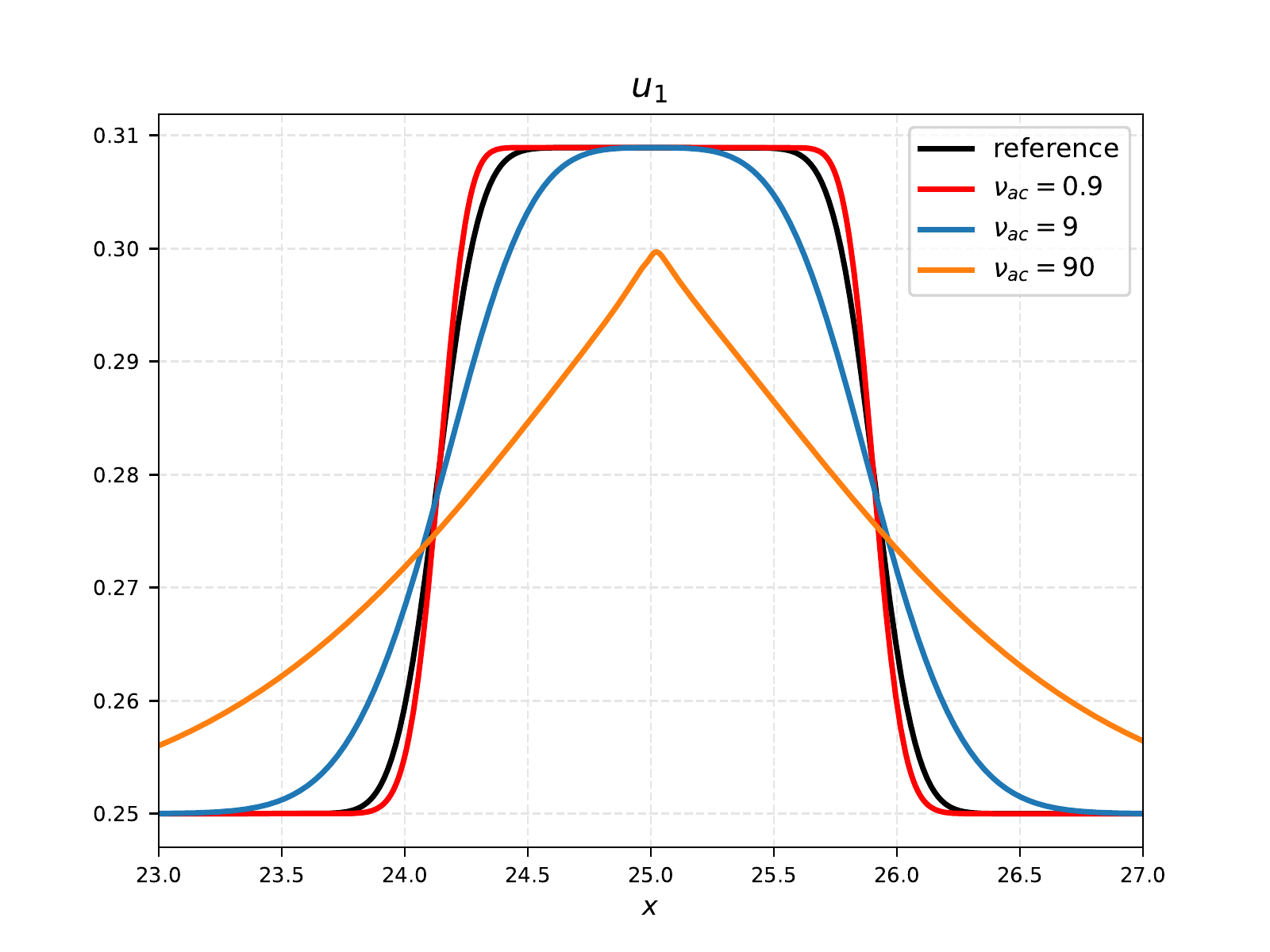}
		\end{center}
		\subcaption{Zoom on acoustic waves associated to phase one. }
		\label{fig:RPMDiffMat}
	\end{subfigure}
	\begin{subfigure}[c]{\textwidth}
		\begin{center}
		\includegraphics[scale=0.4]{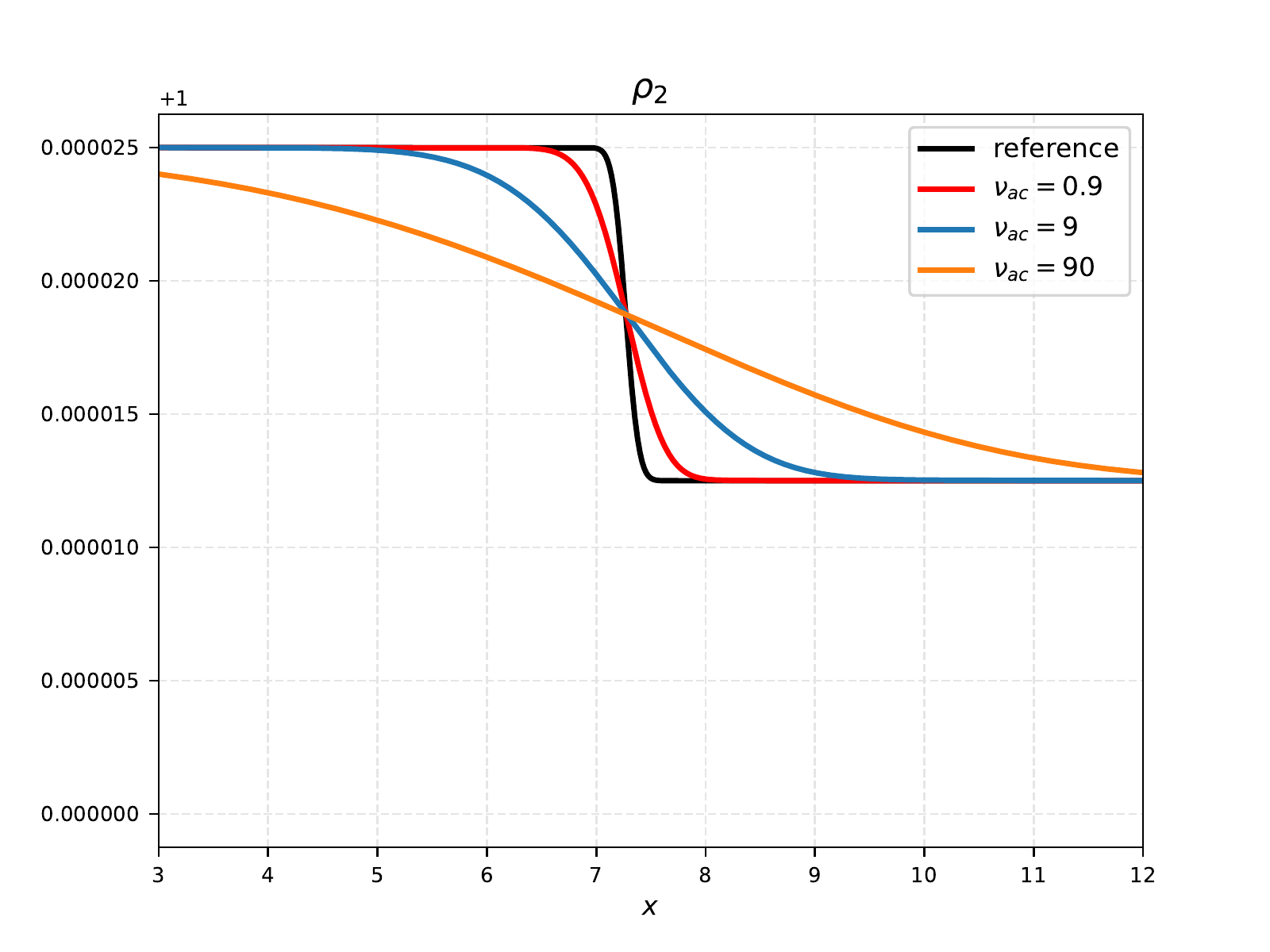}
		\includegraphics[scale=0.4]{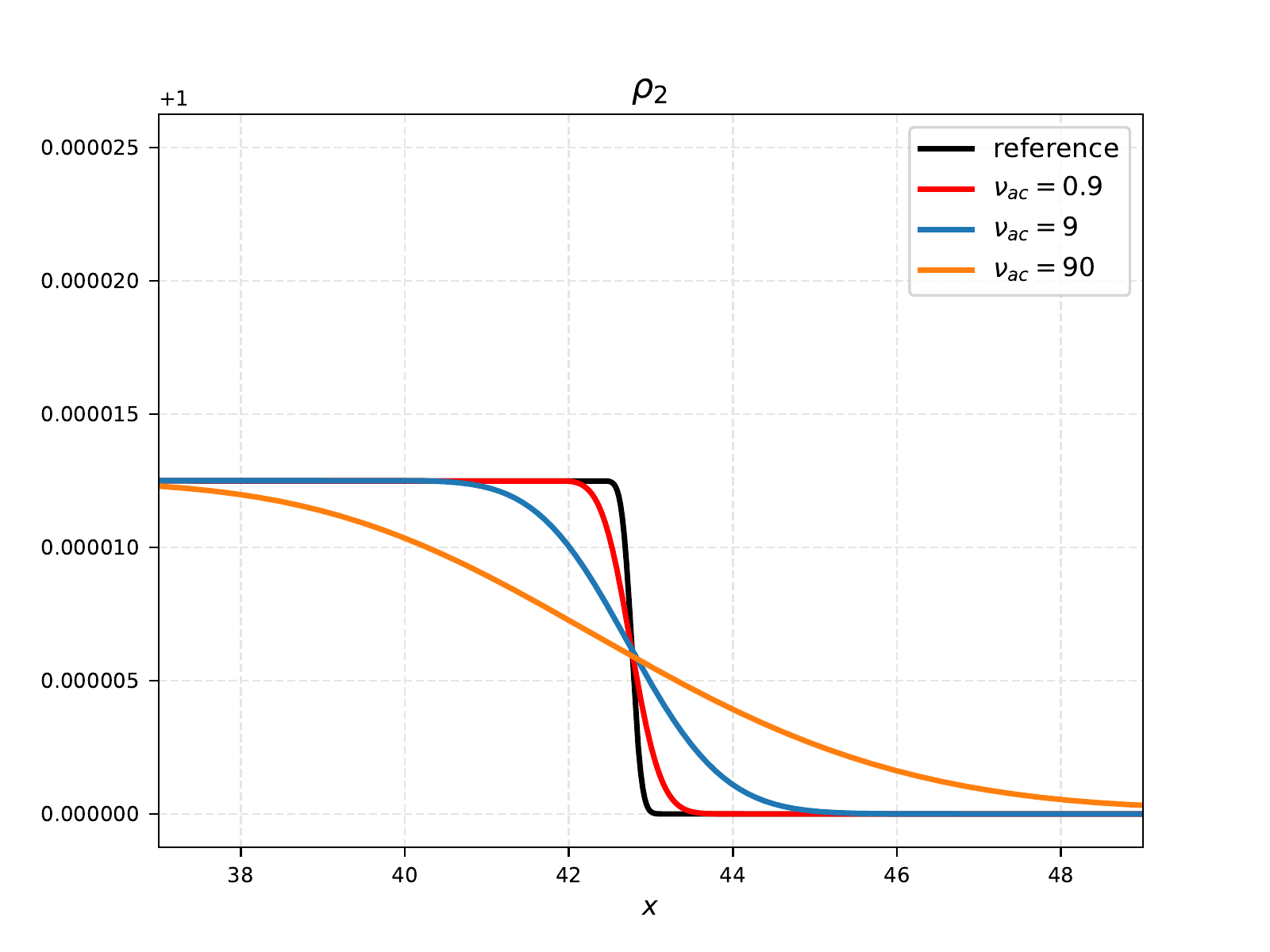}\\
		\includegraphics[scale=0.4]{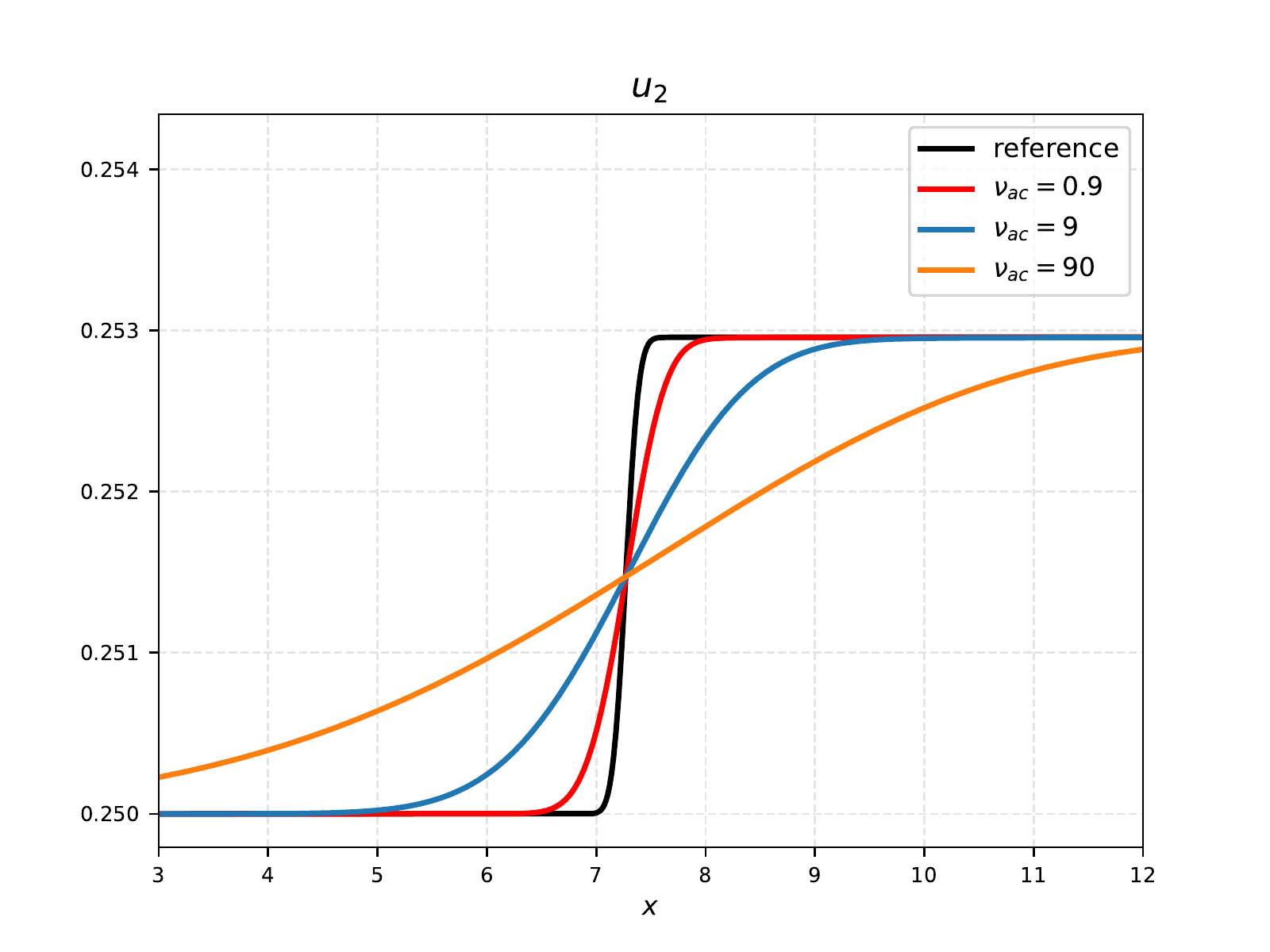}
		\includegraphics*[scale=0.4]{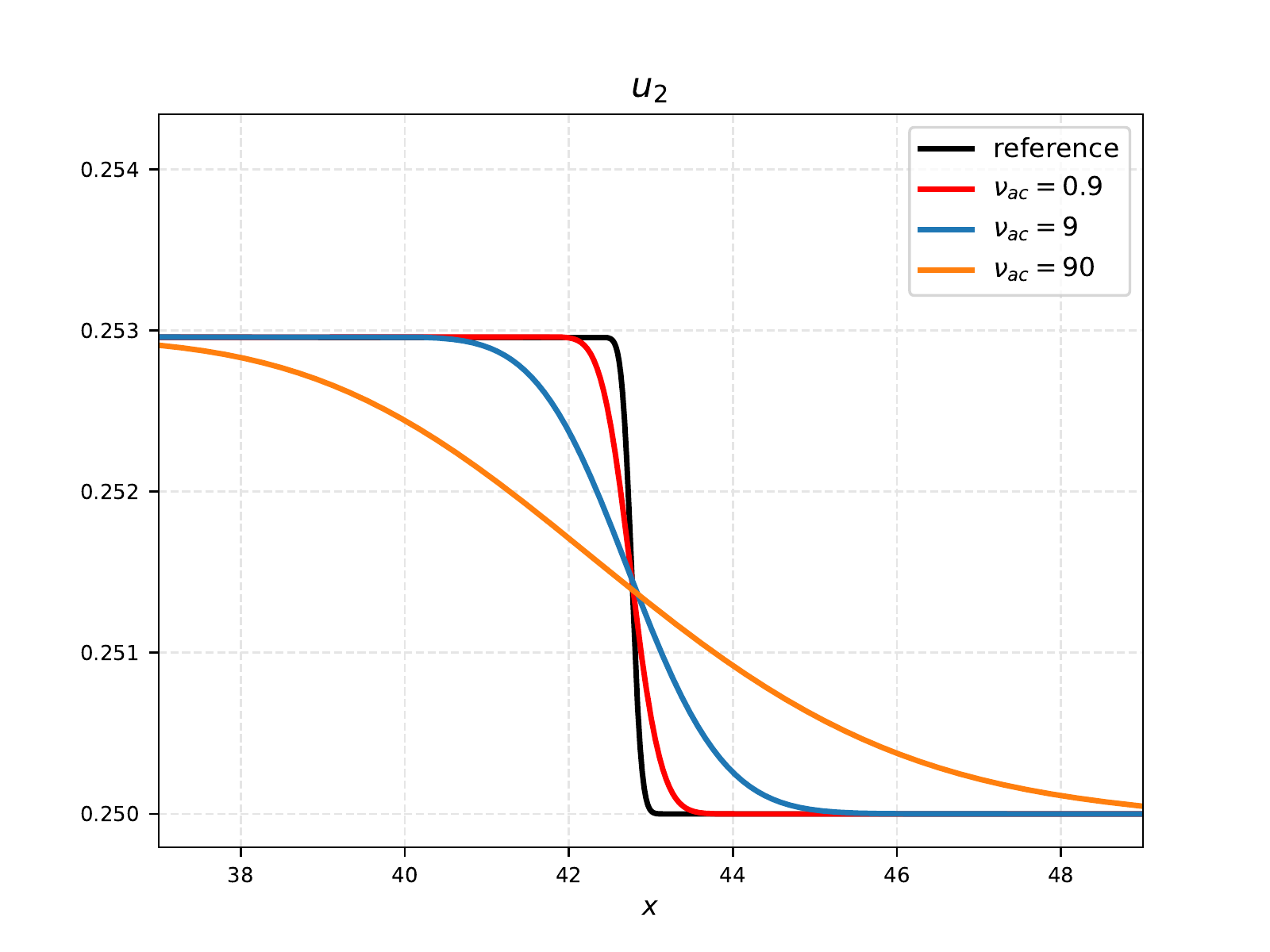}
		\end{center}
		\subcaption{Zoom on acoustic waves associated to phase two.}
		\label{fig:RPMDiffAcc}
	\end{subfigure}
	\caption{Riemann Problem from Section \ref{sec:NumRes:RP} for two phases in different Mach number regimes: Numerical results for $M_1=10^{-2}$ and $M_2=10^{-3}$ at final time  $T_f = 7.5\cdot10^{-2}$ and grid size $\Delta x = 10^{-2}$.}
	\label{fig:RPMDiff}
	\end{figure}

	The numerical results are displayed in Figure \ref{fig:RPMDiff}. 
	Figure \ref{fig:RPMDiffMat} presents a zoom on the acoustic waves of phase one in the domain $[23,27]$ whereas the acoustic waves of phase two have already reached the domains $[3,12]$ and $[37,49]$ as depicted in Figure \ref{fig:RPMDiffAcc}. 
	We run the simulation with three different time steps. 
	The first time step is focused on the resolution of all acoustic waves with $\nu_{ac} = 0.9$ given by $\Delta t = 2.1 \cdot 10^{-5}$. 
	With the second time step given by $\nu_{ac} = 9$ leading to $\Delta t = 2.1\cdot 10^{-4}$, the fast acoustic waves of phase two are smoothened but the acoustic waves of phase one are resolved well. 
	The last time step given by $\nu_{ac} = 90$ leading to $\Delta t = 2.1\cdot 10^{-3}$ is oriented to the material wave neglecting the resolution of all acoustic waves. 

\begin{figure}[t!]
		\begin{subfigure}[c]{\textwidth}
			\begin{center}
			\includegraphics[scale=0.4]{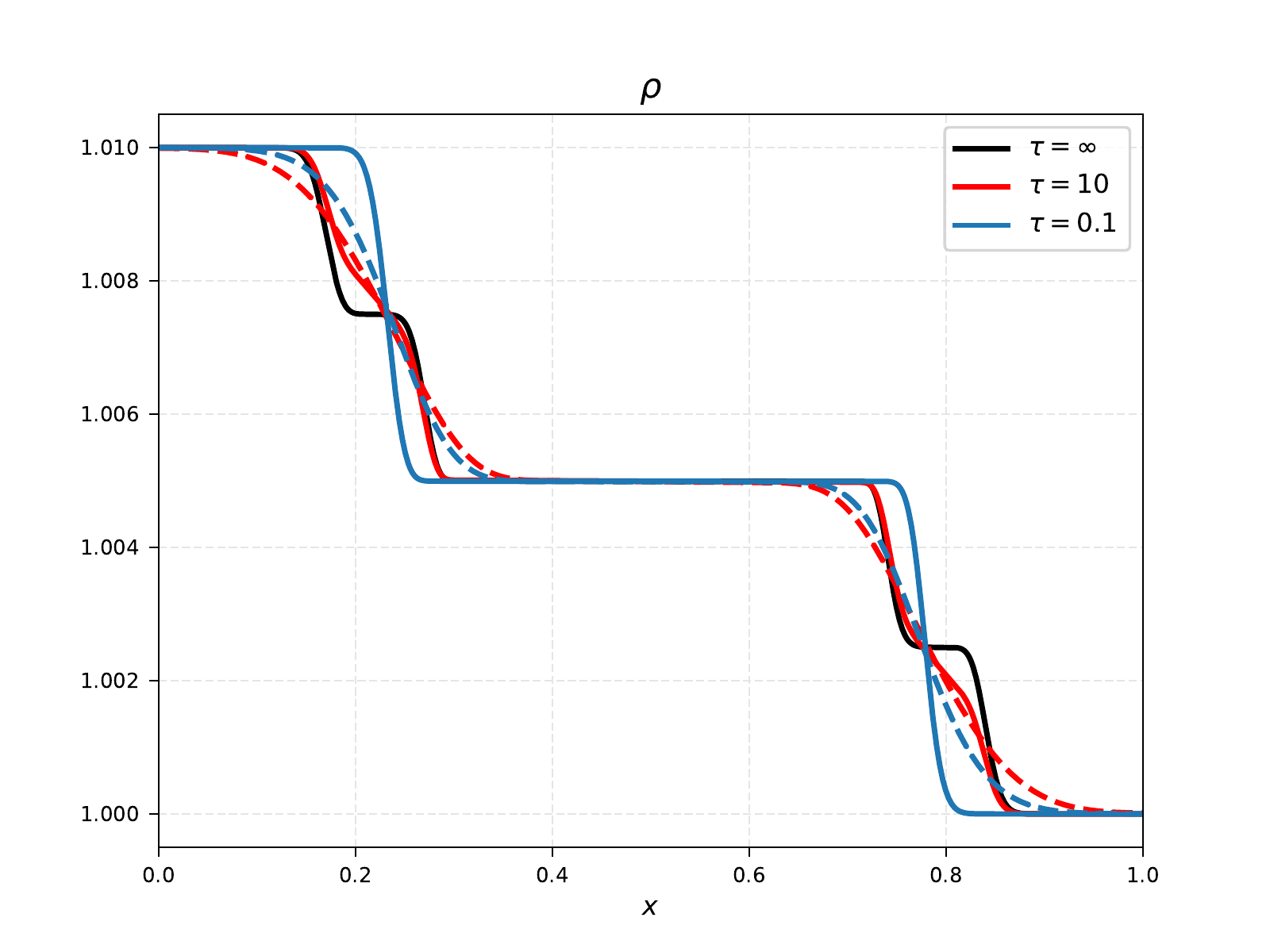}
			\includegraphics[scale=0.4]{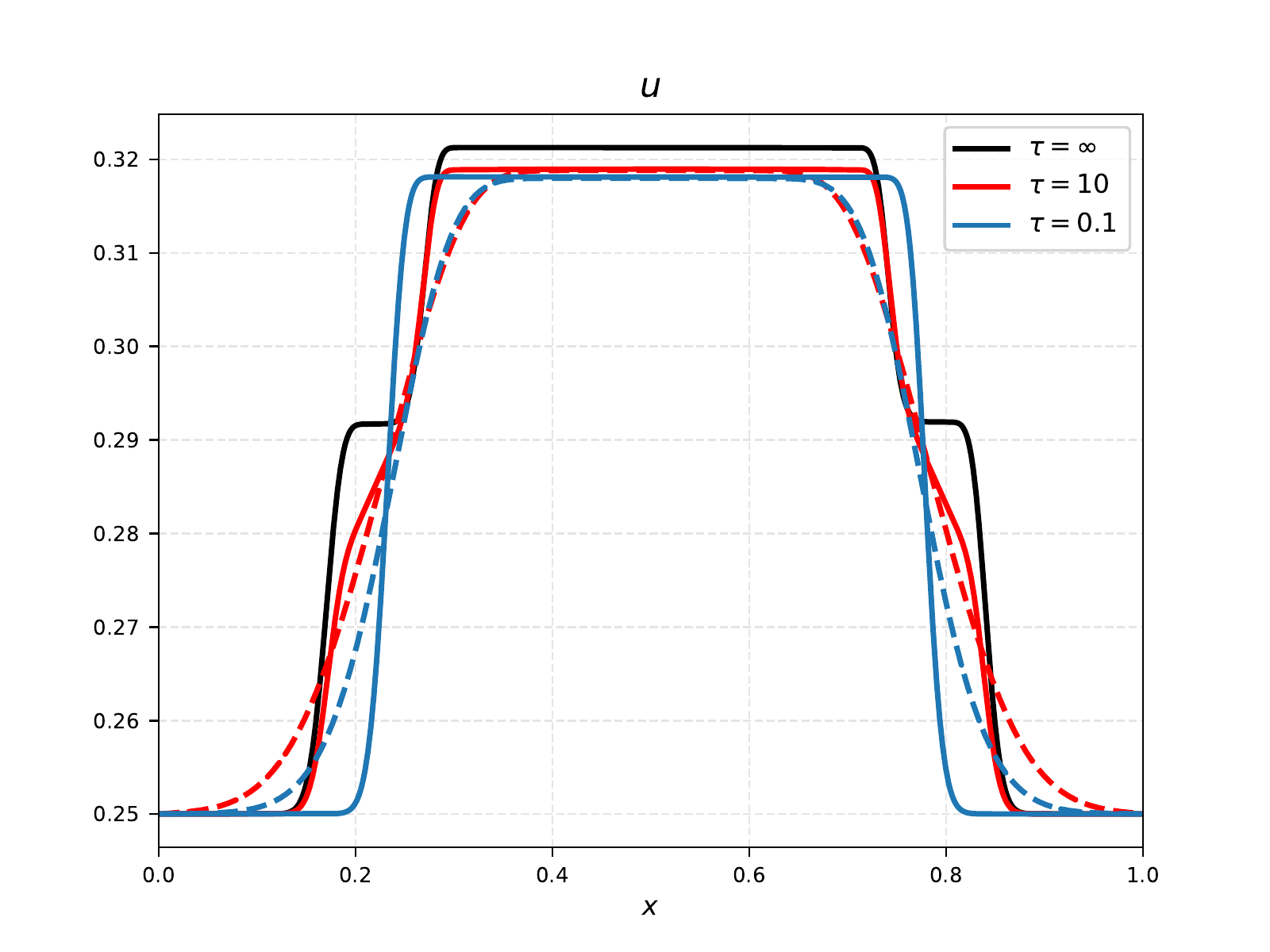}\\
			\includegraphics[scale=0.4]{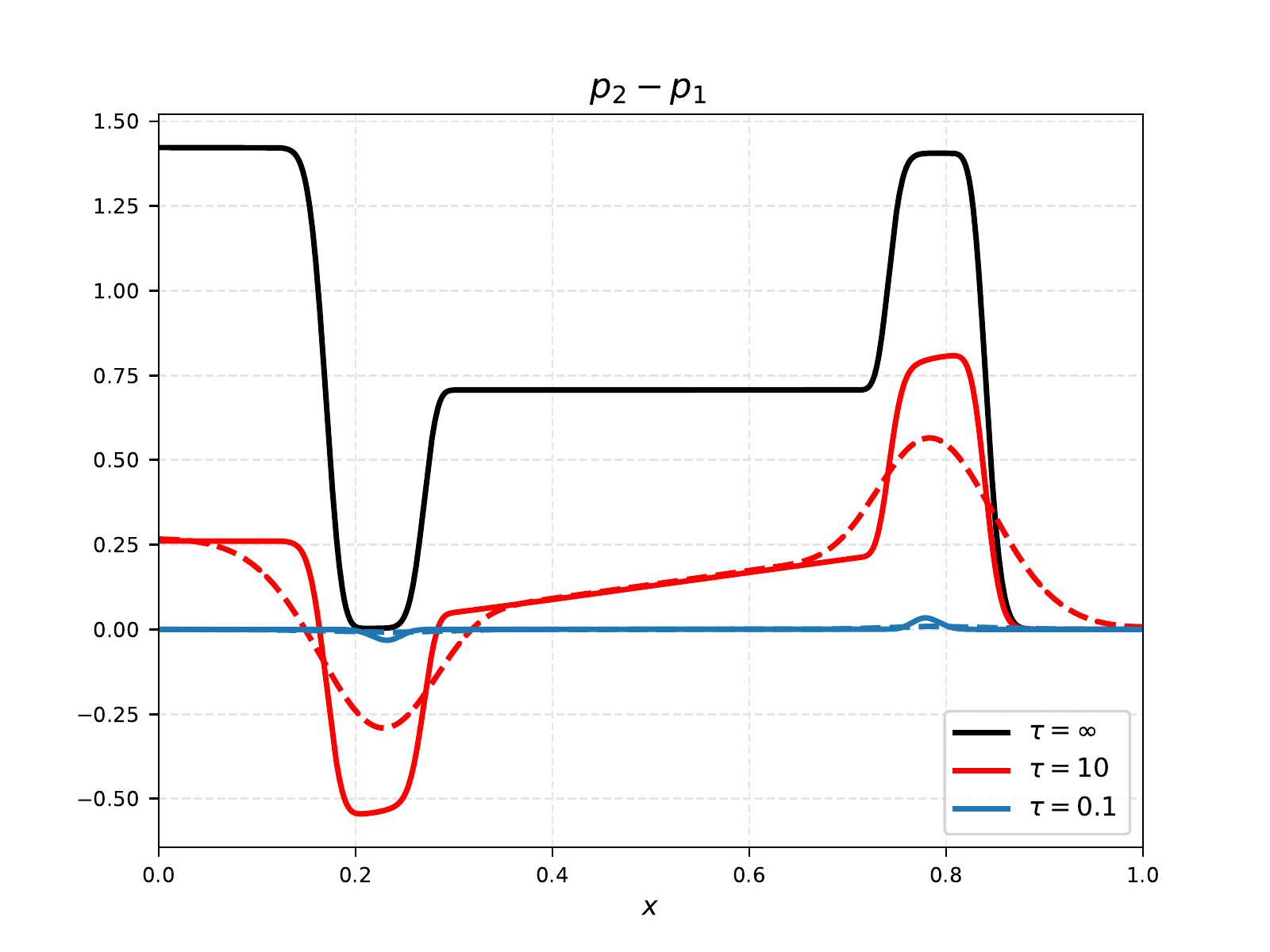}
			\includegraphics[scale=0.4]{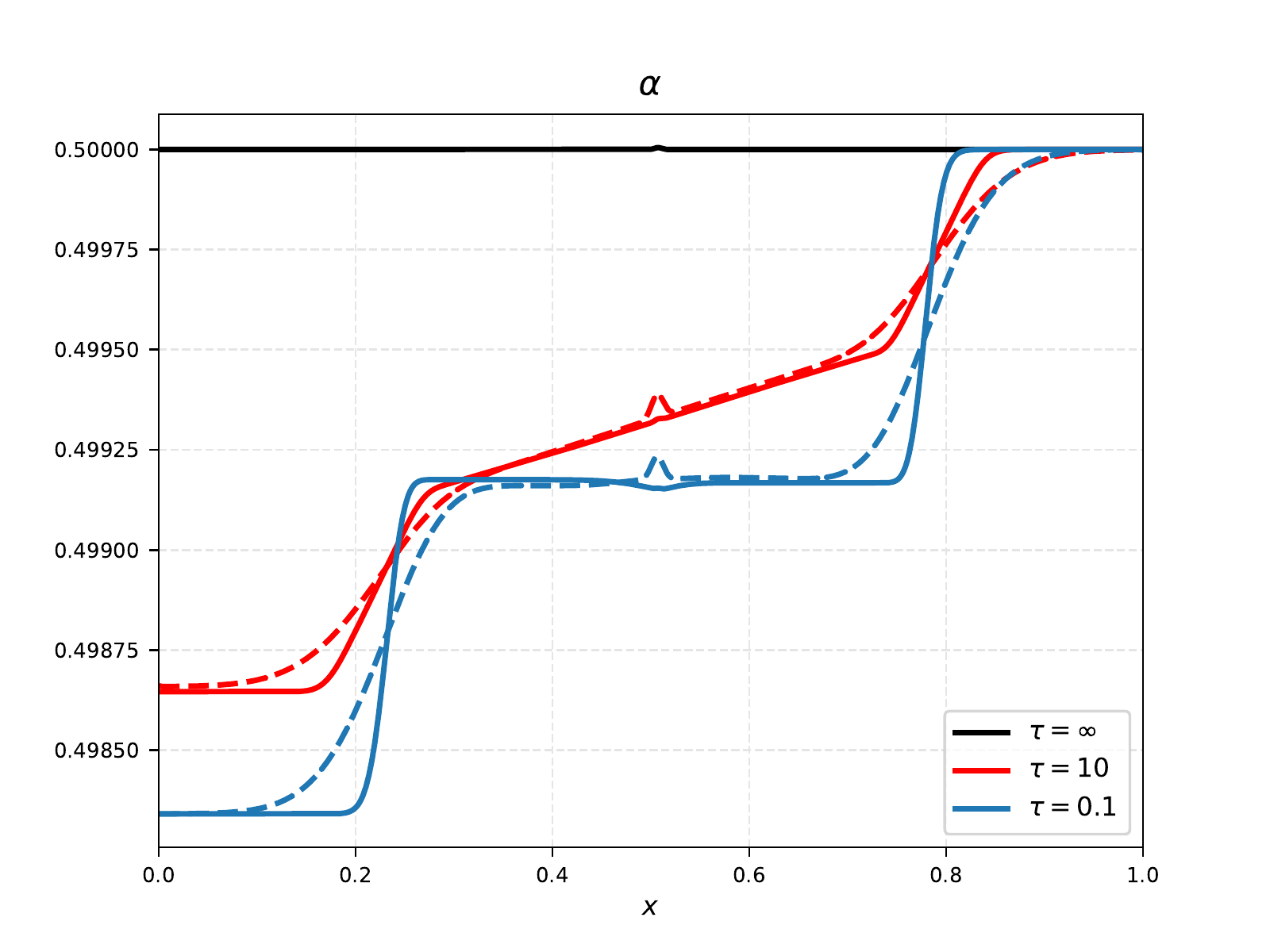}
			\end{center}
			\subcaption{$M=10^{-1}$}
		\end{subfigure}
		\begin{subfigure}[c]{\textwidth}
			\begin{center}
				\includegraphics[scale=0.4]{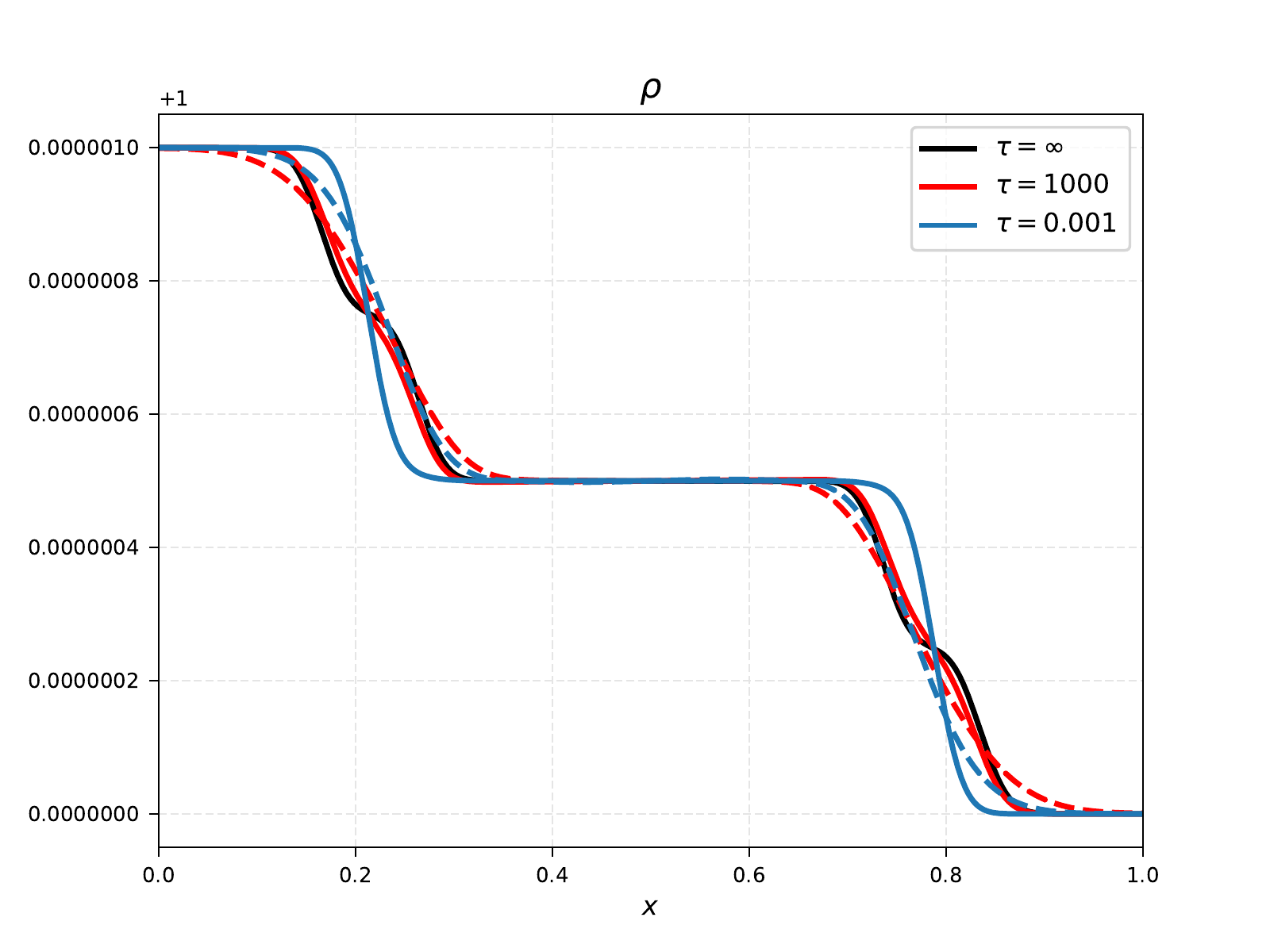}
				\includegraphics[scale=0.4]{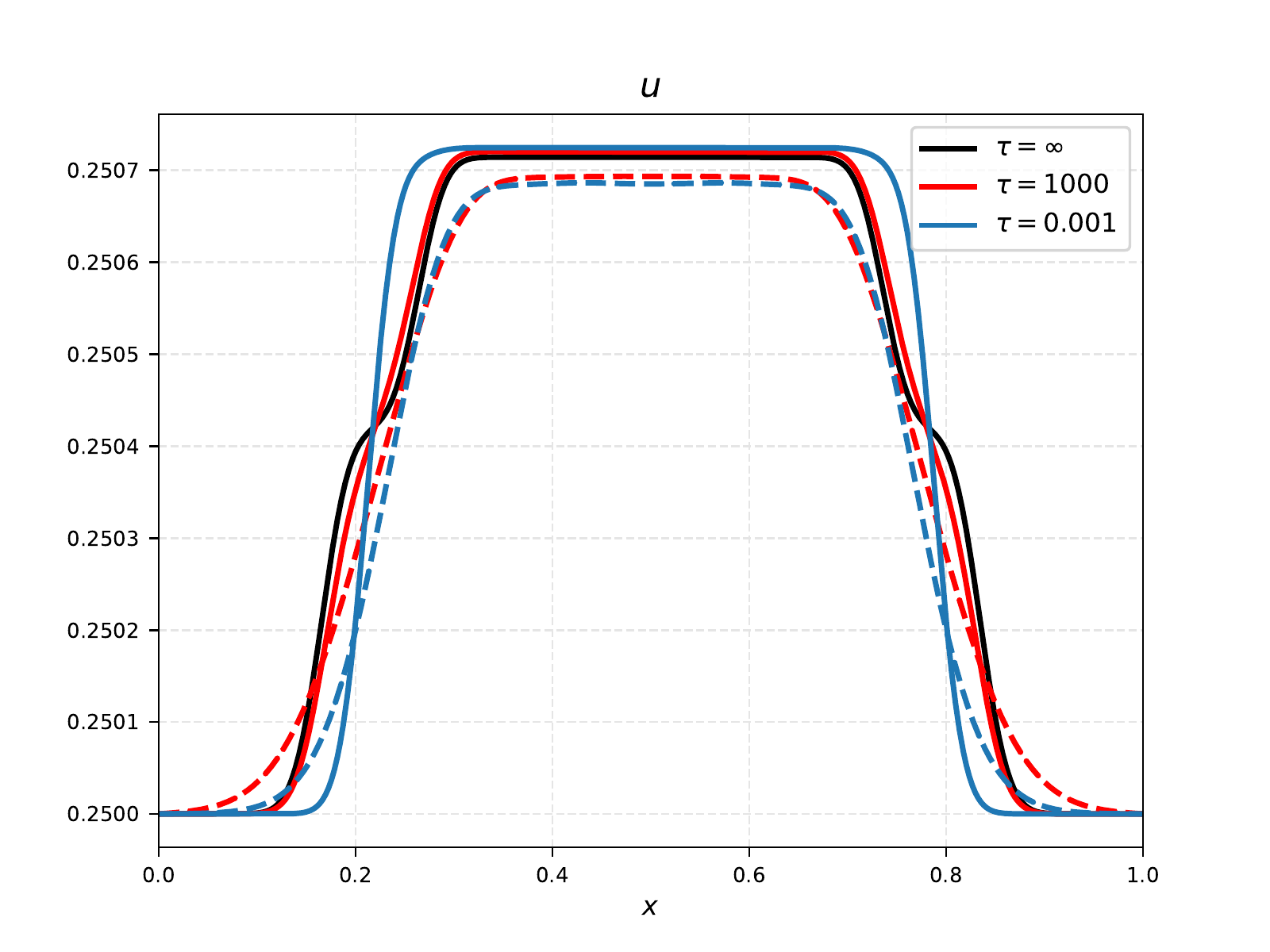}\\
				\includegraphics[scale=0.4]{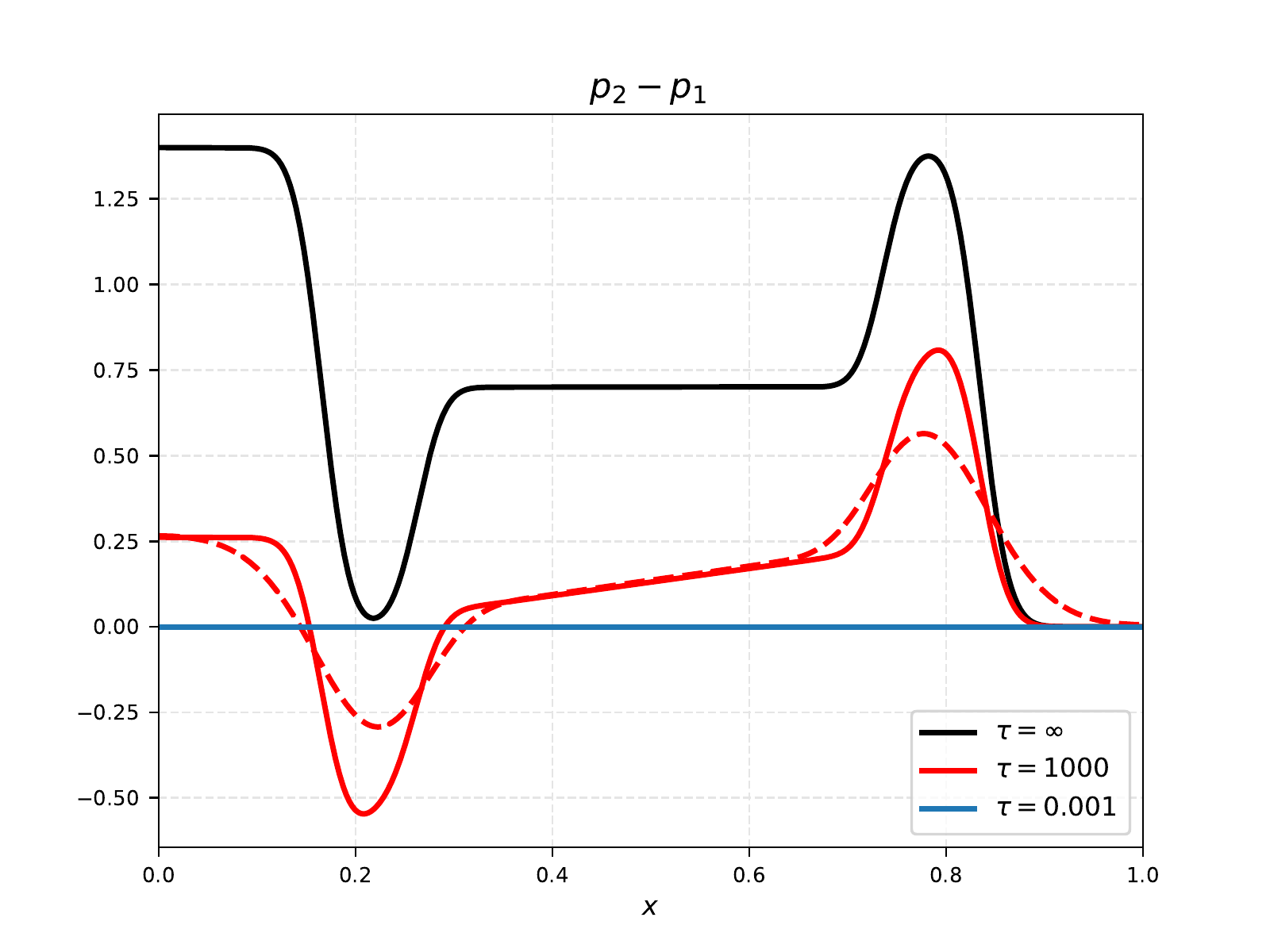}
				\includegraphics[scale=0.4]{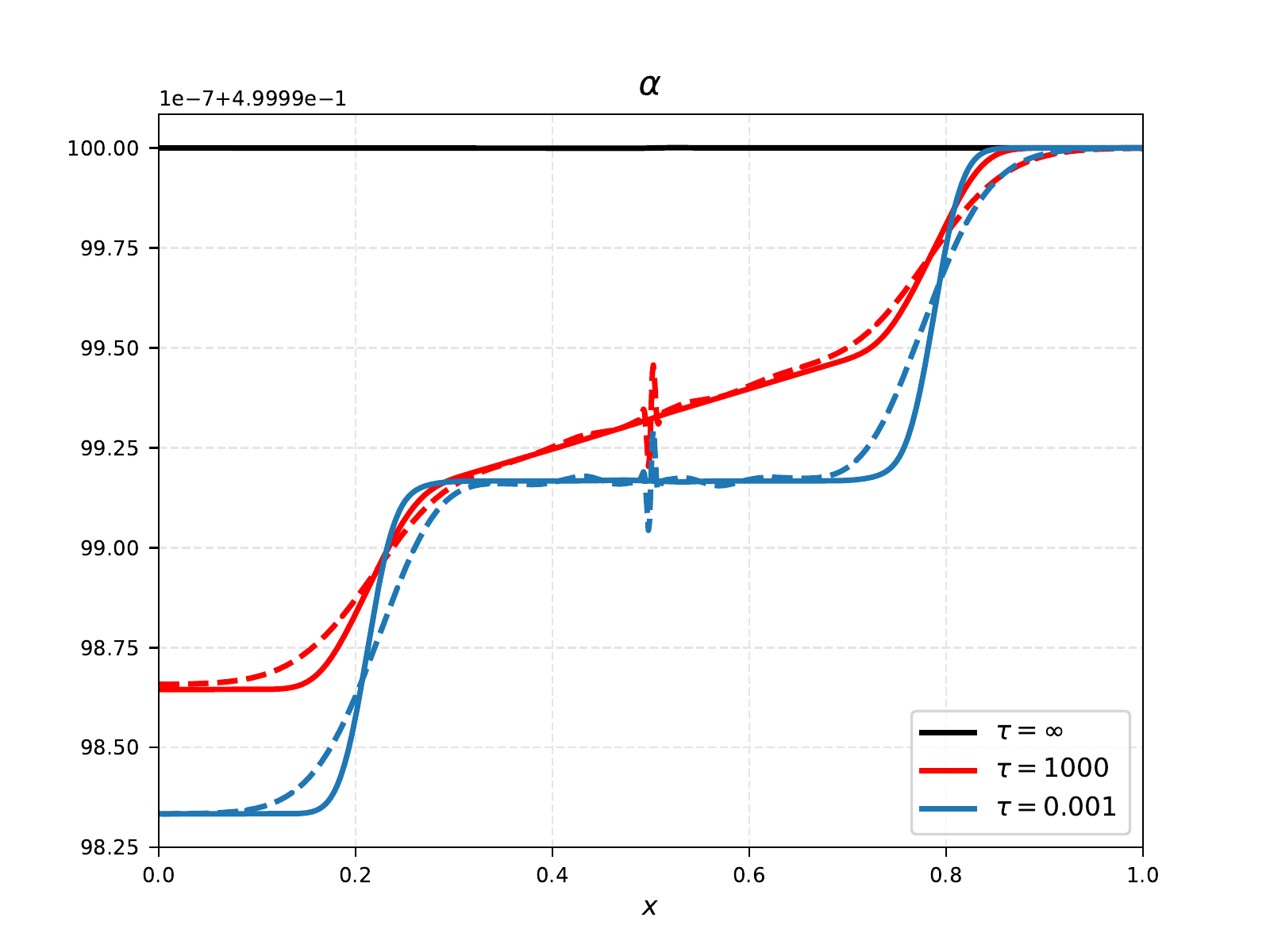}
			\end{center}
			\subcaption{$M=10^{-3}$}
		\end{subfigure}
		\caption{Riemann Problem from Section \ref{sec:NumRes:RP} with pressure relaxation: Numerical results for mixture density $\rho$, mixture velocity $u$, pressure difference $p_2 - p_1$ and volume fraction $\alpha$ for different values of $\tau$ computed with an acoustic time stepping (solid line) and a material time stepping (dashed line).}
		\label{fig:RPMSamePRelax}
	\end{figure}

	\textbf{Influence of pressure relaxation.}
	We consider an initial homogeneous mixture with $M_1 = M = M_2$ and $\alpha =0.5$ with pressure relaxation source term. 
	The initial condition is given by the Riemann problem \eqref{eq:RP}. 
	We consider different values for the relaxation time $\tau = \infty, M, M^{-1}$ for two different Mach number regimes $M = 10^{-1}, 10^{-3}$. 
	As the pressure relaxation acts on the volume fraction $\alpha$, we expect a change in $\alpha$ depending on the relaxation time. 
	The numerical results are presented in Figure \ref{fig:RPMSamePRelax}. 
	We can see that for $\tau = M$, the difference in the pressures goes to zero, leading to a change in $\alpha$ of order $M^2$. 
	
\begin{figure}[t!]
		\begin{subfigure}[c]{\textwidth}
			\includegraphics[scale=0.33]{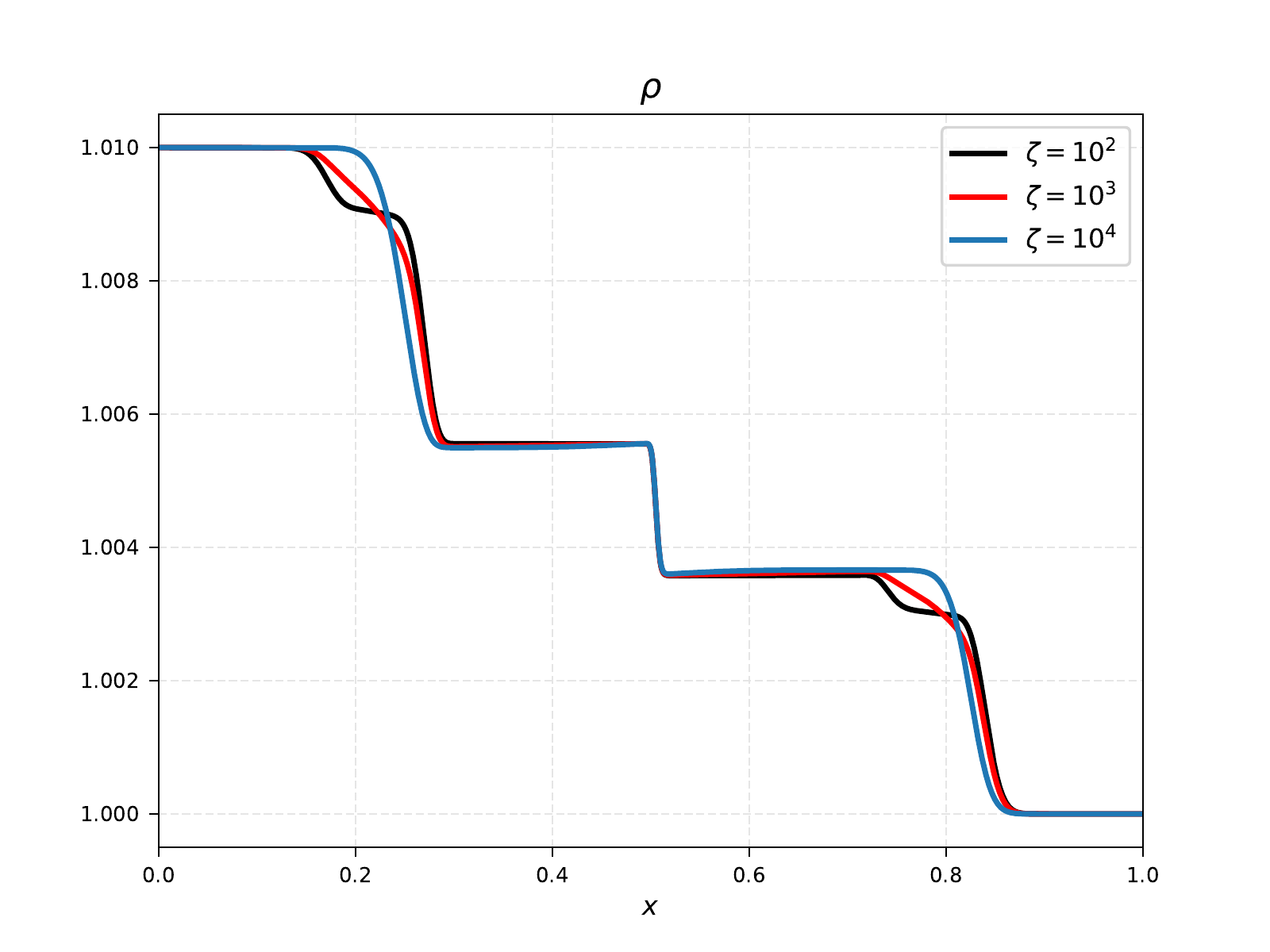}
			\includegraphics[scale=0.33]{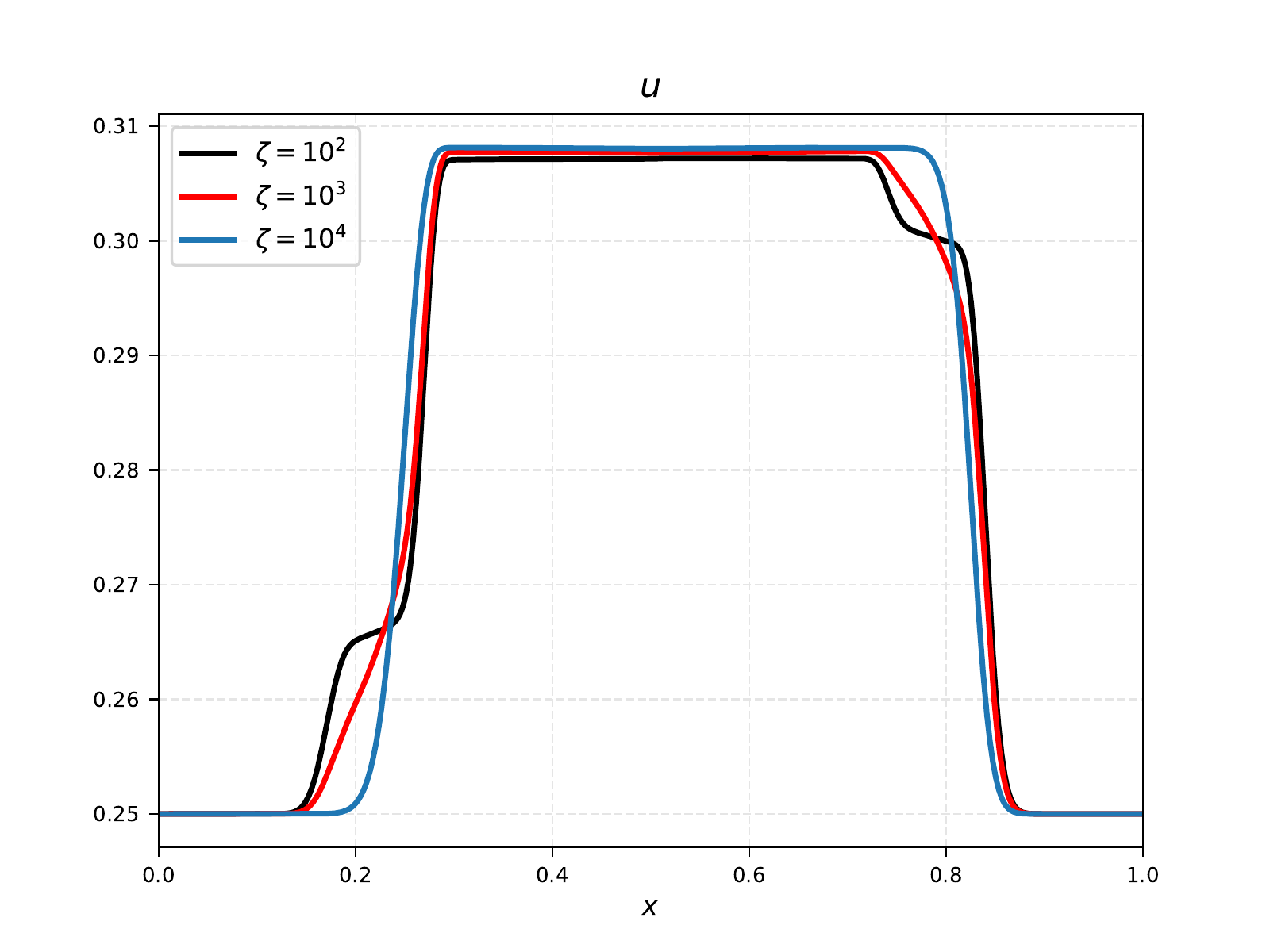}
			\includegraphics[scale=0.33]{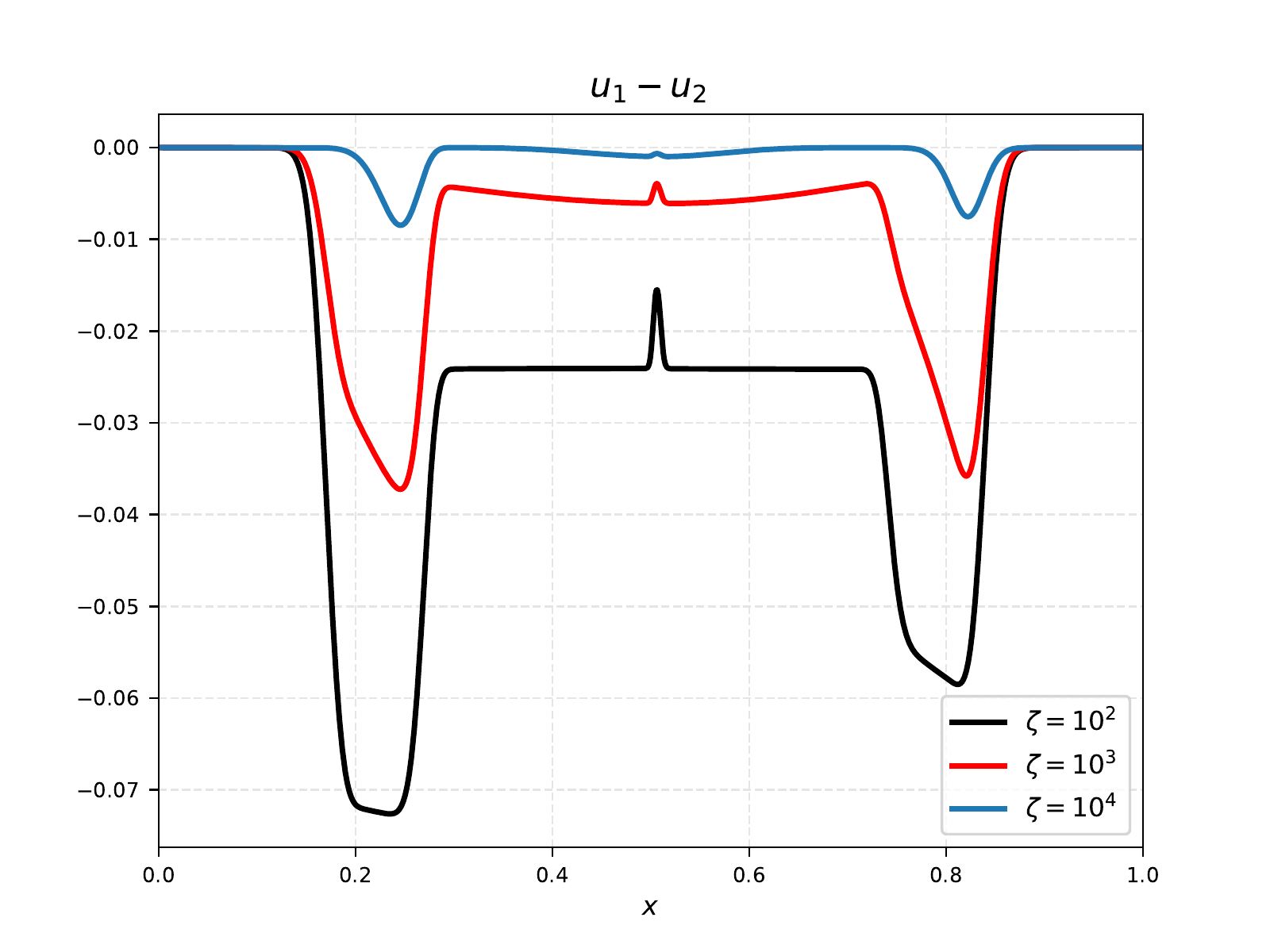}
			\subcaption{$M=10^{-1}$}
		\end{subfigure}
		\begin{subfigure}[c]{\textwidth}
			\includegraphics[scale=0.33]{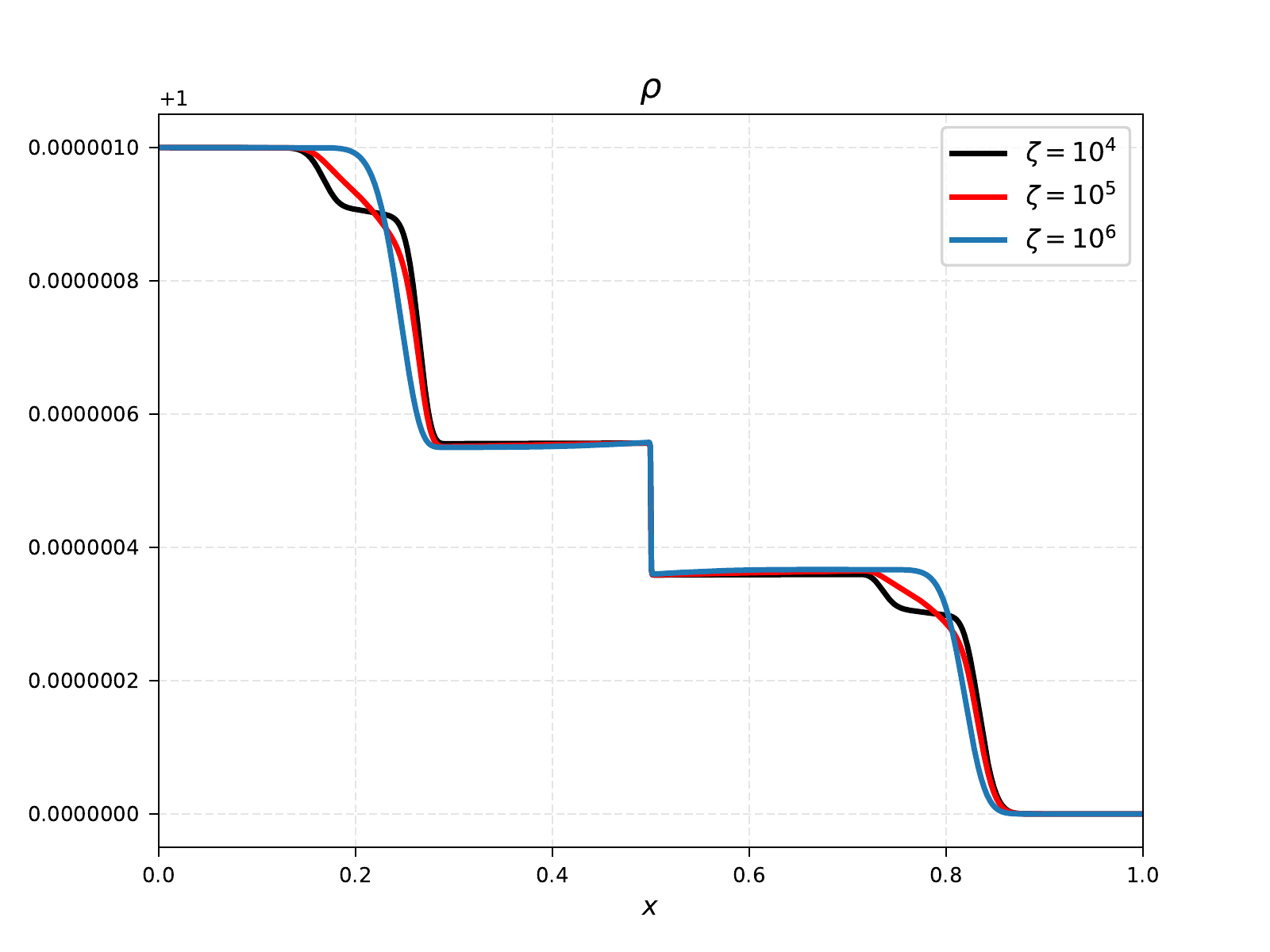}
			\includegraphics[scale=0.33]{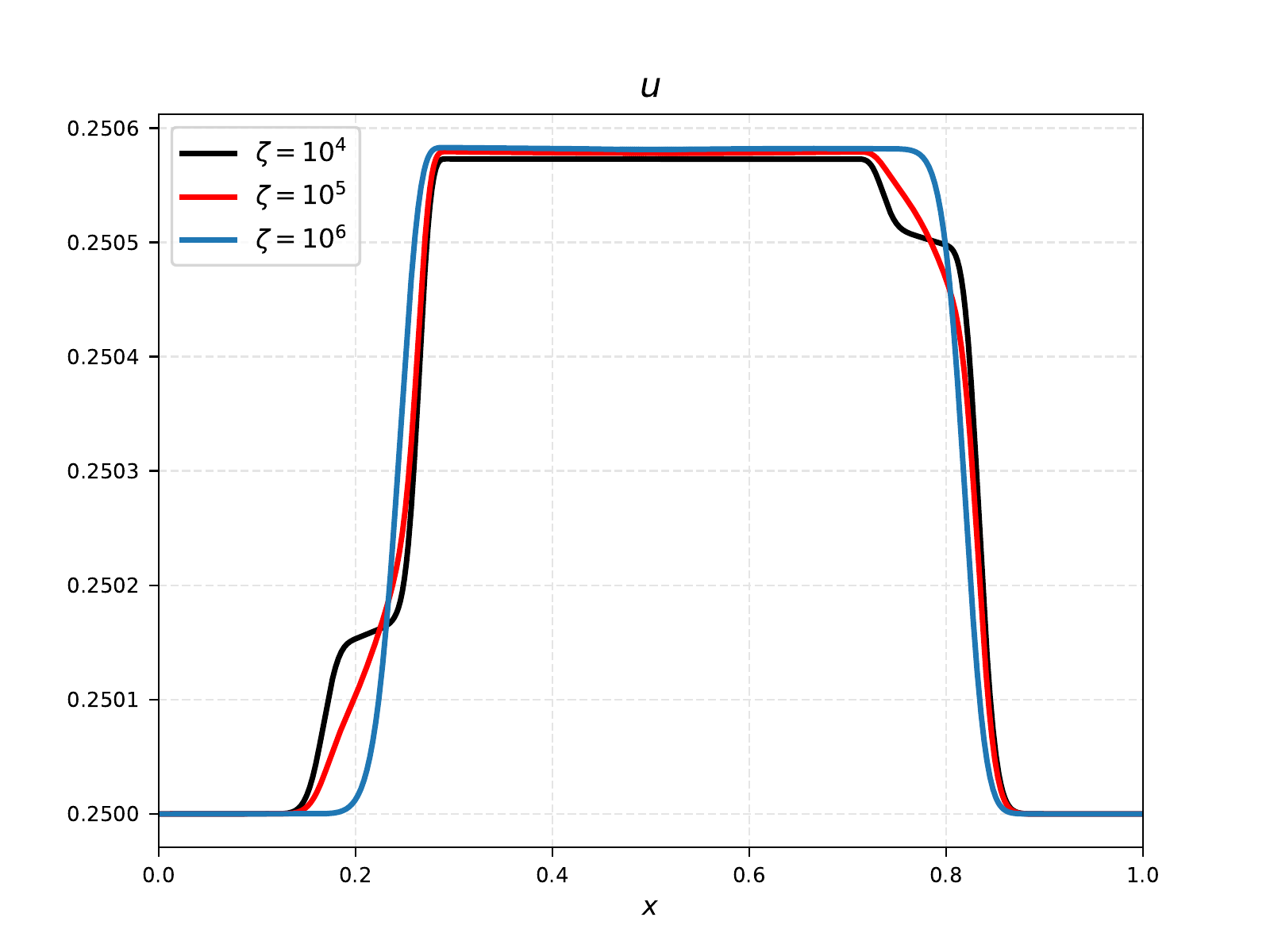}
			\includegraphics[scale=0.33]{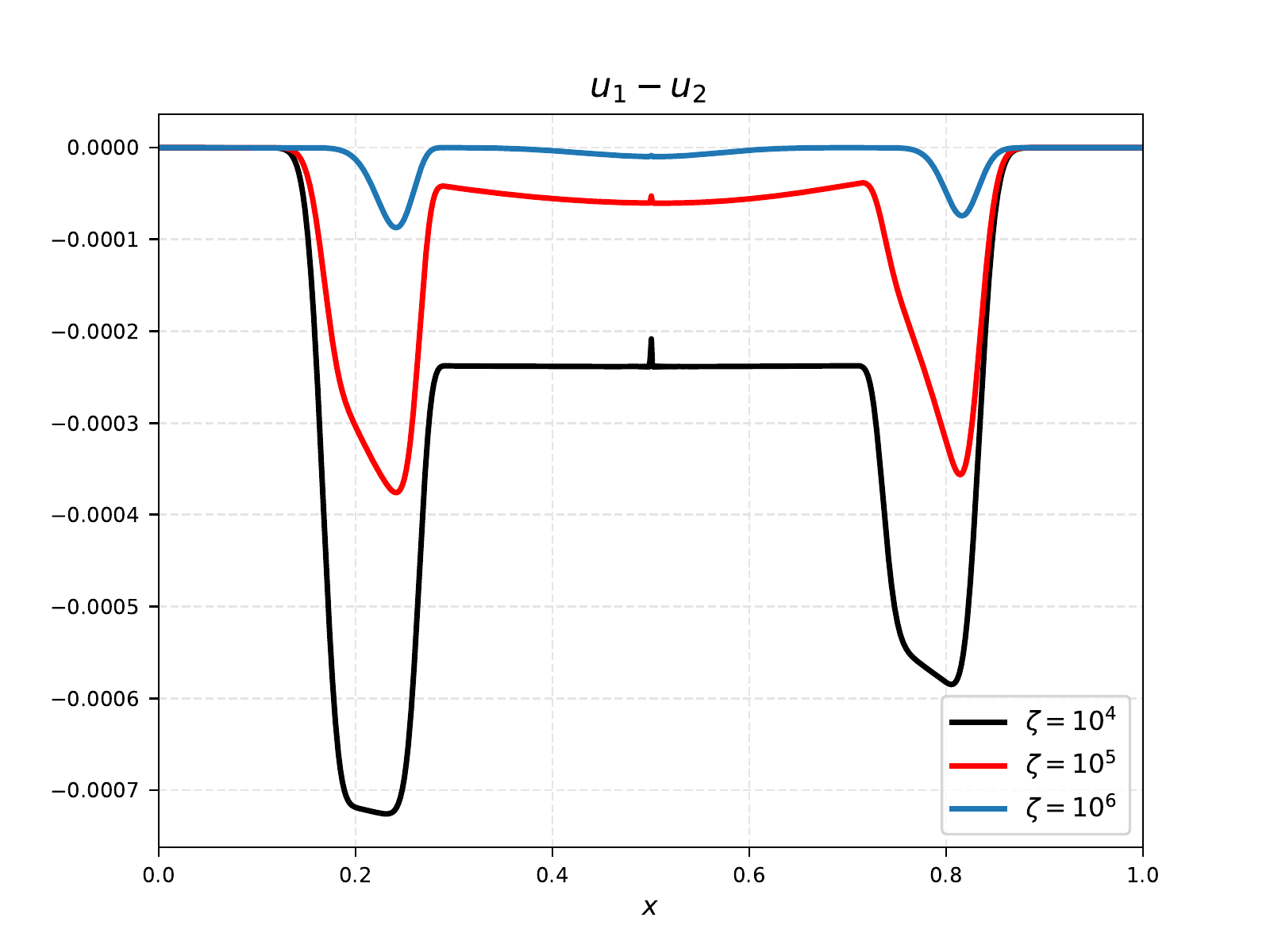}
			\subcaption{$M=10^{-3}$}
		\end{subfigure}
		\caption{Riemann Problem from Section \ref{sec:NumRes:RP} with friction and an initial jump in $\alpha$: Numerical results for mixture density $\rho$, mixture velocity $u$ and relative velocity $u_1 - u_2$ for different values of the friction coefficient $\zeta$ at final time $T_f = 0.2 M$, respectively, and grid size $\Delta x = 10^{-3}$}
		\label{fig:RPMSameJumpAFric}
	\end{figure}

	\textbf{Influence of friction.}
	In the final numerical test we rerun the test case with a jump in $\alpha$ for $M=10^{-1}, M=10^{-3}$ with the friction source term acting on the relative velocity $u_1 - u_2$. 
	We consider different values for the friction coefficent $\zeta$. 
	The numerical results are shown in Figure \ref{fig:RPMSameJumpAFric}. 
	We can see that for large friction coefficients $\zeta$ the relative velocity goes to zero, leading to similar phase velocities.
	For lower Mach number flows a larger friction coefficient is needed to obtain the same effect, as can be seen from the results for $M=10^{-3}$. 
	
	\section{Conclusions}
	\label{sec:Concl}
	We have proposed a first order implicit-explicit numerical method for the simulation of one dimensional isentropic two-phase flow based on the Symmetric Hyperbolic Thermodynamically Compatible model \cite{RomTor2004}.
	The scheme is proved to be consistent with single phase flow, first order accurate and captures accurately material waves in different Mach number regimes. 
	Moreover, the numerical scheme is asymptotic preserving as demonstrated in Theorem \ref{theo:AP1} and \ref{theo:AP2}.
	We have applied a reference solution approach, where stiff non-linear quantities, as pressure and enthalpy, are linearised around a reference state given by the leading order of well-prepared initial data. 
	These data were obtained by an asymptotic analysis of the singular Mach number limits for which the model was reformulated in non-dimensional form. 
	The resulting linear stiff parts were integrated implicitly, whereas the transport terms were treated explicitly leading to the CFL condition that is only restricted by material wave speeds. 
	The final solution contains also explicit nonlinear corrections of pressure and enthalpy yielding an all-speed scheme that performs well also in compressible regimes.
	Due to the complexity of the model, the flux terms were split in three hyperbolic sub-systems.
	This motivates the need to study more complex operator splittings beyond the two subsystems implicit-explicit splittings.
	The question of high order accuracy in the time integration remains open.
		
	A key element of the scheme lies in the reformulation of stiff subsystems in the elliptic form.  
	Even though the system is extremely coupled, we succeeded to solve the homogenous part efficiently with direct or iterative linear solvers.
    Moreover the model contains stiff relaxation source terms describing the interaction of the phases via friction and pressure relaxation processes. 
	In this work, standard methods were used to solve the stiff non-linear pressure relaxation source term implicitly. 
	To improve the implicit treatment of the relaxation source terms, we plan combine the homogeneous all-speed scheme with techniques presented in \cite{ChiMue2020}.
    Therein robust and efficient solvers for relaxation source terms arising in two-phase flows are discussed.
	Further, we aim to extend the scheme to two dimensional problems, as well as to address the full two-phase model given in \cite{RomDriTor2010}.
%	The present AP low Mach number approach can be extended to higher order using 
%	Standard reconstruction techniques keeping the AP property.

	\section*{Acknowledgements}
	A.T. and M.L. have been partially supported by the Gutenberg Research College, JGU Mainz. 
    Further, M.L. is grateful for the support of the Mainz Institute of Multi-Scale Modelling. 
    G.P. is a member of GNCS and acknowledges the support of PRIN2017 and Sapienza, Progetto di Atteneo [RM120172B41DBF3A].
	
	\bibliographystyle{plain}
	\bibliography{lit_2pf.bib}
\end{document}